\documentclass{amsart}
\usepackage{amssymb, amsmath, latexsym, amscd, graphicx, tabularx, color}
\usepackage{rotating}
\usepackage[all]{xy}
\usepackage[dvipdfm]{hyperref}
\usepackage[all]{hypcap}

\newtheorem{theorem}{Theorem}[section]
\newtheorem{lemma}[theorem]{Lemma}
\newtheorem{cor}[theorem]{Corollary}
\newtheorem{definition}[theorem]{Definition}
\newtheorem{proposition}[theorem]{Proposition}
\newtheorem{remark}[theorem]{Remark}
\newtheorem{example}[theorem]{Example}

\def\pagenumber{1}

\begin{document}
\setcounter{page}{\pagenumber}
\newcommand{\T}{\mathbb{T}}
\newcommand{\R}{\mathbb{R}}
\newcommand{\Q}{\mathbb{Q}}
\newcommand{\N}{\mathbb{N}}
\newcommand{\Z}{\mathbb{Z}}
\newcommand{\tx}[1]{\quad\mbox{#1}\quad}
\parindent=0pt
\def\SRA{\hskip 2pt\hbox{$\joinrel\mathrel\circ\joinrel\to$}}
\def\tbox{\hskip 1pt\frame{\vbox{\vbox{\hbox{\boldmath$\scriptstyle\times$}}}}\hskip 2pt}
\def\circvert{\vbox{\hbox to 8.9pt{$\mid$\hskip -3.6pt $\circ$}}}
\def\IM{\hbox{\rm im}\hskip 2pt}
\def\COIM{\hbox{\rm coim}\hskip 2pt}
\def\COKER{\hbox{\rm coker}\hskip 2pt}
\def\TR{\hbox{\rm tr}\hskip 2pt}
\def\GRAD{\hbox{\rm grad}\hskip 2pt}
\def\RANK{\hbox{\rm rank}\hskip 2pt}
\def\MOD{\hbox{\rm mod}\hskip 2pt}
\def\DEN{\hbox{\rm den}\hskip 2pt}
\def\DEG{\hbox{\rm deg}\hskip 2pt}

\title[The Landau's problems.I: The Goldbach's conjecture proved]{THE LANDAU's PROBLEMS.I:\\ THE GOLDBACH's CONJECTURE PROVED}

\author{Agostino Pr\'astaro}
\maketitle
\vspace{-.5cm}
{\footnotesize
\begin{center}
Department SBAI - Mathematics, University of Rome La Sapienza, Via A.Scarpa 16,
00161 Rome, Italy. \\
E-mail: {\tt agostino.prastaro@uniroma1.it}
\end{center}
}
\vskip 0.5cm
\centerline{\em This paper is dedicated to E. Artin, E. Noether, L.S. Pontrjagin and R. Thom.}
\vskip 0.5cm
\begin{abstract}
We give a direct proof of the {\em Goldbach's conjecture}, (GC), in number theory, in the Euler's form. The proof is also constructive, since it gives a criterion to find two prime numbers $\ge 1$, such that their sum gives a fixed even number $\ge 2$. The proof is obtained by recasting the problem in the framework of the Commutative Algebra and Algebraic Topology.
 Even if in this paper we consider $1$ as a prime number, our proof of the GC works also for the {\em restricted Goldbach conjecture}, (RGC), i.e., by excluding $1$ from the set of prime numbers.
\end{abstract}

\vskip 0.5cm

\noindent {\bf AMS Subject Classification:} 11R04; 11T30; 11D99; 11U05; 55N22; 55N20; 81R50; 81T99; 20H15.

\vspace{.08in} \noindent \textbf{Keywords}: Landau's problems; Goldbach's conjecture; Algebraic number theory; Diophantine equations; Quantum algebra; Cobordism groups; Integral cobordism groups of quantum PDEs; Crystallographic groups.

\section[Introduction]{\bf Introduction}
\vskip 0.5cm

\begin{quotation}
\rightline{\footnotesize ``{\em Every even integer is a sum of two primes}."}
\rightline{\footnotesize ``{\em I regard this as a completely certain theorem, although I cannot prove it.}"}
\end{quotation}
\vskip -0.5cm
\rightline{\footnotesize (Euler's letter to Goldbach, June 30, 1742.)}
\vskip 0.5cm

This work in two parts, is devoted to solve the so-called {\em Landau's problems}. These are four well-known problems in Number Theory, listed by Edmund Landau at the 1912 International Congress of Mathematics, remained unsolved up to now. In this first part, the Goldbach's conjecture in Number Theory, is considered. This is the first problem in the Landau's list, and was one of the most famous example of the G\"odel 's incompleteness theorem \cite{GOEDEL1,GOEDEL2,GOEDEL3,HIRZEL}. In this paper we give a direct proof of this conjecture. Some useful applications regarding geometry and quantum algebra are also obtained. (The other three Landau's problems are considered and solved in Part II \cite{PRAS3}.)

Our proof of the  Goldbach's conjecture is motivated by the experimental observation that fixed an even integer, say $2n$, $n\ge 1$, and considered the highest prime number $p_1\in P$, that does not exceed $2n$, the difference $2n-p_1$ is often a prime number, or if not, we can pass to consider the next prime number, say $p_1^{(1)}<p_1$, and find that $2n-p_1^{(1)}$ is just a prime number. (We denote by $P$ the set of prime numbers.) Otherwise, we can continue this process, and after a finite number of steps, obtain that $2n-p_1^{(s)}=p_2^{(s)}$, where $p_2^{(s)}\in P$. This process gives us a practical way to find two primes $p_1^{(s)}$ and $p_2^{(s)}$, such that $2n=p_1^{(s)}+p_2^{(s)}$, hence satisfy the Goldbach's conjecture. In Tab. \ref{criterion-to-find-solution-goldbach-conjecture} are reported some explicit calculations for $2\le 2n\le 998$. Here, in agreement to the original GC, we consider the number $1$ as a prime number. However, our criterion works well also if the number $1$ is excluded by the set of prime numbers. Of course the question is
``{\em Does this phenomenon is a law and why ? }''\footnote{The Goldbach's conjecture formulated in this way is usually called {\em strong GC}. This implies the following {\em weak GC}: ``{\em All odd numbers greater than $7$ are the sum of three odd numbers}.'' Another version of the GC is the following: ``{\em Every integer greater than $5$ can be written as the sum of three primes}.''} The main result of this paper is to prove that this criterion (in the following referred as ``{\em criterion in Tab.~\ref{criterion-to-find-solution-goldbach-conjecture}}''), is mathematically justified.\footnote{Let us emphasize that after this proof the name ``{\em criterion in Tab. \ref{criterion-to-find-solution-goldbach-conjecture}}'' is justified since it allows us to get the goal after a finite number of steps.} For this we recast the problem in the framework of the Commutative Algebra and Algebraic Topology, by showing that to solve the GC is equivalent to understand the algebraic topologic structure of the ring $\mathbb{Z}_{2n}$. In fact, the criterion in Tab. \ref{criterion-to-find-solution-goldbach-conjecture} is encoded by Theorem \ref{goldbach-strong-generators-cyclic-group-2n}. After the proof of this theorem the GC and RGC are simple corollaries.

The paper is organized in an Introduction (Section 1), where we illustrated our criterion to solve the GC, by means of algebraic topologic methods. There is also emphasized by means of a cannot-go theorem (Theorem \ref{a-cannot-go-theorem}) the difficulty to solve the GC by simply looking to the prime numbers in the ring $\mathbb{Z}$ of integers. In  Section 2 we study some fundamental properties of the rings $\mathbb{Z}$ and $\mathbb{Z}_{m}$.\footnote{According to the general mathematical interest of the Goldbach conjecture, this paper has been written in an expository style.} The main result is contained in Theorem \ref{goldbach-strong-generators-cyclic-group-2n} that proves that criterion in Tab. \ref{criterion-to-find-solution-goldbach-conjecture} is justified by some new algebraic topological structures and some properties of commutative algebra applied to suitable rings $\mathbb{Z}_{m}$. Then Corollary \ref{goldbach-conjecture} and Corollary \ref{restricted-goldbach-conjecture} conclude the proof of the GC and RGC too. Corollary \ref{criterion} summarizes above results into a general criterion to find all the Goldbach couples associated to any fixed integer $2n$, $n\ge 1$.\footnote{This criterion agrees with the above one, quoted {\em criterion in Tab. \ref{criterion-to-find-solution-goldbach-conjecture}}, and gives a relation with the structure of suitable ideals of a ring.} In Section 3 are shortly given some applications of the GC respectively in the Euclidean Geometry and in Quantum Algebra and Quantum PDE's, as formulated by A. Pr\'astaro. (For information on this last subject see \cite{PRAS1,PRAS2} and related works quoted therein.) More precisely, in Proposition \ref{goldbach-triangle} we recall a previous application of the GC given by \cite{NAMBIAR} that now is a theorem. This relation is interesting, since it relates the GC to a diophantine equation that, now, after Corollary \ref{goldbach-conjecture} and Corollary \ref{restricted-goldbach-conjecture}, can be considered solved too. Finally Theorem \ref{quantum-algebraic-interpretation-goldbach-conjecture} relates the GC to the quantum algebra and algebraic topology of quantum PDEs, as formulated by A. Pr\'astaro, showing the existence of a canonical homomorphism between the group of even quantum numbers and a suitable group related to a point group of crystallographic groups.

\begin{table}[t]
\caption{Criterion to find a solution to the Goldbach's conjecture: $2n=p_1^{(s)}+p_2^{(s)}$ con $p_1^{(s)},\ p_2^{(s)}\in P$.}
\label{criterion-to-find-solution-goldbach-conjecture}
\scalebox{0.7}{$\begin{tabular}{|c|c|c|l|}
\hline
\hfil{\rm{\footnotesize $n\ge 1$}}\hfil&\hfil{\rm{\footnotesize $2n$}}\hfil&\hfil{\rm{\footnotesize $p_1\in P$}}\hfil&\hfil{\rm{\footnotesize $2n-p_1=p_2\, \Rightarrow\, 2n-p_1^{(1)}=p_2^{(1)}\, \Rightarrow\, \cdots\, 2n-p_1^{(s)}=p_2^{(s)}$}}\hfil\\
\hline\hline
{\rm{\footnotesize $1$}}\hfil&\hfil{\rm{\footnotesize $2$}}&\hfil{\rm{\footnotesize $1$}}&{\rm{\footnotesize $2-1=1$}}\hfill\\
{\rm{\footnotesize $2$}}\hfil&\hfil{\rm{\footnotesize $4$}}&\hfil{\rm{\footnotesize $3$}}&{\rm{\footnotesize $4-3=1$}}\hfill\\
{\rm{\footnotesize $3$}}\hfil&\hfil{\rm{\footnotesize $6$}}&\hfil{\rm{\footnotesize $5$}}&{\rm{\footnotesize $6-5=1$}}\hfill\\
{\rm{\footnotesize $4$}}\hfil&\hfil{\rm{\footnotesize $8$}}&\hfil{\rm{\footnotesize $7$}}&{\rm{\footnotesize $8-7=1$}}\hfill\\
{\rm{\footnotesize $5$}}\hfil&\hfil{\rm{\footnotesize $10$}}&\hfil{\rm{\footnotesize $7$}}&{\rm{\footnotesize $10-7=3$}}\hfill\\
{\rm{\footnotesize $6$}}\hfil&\hfil{\rm{\footnotesize $12$}}&\hfil{\rm{\footnotesize $11$}}&{\rm{\footnotesize $12-11=1$}}\hfill\\
{\rm{\footnotesize $7$}}\hfil&\hfil{\rm{\footnotesize $14$}}&\hfil{\rm{\footnotesize $13$}}&{\rm{\footnotesize $14-13=1$}}\hfill\\
{\rm{\footnotesize $8$}}\hfil&\hfil{\rm{\footnotesize $16$}}&\hfil{\rm{\footnotesize $13$}}&{\rm{\footnotesize $16-13=3$}}\hfill\\
{\rm{\footnotesize $9$}}\hfil&\hfil{\rm{\footnotesize $18$}}&\hfil{\rm{\footnotesize $17$}}&{\rm{\footnotesize $18-17=1$}}\hfill\\
{\rm{\footnotesize $10$}}\hfil&\hfil{\rm{\footnotesize $20$}}&\hfil{\rm{\footnotesize $19$}}&{\rm{\footnotesize $20-19=1$}}\hfill\\
\hline
{\rm{\footnotesize $\cdots$}}\hfil&\hfil{\rm{\footnotesize $\cdots$}}&\hfil{\rm{\footnotesize $\cdots$}}&{\rm{\footnotesize $\cdots$}}\hfill\\
\hline
{\rm{\footnotesize $110$}}\hfil&\hfil{\rm{\footnotesize $220$}}&\hfil{\rm{\footnotesize $211$}}&{\rm{\footnotesize $220-211=9=3\times 3\, \Rightarrow\, 220-199=21=3\times 7$}}\hfill\\
&&&{\rm{\footnotesize $\Rightarrow\, 220-197=23$}}\hfill\\
\hline
{\rm{\footnotesize $173$}}\hfil&\hfil{\rm{\footnotesize $346$}}&\hfil{\rm{\footnotesize $337$}}&{\rm{\footnotesize $346-337=9=3\times 3\, \Rightarrow\, 346-331=15=3\times 5$}}\hfill\\
&&&{\rm{\footnotesize $\Rightarrow\, 346-317=29$}}\hfill\\
\hline
{\rm{\footnotesize $259$}}\hfil&\hfil{\rm{\footnotesize $518$}}&\hfil{\rm{\footnotesize $509$}}&{\rm{\footnotesize $518-509=9=3\times 3\, \Rightarrow\, 518-503=15=3\times 5$}}\hfill\\
&&&{\rm{\footnotesize $\Rightarrow\, 518-499=19$}}\hfill\\
\hline
{\rm{\footnotesize $266$}}\hfil&\hfil{\rm{\footnotesize $532$}}&\hfil{\rm{\footnotesize $523$}}&{\rm{\footnotesize $532-523=9=3\times 3\, \Rightarrow\, 532-521=11$}}\hfill\\
\hline
{\rm{\footnotesize $269$}}\hfil&\hfil{\rm{\footnotesize $538$}}&\hfil{\rm{\footnotesize $523$}}&{\rm{\footnotesize $538-523=15=3\times 5\, \Rightarrow\, 538-521=17$}}\hfill\\
\hline
{\rm{\footnotesize $278$}}\hfil&\hfil{\rm{\footnotesize $556$}}&\hfil{\rm{\footnotesize $547$}}&{\rm{\footnotesize $556-547=9=3\times 3\, \Rightarrow\, 556-541=15=3\times 5$}}\hfill\\
&&&{\rm{\footnotesize $\Rightarrow\, 556-523=33=3\times 11\, \Rightarrow\, 556-521=35=5\times 7$}}\hfill\\
&&&{\rm{\footnotesize $\Rightarrow\, 556-509=47$}}\hfill\\
\hline
{\rm{\footnotesize $298$}}\hfil&\hfil{\rm{\footnotesize $586$}}&\hfil{\rm{\footnotesize $577$}}&{\rm{\footnotesize $586-577=9=3\times 3\, \Rightarrow\, 586-571=15=3\times 5$}}\hfill\\
&&&{\rm{\footnotesize $\Rightarrow\, 586-569=17$}}\hfill\\
\hline
{\rm{\footnotesize $319$}}\hfil&\hfil{\rm{\footnotesize $628$}}&\hfil{\rm{\footnotesize $619$}}&{\rm{\footnotesize $628-619=9=3\times 3\, \Rightarrow\, 628-617=11$}}\hfill\\
\hline
{\rm{\footnotesize $320$}}\hfil&\hfil{\rm{\footnotesize $640$}}&\hfil{\rm{\footnotesize $631$}}&{\rm{\footnotesize $640-631=9=3\times 3\, \Rightarrow\, 640-619=21=3\times 7$}}\hfill\\
&&&{\rm{\footnotesize $\Rightarrow\, 640-617=23$}}\hfill\\
\hline
{\rm{\footnotesize $335$}}\hfil&\hfil{\rm{\footnotesize $670$}}&\hfil{\rm{\footnotesize $661$}}&{\rm{\footnotesize $670-661=9=3\times 3\, \Rightarrow\, 640-659=21=3\times 7$}}\hfill\\
&&&{\rm{\footnotesize $\Rightarrow\, 670-653=17$}}\hfill\\
\hline
{\rm{\footnotesize $350$}}\hfil&\hfil{\rm{\footnotesize $700$}}&\hfil{\rm{\footnotesize $691$}}&{\rm{\footnotesize $700-691=9=3\times 3\, \Rightarrow\, 700-683=17$}}\hfill\\
\hline
{\rm{\footnotesize $309$}}\hfil&\hfil{\rm{\footnotesize $718$}}&\hfil{\rm{\footnotesize $709$}}&{\rm{\footnotesize $718-709=9=3\times 3\, \Rightarrow\, 718-701=17$}}\hfill\\
\hline
{\rm{\footnotesize $391$}}\hfil&\hfil{\rm{\footnotesize $782$}}&\hfil{\rm{\footnotesize $773$}}&{\rm{\footnotesize $782-773=9=3\times 3\, \Rightarrow\, 782-769=13$}}\hfill\\
\hline
{\rm{\footnotesize $393$}}\hfil&\hfil{\rm{\footnotesize $796$}}&\hfil{\rm{\footnotesize $787$}}&{\rm{\footnotesize $796-787=9=3\times 3\, \Rightarrow\, 796-773=23$}}\hfill\\
\hline
{\rm{\footnotesize $403$}}\hfil&\hfil{\rm{\footnotesize $806$}}&\hfil{\rm{\footnotesize $797$}}&{\rm{\footnotesize $806-797=9=3\times 3\, \Rightarrow\, 806-787=19$}}\hfill\\
\hline
{\rm{\footnotesize $410$}}\hfil&\hfil{\rm{\footnotesize $820$}}&\hfil{\rm{\footnotesize $811$}}&{\rm{\footnotesize $820-811=9=3\times 3\, \Rightarrow\, 820-809=11$}}\hfill\\
\hline
{\rm{\footnotesize $419$}}\hfil&\hfil{\rm{\footnotesize $838$}}&\hfil{\rm{\footnotesize $829$}}&{\rm{\footnotesize $838-829=9=3\times 3\, \Rightarrow\, 838-827=11$}}\hfill\\
\hline
{\rm{\footnotesize $424$}}\hfil&\hfil{\rm{\footnotesize $848$}}&\hfil{\rm{\footnotesize $839$}}&{\rm{\footnotesize $848-839=9=3\times 3\, \Rightarrow\, 848-829=19$}}\hfill\\
\hline
{\rm{\footnotesize $436$}}\hfil&\hfil{\rm{\footnotesize $872$}}&\hfil{\rm{\footnotesize $863$}}&{\rm{\footnotesize $872-863=9=3\times 3\, \Rightarrow\, 872-859=13$}}\hfill\\
\hline
{\rm{\footnotesize $448$}}\hfil&\hfil{\rm{\footnotesize $896$}}&\hfil{\rm{\footnotesize $887$}}&{\rm{\footnotesize $896-887=9=3\times 3\, \Rightarrow\, 896-883=13$}}\hfill\\
\hline
{\rm{\footnotesize $451$}}\hfil&\hfil{\rm{\footnotesize $902$}}&\hfil{\rm{\footnotesize $887$}}&{\rm{\footnotesize $902-887=15=3\times 5\, \Rightarrow\, 902-883=19$}}\hfill\\
\hline
{\rm{\footnotesize $464$}}\hfil&\hfil{\rm{\footnotesize $928$}}&\hfil{\rm{\footnotesize $919$}}&{\rm{\footnotesize $928-919=9=3\times 3\, \Rightarrow\, 928-911=17$}}\hfill\\
\hline
{\rm{\footnotesize $481$}}\hfil&\hfil{\rm{\footnotesize $962$}}&\hfil{\rm{\footnotesize $953$}}&{\rm{\footnotesize $962-953=9=3\times 3\, \Rightarrow\, 962-947=15=3\times 5$}}\hfill\\
&&&{\rm{\footnotesize $\Rightarrow\, 962-941=21=3\times 7\, \Rightarrow\, 962-937=25=5\times 5$}}\hfill\\
&&&{\rm{\footnotesize $\Rightarrow\, 962-929=33=3\times 11\, \Rightarrow\, 962-919=43$}}\hfill\\
\hline
{\rm{\footnotesize $486$}}\hfil&\hfil{\rm{\footnotesize $972$}}&\hfil{\rm{\footnotesize $971$}}&{\rm{\footnotesize $972-971=1$}}\hfill\\
\hline
{\rm{\footnotesize $489$}}\hfil&\hfil{\rm{\footnotesize $978$}}&\hfil{\rm{\footnotesize $977$}}&{\rm{\footnotesize $978-977=1$}}\hfill\\
\hline
{\rm{\footnotesize $492$}}\hfil&\hfil{\rm{\footnotesize $984$}}&\hfil{\rm{\footnotesize $983$}}&{\rm{\footnotesize $984-983=1$}}\hfill\\
\hline
{\rm{\footnotesize $496$}}\hfil&\hfil{\rm{\footnotesize $992$}}&\hfil{\rm{\footnotesize $991$}}&{\rm{\footnotesize $992-991=1$}}\hfill\\
\hline
{\rm{\footnotesize $499$}}\hfil&\hfil{\rm{\footnotesize $998$}}&\hfil{\rm{\footnotesize $997$}}&{\rm{\footnotesize $998-997=1$}}\hfill\\
\hline
\multicolumn{4}{l}{\rm{\footnotesize $p_1$ is the highest prime such that $p_1< 2n$.}}\\
\multicolumn{4}{l}{\rm{\footnotesize $p_1^{(i)}$ is the highest prime such that $p_1^{(i)}<p_1^{(i-1)},\, i\ge 1,\ p_1^{(0)}=p_1$.}}\\
\multicolumn{4}{l}{\rm{\footnotesize $p_2^{(s)}$ is the first number in the sequence $i$, $i\ge 1$, such that $p_2^{(s)}\in P$.}}\\
\multicolumn{4}{l}{\rm{\footnotesize $P\subset\mathbb{N}$ is the set of prime numbers of $\mathbb{N}$.}}\\
\end{tabular}$}
\end{table}

Before to pass to the proof of above criterion, i.e., to the proof of the GC, let us emphasize by means of the following theorem the difficulty to prove the GC and RGC by remaining in the Arithmetic framework, namely in the ring $\mathbb{Z}$.

\begin{theorem}[A cannot go theorem]\label{a-cannot-go-theorem}
In general, i.e., for any even integer $2n$, one cannot find two prime integers $p_1$ and $p_2$ satisfying the GC by simply utilizing the primality of these numbers.
\end{theorem}

\begin{proof}
Let us prove that one cannot find two prime integers $p_1,\, p_2\in\mathbb{Z}$, that satisfy the GC simply by using the fact that these numbers must be prime numbers. This can be seen by utilizing the ring structure of $\mathbb{Z}$. In the following lemma we resume some properties of ideals in $\mathbb{Z}$.

\begin{lemma}[Fundamental properties of ideals of $\mathbb{Z}$]\label{fundamental-properties-of-ideals-of-integers-numbers}
One has the following properties for ideals of $\mathbb{Z}$.

{\em 1)} All the ideals of $\mathbb{Z}$ are the principal ideals $n\mathbb{Z}$, $n\ge 0$.\footnote{A principal ideal $\mathfrak{p}$ of a ring $R$, is characterized by the property $x\, y\in \mathfrak{p}\, \Rightarrow x\in \mathfrak{p}$ or $y\in \mathfrak{p}$.} (These are additive subgroups of $\mathbb{Z}$.) One has $n\mathbb{Z}=\mathbb{Z}$ iff $n$ is invertible, i.e., $n=1$.

{\em 2)} $n\mathbb{Z}\subset m\mathbb{Z}$, ($m\ge 1$, $n\ge 1$), iff $n|m$, {\em $m$ divides $n$}, i.e., $n=m\, p$, $p\ge 1$.

{\em 3)} $m\mathbb{Z}$ is a maximal ideal in $\mathbb{Z}$, iff $m$ is prime.

{\em 4)} The principal ideal $m\mathbb{Z}+n\mathbb{Z}=d\mathbb{Z}$, has $d=g.c.d.(m,n)$.

$\bullet$\hskip 2pt Then we can write $d=m\, x+n\, y$, for some $x,\, y\in\mathbb{Z}$.

$\bullet$\hskip 2pt In particular, if $m$ and $n$ are coprimes, then $m\mathbb{Z}+n\mathbb{Z}=\mathbb{Z}$ and $1=m\, x+n\, y$. In such a case $m\mathbb{Z}$ and $n\mathbb{Z}$ are called {\em coprime ideals}.

{\em 5) (Intersection of two ideals)} $m\mathbb{Z}\bigcap n\mathbb{Z}=r\mathbb{Z}$, $r=l.c.m.(m,n)$, $r\ge 1$.

$\bullet$\hskip 2pt Therefore one has $m\mathbb{Z}\bigcap n\mathbb{Z}\not=\varnothing$, and contains $m\, n$.

{\em 6) (Product of two ideals)} $(m\mathbb{Z}) (n\mathbb{Z})=(m\, n)\mathbb{Z}$.

$\bullet$\hskip 2pt Therefore one has $(m\mathbb{Z}) (n\mathbb{Z})=m\mathbb{Z}\bigcap n\mathbb{Z}$ iff $m$ and $n$ are coprimes.

$\bullet$\hskip 2pt $(m\mathbb{Z}+n\mathbb{Z})(m\mathbb{Z}\bigcap n\mathbb{Z})=(m\mathbb{Z})(n\mathbb{Z})$.

{\em 7) (Ideals quotient)} $\mathbb{Z}_n\equiv \mathbb{Z}/n\mathbb{Z}\cong\{0,1,2,\cdots,n-1\}$, $n\ge 1$.

$\bullet$\hskip 2pt If $n$ is prime then $\mathbb{Z}_n$ is the field of the maximal ideal $n\mathbb{Z}\subset\mathbb{Z}$. Then every non-zero element $a\in \mathbb{Z}_n$ is an unit, i.e., $\exists\,  a^{-1}\in \mathbb{Z}_n$, such that $a\, a^{-1}=a^{-1}\, a=1$.

{\em 8)} Let be fixed the positive integers $(n_i)_{1\le i\le n}$. Then one has the canonical ring homomorphism {\em(\ref{ring-homomorphism-product-quotients})}.
\begin{equation}\label{ring-homomorphism-product-quotients}
\left\{
  \phi:\mathbb{Z}\to \prod_{1\le i\le n} \mathbb{Z}_{n_i},\,
  \phi(a)=(a+\mathbb{Z}_{n_i})\right\}.
\end{equation}
$\phi$ is surjective iff $n_i$ and $n_j$ are coprimes for $i\not=j$. $\phi$ is injective iff $\bigcap_{1\le i\le n} n_i\mathbb{Z}=<0>$. This condition is never verified for the ideals of $n_i\mathbb{Z}$, with $n_i\not=0$.

{\em 9)} Let $n=p_1^{r_1}.p_2^{r_2}\cdots p_k^{r_k}$ be the prime factorization of an integer $n\ge 1$. One has the exact commutative diagram reported in {\em(\ref{remainder-theorem-diagram})}.
\begin{equation}\label{remainder-theorem-diagram}
 \xymatrix{&&&0\ar[d]&\\
 0\ar[r]&n\mathbb{Z}\ar@{=}[d]\ar[r]&\mathbb{Z}\ar@{=}[d]\ar[r]^(0.25){j}&{\framebox{$\prod_{1\le i\le k}(\mathbb{Z}_{p_i^{r_i}})$}
 }\ar[d]\ar[r]&0\\
  0\ar[r]&n\mathbb{Z}\ar[r]&\mathbb{Z}\ar[r]&\mathbb{Z}_n\ar[d]\ar[r]&0\\
 &&&0&\\}
\end{equation}
The homomorphism $j$ is given by $j(a)\mapsto(j_i(a))_{1\le i\le k}$, where $j_i:\mathbb{Z}\to \mathbb{Z}_{p_i^{r_i}}$. In other words
$$\phi(a)=(a+p_1^{r_1}\mathbb{Z},\cdots, a+p_k^{r_k}\mathbb{Z}).$$

{\em 10) (Radical of ideal in $\mathbb{Z}$)} The {\em radical} of an ideal $m\mathbb{Z}\subset\mathbb{Z}$ is the ideal
$$\mathfrak{r}(m\mathbb{Z})=\{x\in \mathbb{Z}\, |\, x^n\in\, m\mathbb{Z}\, \hbox{\rm for some}\, n>0\}.$$
Set $\mathfrak{a}=m\mathbb{Z}$. One has the following properties for the radical $\mathfrak{a}$.

{\em(i)} $\mathfrak{r}(\mathfrak{a})\supseteq \mathfrak{a}$. ($\mathfrak{a}$ is called a {\em radical ideal} if   $\mathfrak{r}(\mathfrak{a})= \mathfrak{a}$.)

{\em(ii)} $\mathfrak{r}(\mathfrak{r}(\mathfrak{a}))=\mathfrak{r}(\mathfrak{a})$. Therefore $\mathfrak{r}(\mathfrak{a})$ is a radical ideal.\footnote{For example $\mathfrak{r}(4\mathbb{Z})=2\mathbb{Z}$ and $\mathfrak{r}(2\mathbb{Z})=2\mathbb{Z}$.}

{\em(iii)} $\mathfrak{r}((\mathfrak{a})(\mathfrak{b}))=\mathfrak{r}(\mathfrak{a}\bigcap\mathfrak{b})=
\mathfrak{r}(\mathfrak{a})\bigcap\mathfrak{r}(\mathfrak{b})$.

{\em(iv)} $\mathfrak{r}(\mathfrak{a})=\mathbb{Z}\, \Leftrightarrow\, \mathfrak{a}=<1>$.

{\em(v)} $\mathfrak{r}(\mathfrak{a}+\mathfrak{b})=\mathfrak{r}(\mathfrak{r}(\mathfrak{a})+\mathfrak{r}(\mathfrak{b}))$.

{\em(vi)} If $m$ is prime then $\mathfrak{r}((m\mathbb{Z})^n)=m\mathbb{Z}$, for all $n>0$. ($m\mathbb{Z}$ is an example of radical ideal.)

{\em(vii)} If $m=p_1^{r_1}\cdots p_k^{r_k}$ is the prime factorization of $m$, then $$\mathfrak{r}(m\mathbb{Z})=<p_1,\cdots,p_k>=\bigcap_{1\le i\le k}<p_i>\cong p_1\cdots p_k\mathbb{Z}.$$
Every radical is the intersection of prime ideals containing it.

{\em(viii)} $\mathfrak{r}(m\mathbb{Z})$ and $\mathfrak{r}(n\mathbb{Z})$ are coprime ideals iff $m$ and $n$ are coprime numbers.

\end{lemma}
\begin{proof}
The proof of the propositions of this lemma are standard. (See, e.g., \cite{ATIYAH-MACDONALD, BOURBAKI-1, BOURBAKI-2}.)
\end{proof}

Let us now, take two primes $p_1,\, p_2\in\mathbb{Z}$. From Lemma \ref{fundamental-properties-of-ideals-of-integers-numbers}-4, it follows that
\begin{equation}\label{prime-property}
    p_1\, x+p_2\, y=1
\end{equation}
for some $x,y\in\mathbb{Z}$. Multiplying both sides of equation (\ref{prime-property}) by $2n$, we get
\begin{equation}\label{prime-property-a}
    p_1\, x\, 2n+p_2\, y\, 2n=2n.
\end{equation}
Then from (\ref{prime-property-a}) it should be possible to prove the GC if we should be able to find two prime integers $\bar p_1$ and $\bar p_2$, such that $\bar p_1=p_1\, x\, 2n$ and $\bar p_2=p_2\, y\, 2n$. But this should imply $\bar p_1|p_1$ and $\bar p_2|p_2$.\footnote{Warning. In this paper we adopte the symbol $p|q$ to say that the integer $q$ divides $p$, hence $p=q\cdot r$, for some other integer $r$. This warning is necessary, since in Number Theory one usually adopt the mirror symbol.} This is impossible for prime numbers $\bar p_i$, $i=1,2$. Therefore, the road to find a solution for the GC, simply by starting from two primes, is wrong.
\end{proof}

{\bf Warning.} Let us close this Introduction, by emphasizing that the new algebraic topologic methods used to prove the Goldbach's conjecture are not standard in Number Theory. More precisely let us stress these new methods and their purposes.

{\bf Algebraic Topology -} The first method is focused on new bordism groups ({\em Goldbach-bordism groups}). By means of these mathematical tools it is possible to decide whether, for any fixed positive integer $n$, there exists at least a Goldbach-couple in the interval $(0,2n]$ of integers, i.e., two primes $p$ and $q$ such that $p+q=2n$. Therefore the existence of Goldbach-couples is recast into suitable boundary value problems in Algebraic Topology. (This is the main novelty of our proof. See Lemma \ref{goldbach-bordism-and-goldbach-couples}.)

{\bf Commutative Algebra -} The second method uses commutative algebra and in particular Noether rings, Artin rings, maximal ideals, ... to built all possible Goldbach-couples for any fixed positive integer $n$. These mathematical tools are standard, but the novelty is their use in connection with Goldbach-bordism groups. (In part II \cite{PRAS3} are solved also three other Landau's problems by using the same philosophy, namely by using suitable new bordism groups.)\footnote{This paper is a revised version of the paper posted on arXiv with the same title \cite{PRAS3}. The difference is limited to the Introduction that has been shortly expanded to emphasize the new mathematical methods introduced. This has been made in order to give to the reader an yellow line to follow in this complex paper.}

\section{\bf The Proof}\label{proof-goldbach-conjecture}

In order to build the proof, let us associate to any integer $n\in\mathbb{N}$ the additive group $\mathbb{Z}_n\equiv \mathbb{Z}/n\mathbb{Z}$. Let us consider the following lemmas.

\begin{lemma}\label{cyclic-group-and-zeta-n}
$\bullet$\hskip 2pt Let $G=<a>=\{a=a^1,a^2,\cdots,a^n=e\}$ be a cyclic group of order $n$.\footnote{In a ring $R$, with multiplicative identity element $e$, a {\em root of unity} is any element $a\in R$, of finite multiplicative order, i.e., $a^n=e$. If $\mathbb{F}$ is a {\em Galois field} (i.e., finite field, e.g., $\mathbb{Z}_p$, with $p$ prime) the {\em $n-th$ root of unit} of $\mathbb{F}$, is a solution of the equation $x^n-1=0$ in $\mathbb{F}$.} One has the canonical mapping $G\to \mathbb{Z}_n$, $a^r\mapsto\, [r]$, $1\le r\le n$, that is an isomorphism: $G=<a>\cong \mathbb{Z}_n$.

$\bullet$\hskip 2pt Every group of order $p$ prime is cyclic and abelian.\footnote{A group where every element is of infinite order, is called {\em without torsion}. A {\em group with torsion} is one where every element has finite order. In general every finitely generated abelian group $G$ is a finite direct sum of cyclic subgroups $C_j\cong \mathbb{Z}_{\nu_j}$, $\nu_j\ge 0$. Therefore $G$ has a {\em torsion subgroup} $T\equiv\oplus_{\nu_j>0}C_j=\oplus_{\nu_j>1}C_j$. The {\em free part} of $G$ is $\oplus_{\nu_j=0}C_j\cong G/T$. The number of summand $\mathbb{Z}\cong C_0$ in the free part of $G$ is called the {\em rank} of $G$, and represents the maximal number of linearly independent elements in $G$. The numbers $\nu_j>1$ are called {\em torsion coefficients} of $G$ and can be chosen as powers of prime numbers: $\nu_j=p_j^{\rho_j}$, $p_j\in P$, $\rho_j>0$. Two finitely generated abelian groups are isomorphic iff they have the same rank and the same system of torsion coefficients. (For complementary information see e.g., \cite{BOURBAKI-1, BOURBAKI-2}.)}

$\bullet$\hskip 2pt If $G=<a>$ is a cyclic group of order $n$, the equality $a^\lambda=e$ happens iff $\lambda=q\, n$.

$\bullet$\hskip 2pt Every subgroup of a cyclic group $G=<a>$ is a cyclic group.

$\bullet$\hskip 2pt The subgroup $(a^k)$, $1\le k\le n$, with $a^k\in G$, $G$ cyclic group of order $n$, coincides with $(a^d)$ iff $k=k'\, d$ and $n=n'\, d$. ($d$ divides $k$ and $n$.) Furthermore, the order of $(a^k)$ is $n'=n/d$.

The element $x=a^k$ is a generator of the cyclic group $G=<a>$, of order $n$, iff $k$ and $n$ are coprimes.\footnote{In fact, one has $a^k=a^{k'\, d}\in<a^d>\, \Rightarrow\, <a^k>\subseteq<a^d>$. On the other hand, after the {\em Bezout relation}, $d=\alpha\, n+\beta\, k$, $\alpha,\beta\in\mathbb{Z}$. So we get $a^d=a^{\alpha\, n+\beta\, k}=a^{\beta\, k}\in<a^k>\, \Rightarrow\, (a^d)\subseteq (a^k)$. We can conclude that $(a^d)=(a^k)$.}
\end{lemma}
\begin{lemma}[Euler's totient function and Euler's theorem]\label{cyclic-group-and-number-generators-euler-totient-function}
$\bullet$\hskip 2pt The number of distinct generators of a cyclic group of order $n$ is the {\em Euler's totient function} $\varphi(n)=\sharp\{k\in\mathbb{N}\, |\, g.c.d.(n,k)=1,\, 1\le k< n\}$, i.e., the number of positive prime integers with respect to $n$, in the interval $1\le k< n$.\footnote{For example, the group of units of $\mathbb{Z}_6=\{0,1,2,3,4,5\}$, is $\mathbb{Z}^{\times}_6=\{1,5\}$, hence $\varphi(6)=2$.}

$\bullet$\hskip 2pt {\em(Euler's product formula)} If $n$ admits the prime factorization $n=a_1^{r_1}\cdots a_k^{r_k}$, then one has the relation {\em(\ref{euler-product-formula})} between $\varphi(n)$ and the primes $a_i$, $i=1,\cdots,k$.
\begin{equation}\label{euler-product-formula}
    \varphi(n)=n(1-\frac{1}{a_1})\cdots (1-\frac{1}{a_k})=n\prod_{n|a}(1-\frac{1}{a})
\end{equation}
where the product is over the distinct prime numbers dividing $n$.

$\bullet$\hskip 2pt {\em(Euler's classical formula)} The relation between $n$, its positive divisors $d$ and the Euler's totient function $\varphi$, is given by the formula {\em(\ref{euler-classical-formula})}.
\begin{equation}\label{euler-classical-formula}
\sum_{n|d}\varphi(d)=n.
\end{equation}
where the sum is over the positive divisors $d$ of $n$.

$\bullet$\hskip 2pt {\em(Euler's theorem)} If $a$ is a generator of $\mathbb{Z}_n$, then $a^{\varphi(n)}\, \equiv\, 1\, {\rm mod}\, n$.
\end{lemma}

\begin{lemma}[The ring $\mathbb{Z}_n$ and its authomorphism group]\label{ring-zn-and-its-authomorphism-group} By considering $\mathbb{Z}_n$ a ring, one has the natural ring isomorphism: $$\phi:\mathbb{Z}_n\cong Hom_{Abelian-group}(\mathbb{Z}_n,\mathbb{Z}_n),$$ given by $r\mapsto\phi(r)$, $\phi(r)(p)=p^r=\underbrace{p+\cdots+p}_r$. In particular, if $r$ is coprime with $n$, then $\phi_r:\mathbb{Z}_n\to \mathbb{Z}_n$ is a bijection. Therefore, one has the isomorphism
   $$\mathbb{Z}^{\times}_n\cong Aut_{Abelian-group}(\mathbb{Z}_n),$$
where $\mathbb{Z}^{\times}_n\subset \mathbb{Z}_n$ is the {\em group of units} of the ring $\mathbb{Z}_n$. The elements of $\mathbb{Z}^{\times}_n$ are the generators of $\mathbb{Z}_n$.\footnote{Let us recall that a {\em unit} for an unital commutative ring $R$ is an element a that admits {\em inverse}, i.e., an element $a^{-1}$, such that $aa^{-1}=1$. If $g.c.d.(n,a)=1$ in $\mathbb{Z}$, then, $a$ identifies in $\mathbb{Z}_n$ an unity. In fact, if $a$ is coprime with $n$, then holds the following equation in $\mathbb{Z}$: $x\cdot n+y\cdot a=1$, hence $y\cdot a=1-x\cdot n$, for some $x,\, y\in\mathbb{Z}$. This means that we can write $y\cdot a= 1\, {\rm mod}\, n$, or simply $y\cdot a=1$ in $\mathbb{Z}_n$. Therefore $y=a^{-1}\in \mathbb{Z}^{\times}_n\subset \mathbb{Z}_n$. In Tab. \ref{multiplication-table-in-the-multiplicative-group-integers-mod-10} is reported the multiplication table of $\mathbb{Z}^{\times}_{10}$ and in Tab. \ref{multiplication-table-in-the-multiplicative-group-integers-mod-22} the multiplication table of $\mathbb{Z}^{\times}_{22}$. The group of units of $\mathbb{Z}$ is $\mathbb{Z}^{\times}=\{-1,+1\}$.}
\end{lemma}

\begin{table}[t]
\caption{Multiplication table in $\mathbb{Z}^{\times}_{10}$.}
\label{multiplication-table-in-the-multiplicative-group-integers-mod-10}
\scalebox{0.8}{$\begin{tabular}{|c|c|c|c|c|}
\hline
\hfil{\rm{\footnotesize }}\hfil&\hfil{\rm{\footnotesize $1$}}\hfil&\hfil{\rm{\footnotesize $3$}}\hfil&\hfil{\rm{\footnotesize $7$}}\hfil&\hfil{\rm{\footnotesize $9$}}\hfil\\
\hline
{\rm{\footnotesize $1$}}\hfil&\hfil{\rm{\footnotesize $1$}}\hfil&\hfil{\rm{\footnotesize $3$}}\hfil&\hfil{\rm{\footnotesize $7$}}\hfil&\hfil{\rm{\footnotesize $9$}}\hfil\\
\hline
\hfil{\rm{\footnotesize $3$}}\hfil&\hfil{\rm{\footnotesize $3$}}\hfil&\hfil{\rm{\footnotesize $9$}}\hfil&\hfil{\rm{\footnotesize $1$}}\hfil&\hfil{\rm{\footnotesize $7$}}\hfil\\
\hline
\hfil{\rm{\footnotesize $7$}}\hfil&\hfil{\rm{\footnotesize $7$}}\hfil&\hfil{\rm{\footnotesize $1$}}\hfil&\hfil{\rm{\footnotesize $9$}}\hfil&\hfil{\rm{\footnotesize $3$}}\hfil\\
\hline
\hfil{\rm{\footnotesize $9$}}\hfil&\hfil{\rm{\footnotesize $9$}}\hfil&\hfil{\rm{\footnotesize $7$}}\hfil&\hfil{\rm{\footnotesize $3$}}\hfil&\hfil{\rm{\footnotesize $1$}}\hfil\\
\hline
\multicolumn{5}{l}{\rm{\footnotesize $1^{-1}=1$; $3^{-1}=7$;}}\\
\multicolumn{5}{l}{\rm{\footnotesize $7^{-1}=3$; $9^{-1}=9$.}}\\
\end{tabular}$}
\end{table}

\begin{table}[t]
\caption{Multiplication table in $\mathbb{Z}^{\times}_{22}$.}
\label{multiplication-table-in-the-multiplicative-group-integers-mod-22}
\scalebox{0.8}{$\begin{tabular}{|c|c|c|c|c|c|c|c|c|c|c|}
\hline
\hfil{\rm{\footnotesize }}\hfil&\hfil{\rm{\footnotesize $1$}}\hfil&\hfil{\rm{\footnotesize $3$}}\hfil&\hfil{\rm{\footnotesize $5$}}\hfil&\hfil{\rm{\footnotesize $7$}}\hfil&\hfil{\rm{\footnotesize $9$}}\hfil&\hfil{\rm{\footnotesize $13$}}\hfil&\hfil{\rm{\footnotesize $15$}}\hfil&\hfil{\rm{\footnotesize $17$}}\hfil&\hfil{\rm{\footnotesize $19$}}\hfil&\hfil{\rm{\footnotesize $21$}}\\
\hline
{\rm{\footnotesize $1$}}\hfil&\hfil{\rm{\footnotesize $1$}}&\hfil{\rm{\footnotesize $3$}}&{\rm{\footnotesize $5$}}\hfill&{\rm{\footnotesize $7$}}\hfill&\hfil{\rm{\footnotesize $9$}}\hfil&\hfil{\rm{\footnotesize $13$}}\hfil&\hfil{\rm{\footnotesize $15$}}\hfil&\hfil{\rm{\footnotesize $17$}}\hfil&\hfil{\rm{\footnotesize $19$}}\hfil&\hfil{\rm{\footnotesize $21$}}\hfil\\
\hline
{\rm{\footnotesize $3$}}\hfil&\hfil{\rm{\footnotesize $3$}}&\hfil{\rm{\footnotesize $9$}}&{\rm{\footnotesize $15$}}\hfill&{\rm{\footnotesize $21$}}\hfill&\hfil{\rm{\footnotesize $5$}}\hfil&\hfil{\rm{\footnotesize $17$}}\hfil&\hfil{\rm{\footnotesize $1$}}\hfil&\hfil{\rm{\footnotesize $7$}}\hfil&\hfil{\rm{\footnotesize $13$}}\hfil&\hfil{\rm{\footnotesize $19$}}\hfil\\
\hline
{\rm{\footnotesize $5$}}\hfil&\hfil{\rm{\footnotesize $5$}}&\hfil{\rm{\footnotesize $15$}}&{\rm{\footnotesize $3$}}\hfill&{\rm{\footnotesize $13$}}\hfill&\hfil{\rm{\footnotesize $1$}}\hfil&\hfil{\rm{\footnotesize $21$}}\hfil&\hfil{\rm{\footnotesize $9$}}\hfil&\hfil{\rm{\footnotesize $19$}}\hfil&\hfil{\rm{\footnotesize $7$}}\hfil&\hfil{\rm{\footnotesize $17$}}\hfil\\
\hline
{\rm{\footnotesize $7$}}\hfil&\hfil{\rm{\footnotesize $7$}}&\hfil{\rm{\footnotesize $21$}}&{\rm{\footnotesize $13$}}\hfill&{\rm{\footnotesize $5$}}\hfill&\hfil{\rm{\footnotesize $19$}}\hfil&\hfil{\rm{\footnotesize $3$}}\hfil&\hfil{\rm{\footnotesize $17$}}\hfil&\hfil{\rm{\footnotesize $9$}}\hfil&\hfil{\rm{\footnotesize $1$}}\hfil&\hfil{\rm{\footnotesize $15$}}\hfil\\
\hline
{\rm{\footnotesize $9$}}\hfil&\hfil{\rm{\footnotesize $9$}}&\hfil{\rm{\footnotesize $5$}}&{\rm{\footnotesize $1$}}\hfill&{\rm{\footnotesize $19$}}\hfill&\hfil{\rm{\footnotesize $15$}}\hfil&\hfil{\rm{\footnotesize $7$}}\hfil&\hfil{\rm{\footnotesize $3$}}\hfil&\hfil{\rm{\footnotesize $21$}}\hfil&\hfil{\rm{\footnotesize $17$}}\hfil&\hfil{\rm{\footnotesize $13$}}\hfil\\
\hline
{\rm{\footnotesize $13$}}\hfil&\hfil{\rm{\footnotesize $13$}}&\hfil{\rm{\footnotesize $17$}}&{\rm{\footnotesize $21$}}\hfill&{\rm{\footnotesize $3$}}\hfill&\hfil{\rm{\footnotesize $7$}}\hfil&\hfil{\rm{\footnotesize $15$}}\hfil&\hfil{\rm{\footnotesize $19$}}\hfil&\hfil{\rm{\footnotesize $1$}}\hfil&\hfil{\rm{\footnotesize $5$}}\hfil&\hfil{\rm{\footnotesize $9$}}\hfil\\
\hline
{\rm{\footnotesize $15$}}\hfil&\hfil{\rm{\footnotesize $15$}}&\hfil{\rm{\footnotesize $1$}}&{\rm{\footnotesize $9$}}\hfill&{\rm{\footnotesize $17$}}\hfill&\hfil{\rm{\footnotesize $3$}}\hfil&\hfil{\rm{\footnotesize $19$}}\hfil&\hfil{\rm{\footnotesize $5$}}\hfil&\hfil{\rm{\footnotesize $13$}}\hfil&\hfil{\rm{\footnotesize $21$}}\hfil&\hfil{\rm{\footnotesize $7$}}\hfil\\
\hline
{\rm{\footnotesize $17$}}\hfil&\hfil{\rm{\footnotesize $17$}}&\hfil{\rm{\footnotesize $7$}}&{\rm{\footnotesize $19$}}\hfill&{\rm{\footnotesize $9$}}\hfill&\hfil{\rm{\footnotesize $21$}}\hfil&\hfil{\rm{\footnotesize $1$}}\hfil&\hfil{\rm{\footnotesize $13$}}\hfil&\hfil{\rm{\footnotesize $3$}}\hfil&\hfil{\rm{\footnotesize $15$}}\hfil&\hfil{\rm{\footnotesize $5$}}\hfil\\
\hline
{\rm{\footnotesize $19$}}\hfil&\hfil{\rm{\footnotesize $19$}}&\hfil{\rm{\footnotesize $13$}}&{\rm{\footnotesize $7$}}\hfill&{\rm{\footnotesize $1$}}\hfill&\hfil{\rm{\footnotesize $17$}}\hfil&\hfil{\rm{\footnotesize $5$}}\hfil&\hfil{\rm{\footnotesize $21$}}\hfil&\hfil{\rm{\footnotesize $15$}}\hfil&\hfil{\rm{\footnotesize $9$}}\hfil&\hfil{\rm{\footnotesize $3$}}\hfil\\
\hline
{\rm{\footnotesize $21$}}\hfil&\hfil{\rm{\footnotesize $21$}}&\hfil{\rm{\footnotesize $19$}}&{\rm{\footnotesize $17$}}\hfill&{\rm{\footnotesize $15$}}\hfill&\hfil{\rm{\footnotesize $13$}}\hfil&\hfil{\rm{\footnotesize $9$}}\hfil&\hfil{\rm{\footnotesize $7$}}\hfil&\hfil{\rm{\footnotesize $5$}}\hfil&\hfil{\rm{\footnotesize $3$}}\hfil&\hfil{\rm{\footnotesize $1$}}\hfil\\
\hline
\multicolumn{11}{l}{\rm{\footnotesize $1^{-1}=1$; $3^{-1}=15$; $5^{-1}=9$; $7^{-1}=19$; $9^{-1}=5$;}}\\
\multicolumn{11}{l}{\rm{\footnotesize $13^{-1}=17$; $15^{-1}=3$; $17^{-1}=13$; $19^{-1}=7$;}}\\
\multicolumn{11}{l}{\rm{\footnotesize $21^{-1}=21$.}}\\

\end{tabular}$}
\end{table}

\begin{lemma}\label{cyclic-group-and-subgroups}
Let $H$ be a subgroup of $\mathbb{Z}_n$, of order $b$ and index $c$ in $\mathbb{Z}_n$. Then one has $n=b\, c$ and $H=\frac{c\mathbb{Z}}{n\mathbb{Z}}\cong\, \mathbb{Z}_b\cong\frac{\mathbb{Z}}{b\mathbb{Z}}$. The situation is resumed by the exact commutative diagram {\em(\ref{comm-diagram-sub-group-cyclic-group})}.
\begin{equation}\label{comm-diagram-sub-group-cyclic-group}
   \xymatrix{&&0\ar[d]&0\ar[d]&\\
   0\ar[r]&n{\mathbb{Z}}\ar@{=}[d]\ar[r]&c{\mathbb{Z}}\ar[d]\ar[r]
   &{\framebox{$\frac{c\mathbb{Z}}{n\mathbb{Z}}=\frac{c\mathbb{Z}}{bc\mathbb{Z}}\cong\frac{\mathbb{Z}}{b\mathbb{Z}}=\mathbb{Z}_b$}}\ar[d]\ar[r]&0\\
   0\ar[r]&n{\mathbb{Z}}\ar[r]&{\mathbb{Z}}\ar[d]\ar[r]&{\mathbb{Z}_n}\ar[d]\ar[r]&0\\
   &&{\mathbb{Z}_c}\ar[d]\ar@{=}[r]&{\mathbb{Z}_c}\ar[d]&\\
   &&0&0&\\}
\end{equation}
$\bullet$\hskip 2pt There is a one-to-one correspondence between the ideals $b\mathbb{Z}$ of $\mathbb{Z}$ that contain the ideal $n\mathbb{Z}$ and the ideals of $\mathbb{Z}_n$: $b\mathbb{Z}=\phi^{-1}(\mathbb{Z}_b)$, with $n|b$.

$\bullet$\hskip 2pt For any ideal $n\mathbb{Z}\subset\mathbb{Z}$, $n>1$, there exists a maximal ideal $m\mathbb{Z}\subset\mathbb{Z}$, containing $n\mathbb{Z}$. More precisely, if $n$ admits the following prime factorization $n=p_1^{r_1}\cdots p_k^{r_k}$, then any maximal ideal $p_i\mathbb{Z}$, $i=1,\cdots,k$, contains $n\mathbb{Z}$.

$\bullet$\hskip 2pt Let $r<m$ and $p$ be positive integers, such that $(m-r)|p$, i.e., $m-r=p\, q$, for some positive integer $q\ge 1$. One has the exact commutative diagram {\em(\ref{comm-diagram-sub-group-cyclic-group-a})}.
\begin{equation}\label{comm-diagram-sub-group-cyclic-group-a}
   \xymatrix{&&&0\ar[d]&\\
   0\ar[r]&\mathbb{Z}_p\ar[r]&\mathbb{Z}_{m-r}\ar[r]
   &\mathbb{Z}_{m-r}/\mathbb{Z}_p\ar[d]\ar[r]&0\\
   &&&\mathbb{Z}_q\ar[d]&\\
   &&&0&\\}
\end{equation}
$\bullet$\hskip 2pt Furthermore iff $p$ and $q$ are coprimes then $\mathbb{Z}_{m-r}\cong\mathbb{Z}_p\bigoplus\mathbb{Z}_q$.\footnote{If $n=p_1^{r_1}\cdots p_k^{r_k}$ is the prime factorization of the integer $n$, one has the isomorphism $\mathbb{Z}_n\cong \mathbb{Z}_{p_1^{r_1}} \oplus\cdots\oplus \mathbb{Z}_{p_k^{r_k}}$. (See also \cite{NOETHER}.)}

$\bullet$\hskip 2pt For any couple $(m,p)$ of positive integers, with $p\le m$, one can find another couple $(q,r)$ of positive integers, such that $\mathbb{Z}_p\subset\mathbb{Z}_{m-r}$ and $\mathbb{Z}_{m-r}/\mathbb{Z}_p\cong\mathbb{Z}_q$.

In particular if $m=p$, one has $(q,r)=(1,0)$, hence the following isomorphisms:

$\mathbb{Z}=\mathbb{Z}_p\subset\mathbb{Z}_{m-r}=\mathbb{Z}$ and $\mathbb{Z}_q=0=\mathbb{Z}_{m-r}/\mathbb{Z}_p=\mathbb{Z}/\mathbb{Z}$.
\end{lemma}
\begin{lemma}[Group of units and prime factorization]\label{group-units-and-prime-factorization}
$\bullet$\hskip 2pt If $n$ admits the prime factorization $n=a_1^{r_1}\cdots a_k^{r_k}$, then one has the isomorphism {\em(\ref{isomorphism-group-units-and-prime-factorization})}.
\begin{equation}\label{isomorphism-group-units-and-prime-factorization}
    \mathbb{Z}^{\times}_{n}\cong  \mathbb{Z}^{\times}_{a_1^{r_1}}\times\cdots\times \mathbb{Z}^{\times}_{a_k^{r_k}}.
\end{equation}

$\bullet$\hskip 2pt The multiplicative group $\mathbb{Z}^{\times}_{a^{r}}$ is cyclic for odd primes $a$.

$\bullet$\hskip 2pt The multiplicative group $\mathbb{Z}^{\times}_{n}$ is cyclic iff $\varphi(n)=\lambda(n)$, where $\lambda(n)$ is the {\em Carmichael function} of $n$, i.e., the least common multiple (l.c.m.) of the order of the cyclic groups in the direct product {\em(\ref{isomorphism-group-units-and-prime-factorization})}.\footnote{For example, the group $\mathbb{Z}^{\times}_{10}$ (see Tab. \ref{multiplication-table-in-the-multiplicative-group-integers-mod-10}) is cyclic of order $4$. In fact, one has the splitting $\mathbb{Z}^{\times}_{10}\cong\mathbb{Z}^{\times}_{2}\oplus \mathbb{Z}^{\times}_{5}$. $\varphi(10)=4$, $\varphi(2)=1$, $\varphi(5)=4$ and $\lambda(10)=l.c.m.(1,4)=4$. Thus $\varphi(10)=\lambda(10)$, hence $\mathbb{Z}^{\times}_{10}$ is cyclic. (For any $a\in \mathbb{Z}^{\times}_{10}$, one has $a^4=1$.) Another example is $\mathbb{Z}^{\times}_{22}\cong\mathbb{Z}^{\times}_{2}\oplus \mathbb{Z}^{\times}_{11}$, (see Tab. \ref{multiplication-table-in-the-multiplicative-group-integers-mod-22}), where $\varphi(22)=10=\lambda(22)=l.c.m.(1,10)=10$. So the multiplicative group $\mathbb{Z}^{\times}_{22}$ is cyclic of order $10$. (For any element $a\in \mathbb{Z}^{\times}_{22}$, one has $a^{10}=1$.)}
\end{lemma}
\begin{proof}
It is a consequence of Lemma \ref{fundamental-properties-of-ideals-of-integers-numbers}(8) and of the fact that under multiplication the congruence classes modulo $n$ which are relatively primes to $n$ satisfy the axioms for an abelian group.
\end{proof}

\begin{lemma}[Group of units and primality]\label{group-units-and-primality}
A positive integer $m>1$ is prime iff $\varphi(m)=m-1$.
\end{lemma}

\begin{proof}
This follows from Euler's product formula (Lemma \ref{euler-product-formula}) and according to Lemma \ref{group-units-and-prime-factorization}. In fact, $m$ is prime iff $m=a\in P$, hence $|\mathbb{Z}^{\times}_{m}|=|\mathbb{Z}_m|-1=m-1$.
\end{proof}
\begin{lemma}[Maximal ideals in $\mathbb{Z}$]\label{maximal-ideals-in-the-ring-integers}
In the set $\Sigma$ of all ideals, $\not=<1>$, of $\mathbb{Z}$ any chain has at least a maximal ideal.

In particular, any chain in $\Sigma$, that ends with a prime ideal $d\mathbb{Z}$, i.e., $d$ is a prime number, has this ideal as maximal ideal.
\end{lemma}
\begin{proof}
Let order $\Sigma$ by inclusion. Let apply Zorn's lemma to $\Sigma$, i.e., let us show that every chain in $\Sigma$ has an upper bound in $\Sigma$. In fact, let $(n_\alpha\mathbb{Z})_{1\le \alpha\le r}$ be a chain of ideals such that $n_i\mathbb{Z}\subseteq n_{i+1}\mathbb{Z}$. Set $\mathfrak{a}=\bigcup_{1\le\alpha\le r}n_\alpha\mathbb{Z}$. Then $\mathfrak{a}$ is an ideal and $1\not\in\mathfrak{a}$ because $1\not\in\, n_\alpha\mathbb{Z}$ for all $\alpha$. Hence $\mathfrak{a}\in\Sigma$, and $\mathfrak{a}$ is an upper bound of the chain. From Zorn's lemma $\Sigma$ must have at least a maximal element. In fact, from Lemma \ref{cyclic-group-and-subgroups} it follows that in order to be satisfied the condition on the chain, must be $n_i|n_{i+1}$. On the other hand, since all ideals in $\mathbb{Z}$ are principal, must there exist a positive integer $d$, such that $\mathfrak{a}=d\mathbb{Z}$. Really if $n_r=p_1^{s_1}\cdots p_k^{s_k}$ is the prime factorization of $n_r$, we can see that $\mathfrak{a}$ can coincide with any of the following maximal ideals $p_i\mathbb{Z}$, $i=1,\cdots,k$. In particular if $k=1$, i.e., $n_r$ is a prime number, there exists only one maximal ideal of the chain.
\end{proof}

\begin{lemma}[Maximal ideals in $\mathbb{Z}_n$]\label{maximal-ideals-in-the-ring-integers-modulo-n}
The maximal ideals in $\mathbb{Z}_n$ are $p_i\mathbb{Z}/n\mathbb{Z}$, $i=1,\cdots,k$, if $n=p_1^{r_1}\cdots p_k^{r_k}$ is the prime factorization of $n$.
\end{lemma}
\begin{proof}
The proof follows directly from Lemma \ref{cyclic-group-and-subgroups} and Lemma \ref{maximal-ideals-in-the-ring-integers}.
\end{proof}

\begin{lemma}[Jacobson radical of the ring $\mathbb{Z}$]\label{jacobson-radical-of-the-integers-ring}
The {\em Jacobson radical} $J(\mathbb{Z})$ of $\mathbb{Z}$, is for definition, the intersection of the maximal ideals of $\mathbb{Z}$, hence $J(\mathbb{Z})=\{0\}$.

The Jacobson radical of the ring $\mathbb{Z}_n$ is $J(\mathbb{Z}_n)=p_1\cdots p_k\mathbb{Z}/n\mathbb{Z}$, if $n=p_1^{r_1}\cdots p_k^{r_k}$ is the prime factorization of $n$. ($J(\mathbb{Z}_n)$ coincides with the {\em nilradical} of $\mathbb{Z}_n$.)\footnote{For example $J(\mathbb{Z}_{15})=\{0\}$. Instead $J(\mathbb{Z}_{12})=6\mathbb{Z}/12\mathbb{Z}$.}

$\bullet$\hskip 2pt $\mathbb{Z}_n/J(\mathbb{Z}_n)$, is a semiprimitive ring.\footnote{A semiprimitive ring $R$ is one where $J(R)=\{0\}$. It is always semiprimitive the quotient ring $R/J(R)$, i.e., $J(R/J(R))=\{0\}$.}

$\bullet$\hskip 2pt $\mathbb{Z}_n$, with $n$ prime is a semiprimitive ring, (since it is a field).
\end{lemma}

\begin{lemma}[Local rings and semi-local rings]\label{local-rings-and-semi-local-rings}

$\bullet$\hskip 2pt $\mathbb{Z}$ is not a local ring and neither a semi-local ring.

$\bullet$\hskip 2pt $\mathbb{Z}_n$ is a semi-local ring. If $n$ is prime $\mathbb{Z}_n$ becomes a local ring with $\{0\}=J(\mathbb{Z}_n)$ the unique maximal ideal. Therefore $J(\mathbb{Z}_n)=\mathbb{Z}_n\setminus\mathbb{Z}_n^{\times}=\{0\}$, since $\mathbb{Z}_n$ is a field, hence semiprimitive.
\end{lemma}
\begin{proof}
These are direct consequences of the following definitions and results in commutative algebra. A {\em local ring} is a ring with exactly one maximal ideal. A {\em semi-local ring} is a ring with a finite number of maximal ideals. In a local ring $R$, $J(R)=R\setminus R^{\times}$, i.e., the Jacobson radical coincides with the non-units of $R$.

In a local ring $R$, $R/J(R)\cong R/\mathfrak{m}$ is a field, hence semiprimitive. Here $\mathfrak{m}$ is the unique maximal ideal of $R$.
\end{proof}
\begin{lemma}[Nilradical]\label{nil-radical}
The {\em nilradical} $\mathfrak{n}(\mathbb{Z})$ of $\mathbb{Z}$ coincides with $J(\mathbb{Z})$: $\mathfrak{n}(\mathbb{Z})=J(\mathbb{Z})=\{0\}$.

The same happens for the ring $\mathbb{Z}_n$: $\mathfrak{n}(\mathbb{Z}_n)=J(\mathbb{Z}_n)$.
\end{lemma}
\begin{proof}
Let us recall that the nilradical of a ring $R$ is the ideal $\mathfrak{n}(R)$ of its elements $x\in R$, such that $x^n=0$, for some integer $n>0$. $\mathfrak{n}(R)$ is obtained by intersection of all prime ideals of $R$. $\mathfrak{n}(R)$ can be considered the radical of the zero-ideal: $\mathfrak{n}(R)=\mathfrak{r}(<0>)$. In general any maximal ideal is prime, but the converse is not true. In fact the ring $\mathbb{Z}$, has as prime ideals $<m>$, with $m=0$ or $m$ a prime number $\not=1$. The maximal ideal are only the ones with $m$ prime, $\not=1$. However the intersection of all maximal ideals coincides with the ones of all prime ideals and it is just $\{0\}$. Similar considerations hold for the ring $\mathbb{Z}_n$.
\end{proof}

\begin{lemma}[Non-units and maximal ideals]\label{non-units-and-maximal-ideals}
$\bullet$\hskip 2pt Every non-unit of $\mathbb{Z}$ is contained into a maximal ideal.

$\bullet$\hskip 2pt Every non-unit of $\mathbb{Z}_n$ is contained into a maximal ideal.
\end{lemma}
\begin{proof}
The proof can be considered as an application of a similar statement for rings. However, let us see a direct proof. Let us start with the ring $\mathbb{Z}$. Let $n\in\mathbb{Z}\setminus \{-1,1\}$. Since $n\in n\mathbb{Z}\subseteq p\mathbb{Z}$, where $p$ is any prime such that $n|p$. Therefore $n$ belongs to the maximal ideal $p\mathbb{Z}$.

Let us consider the case of the ring $\mathbb{Z}_n$. Then if $a\in \mathbb{Z}_n\setminus\mathbb{Z}_n^{\times}$, it follows that we can write $a$, considered as belonging to $\mathbb{Z}_n$, as $a+ n\mathbb{Z}$. Since $a$ necessarily divides $n$, we can write $a=p.q$, for some prime $p$, such that it appears in the prime factorization of $n$. Therefore we can write $a=p.q+p.q'\mathbb{Z}$, where $n=p.q'$. As a by product we get $a=p(q+q'\mathbb{Z})$. On the other hand $(p\mathbb{Z}/n\mathbb{Z})\cong (\mathbb{Z}/q'\mathbb{Z})=\mathbb{Z}_{q'}$, and since $p\mathbb{Z}/n\mathbb{Z}$ is a maximal ideal in $\mathbb{Z}_n$, it follows that belongs to a maximal ideal in $\mathbb{Z}_n$. As a consequence one has also that $a=p(q+q'\mathbb{Z})$ belongs to the same maximal ideal $\mathbb{Z}_{q'}$ in $\mathbb{Z}_n$, since $p.q\in\mathbb{Z}_{q'}$.
\end{proof}

\begin{lemma}[The rings $\mathbb{Z}$ and $\mathbb{Z}_n$ as  $\mathbb{Z}$-modules]\label{rings-of-integers-and-modulo-n-as-z-molules}
$\bullet$\hskip 2pt The ring $\mathbb{Z}$ has a canonical structure of finitely generated free $\mathbb{Z}$-modulo by means of the following short exact sequence:
$$\xymatrix{0\ar[r]&\mathbb{Z}\ar[r]^{\phi=1}&\mathbb{Z}\ar[r]&0\\}$$

$\bullet$\hskip 2pt The ring $\mathbb{Z}_n$ has a natural structure of finitely generated $\mathbb{Z}$-module by means of the following short exact sequence:
$$\xymatrix{0\ar[r]&H_n\ar[r]&\mathbb{Z}^{\varphi(n)}\ar[r]^{\phi}&\mathbb{Z}_n\ar[r]&0\\}$$
where $\phi$ is defined by $$x=\phi(x^1,\cdots,x^{\varphi(n)})=\sum_{1\le k\le \varphi(n)} x^k\, a_k,\, x^k\in\mathbb{Z}$$
and $\{a_i\}_{1\le i\le\varphi(n)}$ is a set of generators of $\mathbb{Z}_n$. Furthermore $H_n=\ker(\phi)$, is defined by the linear equation (in $\mathbb{Z}_{n}$): $ \sum_{1\le k\le \varphi(n)} x^k\, a_k=0$. One has the isomorphims: $\mathbb{Z}_n\cong \mathbb{Z}^{\varphi(n)}/H_n$.
\end{lemma}

\begin{lemma}[$\mathbb{Z}$ as a Noetherian ring]\label{ring-of-integers-as-noetherian-ring}
$\bullet$\hskip 2pt $\mathbb{Z}$ is a Noetherian ring, i.e., any ascending chain of ideals in $\mathbb{Z}$, terminates (or stabilizes) after a finite number of steps. The maximal ideal of the chain is a prime ideal, i.e., an ideal of the type $p\mathbb{Z}$ with $p$ a positive prime number.

$\bullet$\hskip 2pt $\mathbb{Z}$ is not an Artinian ring.

$\bullet$\hskip 2pt {\em(Dimension)}, $dim(\mathbb{Z})=1$, where $dim(\mathbb{Z})$ is the supremum of the lengths of chains of prime ideals, in $\mathbb{Z}$.

$\bullet$\hskip 2pt The {\em prime spectrum} $Spec(\mathbb{Z})$ of $\mathbb{Z}$ is a topological space (with the {\em Zariski topology}).
\end{lemma}
\begin{proof}
In $\mathbb{Z}$ ideals are principal ideals, of the type $m\mathbb{Z}=<m>$, where $m$ are positive numbers. Moreover, $m\mathbb{Z}\subseteq n\mathbb{Z}$, iff $m|n$. Therefore, a chain  $m\mathbb{Z}\subseteq n\mathbb{Z}\subseteq p\mathbb{Z}\subseteq\cdots$, must necessarily terminates after a finite number of steps, since the possible positive numbers that divide $m$ cannot exceed $m$. Furthermore, taking into account the prime factorization of $m$ it is clear that the maximal ideal in the chain is a prime ideal.

In $\mathbb{Z}$ any descending chain of ideals is of the type
$$m\mathbb{Z}\supseteq p_1m\mathbb{Z}\supseteq p_2p_1m\mathbb{Z}\supseteq \cdots$$
where $m$, $p_i$ are positive numbers $>1$. Such chains cannot stabilize after a finite number of steps, since we can always find ideals $k\mathbb{Z}$, with $k$ a multiple of the previous one in the chain. The intersection of all such ideals is the trivial ideal $<0>$.

The strictly increasing chains of prime ideals in $\mathbb{Z}$ are of the type $\mathfrak{p}_0=<0>\subset \mathfrak{p}_1=p\mathbb{Z}$, or $\mathfrak{p}_0=p\mathbb{Z}$, with $p>1$ prime. Therefore, the supremum of the lengths of such chains is $1$. This is also the dimension of $\mathbb{Z}$.

The set $Spec(\mathbb{Z})$ of all prime ideals in $\mathbb{Z}$ is a topological space with Zariski topology, i.e., generated by closed subsets, defined by $V(X)$, for any subset $X\subset \mathbb{Z}$, as the set of all prime ideals of $\mathbb{Z}$ that contain $X$. $V(X)$ satisfy the following properties.

(i) If $\mathfrak{a}=<X>\subset\mathbb{Z}$, is the ideal generated by $X$, then $V(X)=V(\mathfrak{a})=V(\mathfrak{r}(\mathfrak{a}))$. (If $a\in \mathbb{Z}$, then $\mathfrak{a}=a\mathbb{Z}$.)

(ii) $V(0)=Spec(\mathbb{Z})$.

(iii)  $V(1)=\varnothing$.

(iv) If $(X_i)_{i\in I}$ is any family of subsets of $\mathbb{Z}$, then $V(\bigcup_{i\in I}X_i)=\bigcap_{i\in I}V(X_i)$.

(v) $V(m\mathbb{Z}\bigcap n\mathbb{Z})=V(mn\mathbb{Z})=V(m\mathbb{Z})\bigcup V(n\mathbb{Z})$, for any ideal $m\mathbb{Z}$ and $n\mathbb{Z}$ of $\mathbb{Z}$.

(vi) $V(\sum_{i}\mathfrak{a}_i)=\bigcap_{i}V(\mathfrak{a}_i)$.
The {\em basic open sets} of $Spec(\mathbb{Z})$ is made by sets $X_a=Spec(\mathbb{Z})\setminus V(a)$, for any $a\in Spec(\mathbb{Z})$. The sets $X_a$ are open sets in the Zariski topology of $Spec(\mathbb{Z})$, and satisfy to the following properties.

(vii) $X_a\bigcap X_b=X_{ab}$.

(viii) $X_a=\varnothing\, \Leftrightarrow\, a$ is nilpotent.

(ix) $X_a=Spec(\mathbb{Z})\, \Leftrightarrow\, a$ is a unit.

(x) $X_a=X_b\, \Leftrightarrow\, \mathfrak{r}(<a>)=\mathfrak{r}(<b>)$.

(xi) $Spec(\mathbb{Z})$ is quasi-compact (that is, every open covering of $Spec(\mathbb{Z})$ has a finite subcovering).\footnote{''Quasi-compact'' means ''compact but not necessarily Hausdorff''.}

(xii) Each $X_a$ is quasi-compact.

(xiii) An open subset of $Spec(\mathbb{Z})$ is quasi-compact iff it is a finite union of sets $X_a$.

(xiv) Let $<x>\in Spec(\mathbb{Z})$, be a point of the prime spectrum of $\mathbb{Z}$, i.e., $x$ prime. Then $<x>\equiv x\mathbb{Z}$ is closed in the Zariski topology of $\mathbb{Z}$ iff $x\mathbb{Z}$ is maximal. On the other hand all prime ideals in $\mathbb{Z}$ are maximal ones, hence any point $<x>$ is closed in $Spec(\mathbb{Z})$. Therefore, $Spec(\mathbb{Z})$ is a {\em $T_0$-space}, i.e., if $<x>$ and $<y>$ are distinct points of $Spec(\mathbb{Z})$, then either there is a neighborhood of $<x>$ which does not contain $<y>$, or else there is a neighborhood of $<y>$ which does not contain $<x>$.

(xv) $Spec(\mathbb{Z})$ is an {\em irreducible space}, i.e., any pair of non-empty open sets in the Zariski topology, intersect, or equivalently every non-empty open set is dense in $Spec(\mathbb{Z})$. This is equivalent to say that $\mathfrak{n}(\mathbb{Z})=<0>$.

(xvi) $Spec(\mathbb{Z})=\{\mathfrak{p}\, |\, \mathfrak{p}\subset \mathbb{Z}\, {\rm prime ideal}\}\bigcup\{<0>\}$. Every prime ideal is closed in $Spec(\mathbb{Z})$, except $<0>$, whose closure is $V(0)=Spec(\mathbb{Z})$.
\end{proof}

\begin{lemma}[$\mathbb{Z}_n$ as a Noetherian and Artinian ring]\label{ring-of-integers-modulo-n-as-noetherian-artinian-ring}
$\bullet$\hskip 2pt $\mathbb{Z}_n$ is a Noetherian and Artinian ring.
\end{lemma}
\begin{proof}
Since $\mathbb{Z}_n$ is a finitely generated commutative ring, it is a Noetherian ring.\footnote{Another, way to prove that $\mathbb{Z}_n$ is a Noetherian ring, is to use the following theorem: If $R$ is a Noetherian ring, and $\mathfrak{a}$ is an ideal of $R$, then $R/\mathfrak{a}$ is a Noetherian ring too.\cite{ATIYAH-MACDONALD} In fact, it is enough to take $R=\mathbb{Z}$ and $\mathfrak{a}=n\mathbb{Z}$. This agrees with the epimorphism $\pi:\mathbb{Z}\to \mathbb{Z}_n$, since $\mathbb{Z}_n\cong \mathbb{Z}/\ker(\pi)$.} More precisely, any ascending chain of ideals in $\mathbb{Z}_n$ is of the type:
$$p\mathbb{Z}/n\mathbb{Z}\cong \mathbb{Z}_q\subseteq r\mathbb{Z}/q\mathbb{Z}\cong\mathbb{Z}_s\subseteq\cdots$$
where $n=pq$, $q=rs$, etc. This chain necessarily stops after a finite number of step since the numbers $q$, $r$, etc. all divide $n$, hence the steps in the chain cannot be more than $n$. Furthermore, the last ideal in the chain must be corresponding to a prime number, that results a maximal ideal.

To prove that $\mathbb{Z}_n$ is Artinian, it is enough to prove that $\dim(\mathbb{Z}_n)=0$. In fact, any Noetherian ring is an Artinian ring iff its dimension is zero. \cite{ATIYAH-MACDONALD} On the other hand all the prime ideals of $\mathbb{Z}_n$ are of the type $\mathbb{Z}_p$, where $p$ is a prime number such that $n|p$. Therefore, any strictly increasing chain of prime ideals in $\mathbb{Z}_n$ can be made by only one ideal: $\mathfrak{p}_0=\mathbb{Z}_p$ wit $p$ a prime number, $n|p$, hence the dimension of the ring $\mathbb{Z}_n$ must necessarily be $0$. Therefore, $\mathbb{Z}_n$ is an Artinian ring.

This means that any descending chain of ideals in $\mathbb{Z}_n$
$$\mathfrak{p}_0\supseteq \mathfrak{p}_1\supseteq \mathfrak{p}_2\supseteq \cdots$$
stops (or stabilizes) after a finite number of steps. Now, after above considerations it results that any ascending chain of ideals in $\mathbb{Z}_n$ is of the type
$$\mathbb{Z}_a\supseteq \mathbb{Z}_b\supseteq\mathbb{Z}_c\supseteq\cdots$$
with $a=pb$, $b=rc$, etc. Therefore, since $a$ must be a multiple of any of the numbers $b$, $c$ etc., it follows that such a chain must stop after a finite number of steps, since $a$ is a fixed number. More precisely, the chain stabilizes at an ideal $\mathbb{Z}_x$, where $x$ is a prime number entering in the prime factorization of $a$.
\end{proof}
\begin{remark}
Let us emphasize that after Lemma \ref{cyclic-group-and-number-generators-euler-totient-function} one can understand that the numbers $p_1^{(s)}$ and $p_2^{(s)}$ considered in our criterion to find a solution to the Goldbach's conjecture, are just generators of $\mathbb{Z}_{2n}$. However, they are, in a sense, distinguished generators since they are not only prime with respect to $2n$, but are just prime numbers.
\end{remark}

\begin{definition}[Strong generators in $\mathbb{Z}_m$]\label{strong-generators-in-cyclic-group}
We call {\em strong generators} in $\mathbb{Z}_m$ the generators that are identified by prime numbers. Let us denote by $\mathbb{Z}_m^{\blacksquare}$ the set of strong generators of $\mathbb{Z}_m^{\blacksquare}$. One has the natural inclusions: $$\mathbb{Z}_m^{\blacksquare}\subset\mathbb{Z}_m^{\times}\subset\mathbb{Z}_m.$$
\end{definition}

\begin{proposition}[Existence of strong generators in a cyclic group]\label{existence-strong-generators-in-cyclic-group}
In $\mathbb{Z}_{2n}$, $n\ge 1$, there exist strong generators. When $n>1$, $\mathbb{Z}_{2n}^{\blacksquare}\supset \{1\}$.
\end{proposition}
\begin{proof}
In fact in the set of generators of $\mathbb{Z}_{2n}$ there exists always $1$, for any positive number $n\ge 1$. However, when $n>1$, $\mathbb{Z}_{2n}^{\blacksquare}$ properly contains $1$. Let us denote respectively by $p_k$ the primes entering in the factorization of $2n$, $a_i$ the units that are not primes and by $b_j$ the units that are primes. The prime factorization of $a_k$ must be of the type $a_k=b_1^{m_1}\cdots b_h^{m_k}$, since $a_k$ are coprimes with $2n$. Then any $c\in\mathbb{Z}_{2n}$ can be written in the form $c=x^ka_k+y^jb_j$, $x^k,\, y^j\in\mathbb{Z}$. If we assume that with $n>1$, $\mathbb{Z}_{2n}^{\blacksquare}=\{1\}$, then also the units $a_k$ should reduce to $1$, and any $c\in \mathbb{Z}_{2n}$, should be written $c=x\cdot 1$. This can be happen iff $\mathbb{Z}_{2n}=\mathbb{Z}_2$, hence $n=1$, in contrast with the assumption that $n>1$. This just means that for $n>1$, $\mathbb{Z}_{2n}^{\blacksquare}$ is larger than $\{1\}$.
\end{proof}

\begin{example}
In Tab. \ref{example-strong-generators-in-cyclic-groups} we report generators and strong generators, with respect to examples just considered in Tab. \ref{criterion-to-find-solution-goldbach-conjecture}. There we can verify that some couples of generators satisfy equation $2n=a+b$, but these do not necessitate to be strong generators in $\mathbb{Z}_{2n}$.
\end{example}

\begin{table}[t]
\caption{Examples of strong generators in $\mathbb{Z}_{2n}$.}
\label{example-strong-generators-in-cyclic-groups}
\scalebox{0.75}{$\begin{tabular}{|c|l|l|c|c|c|}
\hline
\hfil{\rm{\footnotesize $\mathbb{Z}_{2n}$}}\hfil&\hfil{\rm{\footnotesize Goldbach's couples $(\clubsuit)$}}\hfil&\hfil{\rm{\footnotesize $\mathbb{Z}^{\times}_{2n}$ (Generators group or group of units)}}\hfil&\hfil{\rm{\footnotesize $\varphi(2n)$}}\hfil&\hfil{\rm{\footnotesize Quasi-Goldbach's couples $(\spadesuit)$}}\hfil\\
\hline\hline
{\rm{\footnotesize $\mathbb{Z}_{2}$}}\hfil&\hfil{\rm{\footnotesize $(1,1)^\bigstar$}}&\hfil{\rm{\footnotesize $\{1\}$}}&{\rm{\footnotesize $1$}}\hfill&{\rm{\footnotesize $-$}}\hfill\\
{\rm{\footnotesize $\mathbb{Z}_{4}$}}\hfil&\hfil{\rm{\footnotesize $(1,3)^\bigstar$}}&\hfil{\rm{\footnotesize $\{1,3\}$}}&{\rm{\footnotesize $2$}}\hfill&{\rm{\footnotesize $-$}}\hfill\\
{\rm{\footnotesize $\mathbb{Z}_{6}$}}\hfil&\hfil{\rm{\footnotesize $(1,5)^\bigstar$}}&\hfil{\rm{\footnotesize $\{1,5\}$}}&{\rm{\footnotesize $2$}}\hfill&{\rm{\footnotesize $-$}}\hfill\\
{\rm{\footnotesize $\mathbb{Z}_{8}$}}\hfil&\hfil{\rm{\footnotesize $(1,7)^\bigstar;\, (3,5)$}}&\hfil{\rm{\footnotesize $\{1,3,5,7\}$}}&{\rm{\footnotesize $4$}}\hfill&{\rm{\footnotesize $-$}}\hfill\\
{\rm{\footnotesize $\mathbb{Z}_{10}$}}\hfil&\hfil{\rm{\footnotesize $(3,7)^\bigstar$}}&\hfil{\rm{\footnotesize $\{1,3,7,(9)\}$}}&{\rm{\footnotesize $4$}}\hfill&{\rm{\footnotesize $(1,9)$}}\hfill\\
{\rm{\footnotesize $\mathbb{Z}_{12}$}}\hfil&\hfil{\rm{\footnotesize $(1,11)^\bigstar;\, (5,7)$}}&\hfil{\rm{\footnotesize $\{1,5,7,11\}$}}&{\rm{\footnotesize $4$}}\hfill&{\rm{\footnotesize $-$}}\hfill\\
{\rm{\footnotesize $\mathbb{Z}_{14}$}}\hfil&\hfil{\rm{\footnotesize $(1,13)^\bigstar;\, (3,11)$}}&\hfil{\rm{\footnotesize $\{1,3,5,(9),11,13\}$}}&{\rm{\footnotesize $6$}}\hfill&{\rm{\footnotesize $(5,9)$}}\hfill\\
{\rm{\footnotesize $\mathbb{Z}_{16}$}}\hfil&\hfil{\rm{\footnotesize $(3,13)^\bigstar;\, (5,11)$}}&\hfil{\rm{\footnotesize $\{1,3,5,7,(9),11,13,(15)\}$}}&{\rm{\footnotesize $8$}}\hfill&{\rm{\footnotesize $(1,15);\, (7,9)$}}\hfill\\
{\rm{\footnotesize $\mathbb{Z}_{18}$}}\hfil&\hfil{\rm{\footnotesize $(1,17)^\bigstar;\, (5,13);\, (7,11)$}}&\hfil{\rm{\footnotesize $\{1,5,7,11,13,17\}$}}&{\rm{\footnotesize $6$}}\hfill&{\rm{\footnotesize $-$}}\hfill\\
{\rm{\footnotesize $\mathbb{Z}_{20}$}}\hfil&\hfil{\rm{\footnotesize $(1,19)^\bigstar;\, (3,17);\, (7,13)$}}&\hfil{\rm{\footnotesize $\{1,3,7,(9),11,13,17,19\}$}}&{\rm{\footnotesize $8$}}\hfill&{\rm{\footnotesize $(9,11)$}}\hfill\\
\hline
{\rm{\footnotesize $\mathbb{Z}_{22}$}}\hfil&\hfil{\rm{\footnotesize $(3,19)^\bigstar;\, (5,17)$}}&\hfil{\rm{\footnotesize $\{1,3,5,7,(9),11,13,(15),17,19,21\}$}}&{\rm{\footnotesize $10$}}\hfill&{\rm{\footnotesize $(1,21);\, (7,15);\, (9,13)$}}\hfill\\
\hline
{\rm{\footnotesize $\mathbb{Z}_{28}$}}\hfil&\hfil{\rm{\footnotesize $(5,23)^\bigstar;\, (11,17)$}}&\hfil{\rm{\footnotesize $\{1,3,5,(9),11,13,(15),17,19,23,(25),(27)\}$}}&{\rm{\footnotesize $12$}}\hfill&{\rm{\footnotesize $(1,27);\, (3,25);\, (9,19);\, (13,15)$}}\hfill\\
\hline
\multicolumn{5}{l}{\rm{\footnotesize The Goldbach's couples marked by  $()^\bigstar$ are ones obtained by criterion in Tab. \ref{criterion-to-find-solution-goldbach-conjecture}.}}\\
\multicolumn{5}{l}{\rm{\footnotesize The set of strong generators is obtained by the ones of generators, by forgetting the numbers between brackets $()$ in $\mathbb{Z}^{\times}_{2n}$.}}\\
\multicolumn{5}{l}{\rm{\footnotesize $\mathbb{Z}^{\times}_{2n}=\{k\in \mathbb{Z}_{2n}\, |\, g.c.d.(2n,k)=1,\, 1\le k<2n\}$ is also called the {\em multiplicative group of integers} (mod $2n$).}}\\
\multicolumn{5}{l}{\rm{\footnotesize $(\clubsuit)$ Warn ! In this table, except for the case $n=1$, do no appear {\em trivial Goldbach couples}, ($(n,n)$ with $n$ prime).}}\\
\multicolumn{5}{l}{\rm{\footnotesize $(\clubsuit)$ Except in the case $n=1$, trivial Goldbach couples are never identified by units in $\mathbb{Z}_{2n}$. (See Lemma \ref{existence-trivial-oldbach-couples}.)}}\\
\multicolumn{5}{l}{\rm{\footnotesize $(\spadesuit)$ In this table are reported the quasi-Goldbach couples that are not Goldbach couples.}}\\
\end{tabular}$}
\end{table}

So, in order to prove GC, we are conduced to prove Theorem \ref{goldbach-strong-generators-cyclic-group-2n}.
\begin{theorem}[Goldbach's couples in $\mathbb{Z}_{2n}$]\label{goldbach-strong-generators-cyclic-group-2n}

$\bullet$\hskip 2pt In the group $\mathbb{Z}_{2n}$ there exist two strong generators identified by positive primes $a$ and $b$ that satisfy the condition  {\em(\ref{equation-goldbach-strong-generators-cyclic-group-2n})}.\footnote{In this paper we denote by $P$ the subset of $\mathbb{N}$ given by all prime natural numbers. It is well known that $P$ is infinite, ({\em Euclide's theorem}), and therefore $P$ has the same cardinality of $\mathbb{N}$: $\sharp(P)=\sharp(\mathbb{N})=\aleph_0$. Recall, also, that $\sharp(\mathbb{Z})=\aleph_0$, with bijection $f:\mathbb{N}\to\mathbb{Z}$, given by $\{f(1)=0,\, f(2n)=n,\, f(2n+1)=-n\}_{n\ge 1}$.}
\begin{equation}\label{equation-goldbach-strong-generators-cyclic-group-2n}
   2n=a+b,\, a,\, b\in P.
\end{equation}

$\bullet$\hskip 2pt We call {\em Goldbach's couples} in $\mathbb{Z}_{2n}$, couples of strong generators of $\mathbb{Z}_{2n}$, identified by two positive primes $a$ and $b$ that satisfy the condition {\em(\ref{equation-goldbach-strong-generators-cyclic-group-2n})}.\footnote{In the following we shall often use the same symbol to denote a number $a\in\mathbb{Z}$ and its projection $\pi(a)\in\mathbb{Z}_m$, via the canonical projection $\pi:\mathbb{Z}\to\mathbb{Z}_m$. In fact, from the context it will be clear what is the right interpretation !}

$\bullet$\hskip 2pt We call also {\em quasi-Goldbach's couples} in $\mathbb{Z}_{2n}$, couples of generators $(a, b)$ of $\mathbb{Z}_{2n}$, that satisfy the condition $2n=a+b$, but where one of the numbers $a$ or $b$ does not necessitate to be prime. (All Goldbach couples are also quasi-Goldbach couples.)

$\bullet$\hskip 2pt Goldbach's couples do not necessitate to be unique in  $\mathbb{Z}_{2n}$, for any $n> 3$.

$\bullet$\hskip 2pt We call {\em canonical Goldbach couple} of $2n$, the first obtained by applying the criterion in Tab. \ref{criterion-to-find-solution-goldbach-conjecture}.

$\bullet$\hskip 2pt We call {\em Noether-Goldbach's couple} in $\mathbb{Z}_{2n}$, the quasi-Goldbach couple $(1,2n-1)$, when it is also a Goldbach couple. If there exists the Noether-Goldbach couple, this is the canonical one too.
\end{theorem}
\begin{proof}
Let us consider the following lemmas.
\begin{lemma}\label{strong-generators-and-criterion-1}
The strong generators of $\mathbb{Z}_{2n}$ satisfy the following properties.

{\em(i)} Each strong generator of $\mathbb{Z}_{2n}$, generates all $\mathbb{Z}_{2n}$.

{\em(ii)} If $p_1\in\mathbb{Z}_m^{\blacksquare}$ then $2n-p_1=p_2$ is a generator of $\mathbb{Z}_{2n}$, i.e., $p_2\in\mathbb{Z}_{2n}^{\times}$. Then $p_2$ has the prime factorization {\em(\ref{prime-factorization-quasi-goldbach-couple})}.
\begin{equation}\label{prime-factorization-quasi-goldbach-couple}
    p_2=b_1^{u_1}\cdots b_s^{u_s}
\end{equation}
where $b_i$ identify strong generators in $\mathbb{Z}_{2n}$. Therefore $p_2$ is coprime with $p_1$ iff $b_i\not=p_1$, $i=1,\cdots, s$.
\end{lemma}
\begin{proof}
The first proposition follows from the fact that a strong generator is a unit of $\mathbb{Z}_{2n}$.

The second proposition follows from the prime factorization of $2n=a_1^{r_1}\cdots a_k^{r_k}$. In fact these primes numbers cannot coincide with $p_1$, since this last is a unit, hence $g.c.d.(2n,p_1)=1$. Therefore the number $2n-p_1=p_2$ cannot be factorized as $a_s^{m}\, q$, with $a_s$ coinciding with a prime number $a_i$, appearing in the prime factorization of $2n$. In other words $g.c.d.(2n,p_2)=1$, hence $p_2\in \mathbb{Z}_{2n}^{\times}$. Furthermore, if $p_2\in \mathbb{Z}_{2n}^{\times}$ then in its prime factorization $p_2=b_1^{u_1}\cdots b_s^{u_s}$ cannot appear the prime numbers of the prime factorization of $2n=a_1^{r_1}\cdots a_k^{r_k}$. This proves the factorization (\ref{prime-factorization-quasi-goldbach-couple}), hence the condition in order $p_2$ should be coprime with $p_1$.
\end{proof}

\begin{lemma}\label{strong-generators-and-criterion-2}
Let $p_2\in\mathbb{Z}_{2n}^{\times}\subset \mathbb{Z}_{2n}$, as defined in Lemma \ref{strong-generators-and-criterion-1}. $p_2$ is prime iff it identifies a strong generator in $\mathbb{Z}_{2n}$, i.e., $p_2$ (or more precisely its projection in $\mathbb{Z}_{2n}$) belongs to  $\mathbb{Z}_{2n}^{\blacksquare}\subset \mathbb{Z}_{2n}^{\times}\subset\mathbb{Z}_{2n}$.
\end{lemma}
\begin{proof}
This follows directly from prime factorization (\ref{prime-factorization-quasi-goldbach-couple}).
\end{proof}

\begin{lemma}[Existence of trivial Goldbach couples]\label{existence-trivial-oldbach-couples}
If $2n$ admits the prime factorization $2n=a_1^{r_1}\cdots a_k^{r_k}$, then one has $2n-a_i\in P[1,2n]$ iff $2n=2a_i$. Here $P[1,2n]$ denotes the set of primes in the interval $[1,2n]$. Thus, except in the {\em trivial cases}, i.e., where $n$ is a prime $a_i$, to the primes $a_i\in P[1,2n]$, entering in the prime factorization of $2n$, cannot be associated Goldbach couples of $2n$. Therefore, in non-trivial cases, a necessary condition for the Goldbach couples $(p_1,p_2)$ of $2n$ is that $p_2=2n-b$, with $b$ identifying in $\mathbb{Z}_{2n}$ strong generators, namely $b\in\mathbb{Z}_{2n}^{\blacksquare}$.
\end{lemma}
\begin{proof}
Let us note that in $P[1,2n]$ admits the following partition in two disjoint sets: $P[1,2n]=P[1,2n]^{\square}\sqcup P[1,2n]^{\blacksquare}$, where $P[1,2n]^{\square}$ denotes the primes entering in the factorization of $2n$ and $P[1,2n]^{\blacksquare}$ are the other primes that identify strong generators in $\mathbb{Z}_{2n}$. If we denotes by $\mathbf{P}[1,2n]$ the projection of $P[1,2n]$ into $\mathbb{Z}_{2n}$, by means of the canonical epimorphism $\pi:\mathbb{Z}\to\mathbb{Z}_{2n}$, then we induce the following partition in $\mathbf{P}[1,2n]$:
$\mathbf{P}[1,2n]=\mathbb{Z}_{2n}^{\square}\sqcup \mathbb{Z}_{2n}^{\blacksquare}$, where $\mathbb{Z}_{2n}^{\square}=\pi(P[1,2n]^{\square})\subset\mathbb{Z}_{2n}$. In order to see under which conditions $2n-a_i$ is prime, let us represent this number with respect the prime factorization of $2n$.
\begin{equation}\label{representation-with-respect-prime-factorization-2n-2n-ai}
   \left\{
   \begin{array}{ll}
     2n-a_i&=a_1^{r_1}\cdots a_k^{r_k}-a_i \\
     & =a_i(a_1^{r_1}\cdots a_i^{r_i-1}\cdots a_k^{r_k}-1)
   \end{array}
   \right.
\end{equation}
Therefore, $2n-a_i$ is prime iff $a_1^{r_1}\cdots a_i^{r_i-1}\cdots a_k^{r_k}-1=1$, namely  $a_1^{r_1}\cdots a_i^{r_i-1}\cdots a_k^{r_k}=2$. This condition can be verified iff $2n=2a_i$.
\end{proof}

\begin{lemma}[Maximal ideals and $2n-1$]\label{maximal-ideals-and-2n-1}
$\bullet$\hskip 2pt If $2n$ admits the prime factorization $2n=a_1^{r_1}\cdots a_k^{r_k}$, $a_i\in P[1,2n]$, then $2n-1$ cannot belong to the ideal  $\mathbb{Z}_{a_i}\subset\mathbb{Z}_{2n}$.

$\bullet$\hskip 2pt The same holds for any element of $\mathbb{Z}_{2n}^{\times}$.

$\bullet$\hskip 2pt Any element of $\mathbb{Z}_{2n}^{\times}$ cannot belong to the Jacobson radical $J(\mathbb{Z}_{2n})\cong a_1\cdots a_k\mathbb{Z}/2n\mathbb{Z}$, namely the intersection of all maximal ideals of $\mathbb{Z}_{2n}$.
\end{lemma}
\begin{proof}
In fact $2n-1$ is coprime with $2n$, hence identifies an element of $\mathbb{Z}_{2n}^{\times}$. Therefore, it cannot be contained into a maximal ideal of $\mathbb{Z}_{2n}$. These are of the type $\mathbb{Z}_{a_i}$, where $a_i$ is a prime entering in the prime factorization of $2n$.

The other propositions are direct consequences of above properties.
\end{proof}

\begin{lemma}[Mirror symmetry in $\mathbb{Z}_{2n}^{\times}$]\label{mirror-symmetry-in-units-group}
The integers $a_i$ in the interval $[1,2n]$, that identify units in $\mathbb{Z}_{2n}$, are symmetrically distributed around the middle. Therefore, in $\mathbb{Z}_{2n}^{\times}$ the order is always even: $\varphi(2n)=2d$.\footnote{See Tab. \ref{example-strong-generators-in-cyclic-groups} for some examples.}
\end{lemma}
\begin{proof}
In fact, from any of such $a_i$ we can see that $2n-a_i=a_j$ where $a_j$ identifies another unit in $\mathbb{Z}_{2n}$. The proof is similar to the one considered for the Lemma \ref{strong-generators-and-criterion-1}.
\end{proof}
\begin{lemma}[Mirror symmetry in Goldbach couples and existence of Goldbach couples and Noether-Goldbach couples]\label{mirror-symmetry-and-noether-goldbach-couples}
$\bullet$\hskip 2pt For any fixed even integer $2n$, $n\ge 1$, the Goldbach couples are symmetrically distributed around the middle in the interval $[1,2n]$.

$\bullet$\hskip 2pt Goldbach couples are identified by $b_i\in\mathbb{Z}^{\blacksquare}_{2n}$, $b_i\ge n$ iff there exists a strong generator $b_j$, symmetric to $b_i$ with respect to the middle, or equivalently $b_i-n=n-b_j$.

$\bullet$\hskip 2pt In the case that $b_i=n$, then there exists the trivial Goldbach couple $(n,n)$.

$\bullet$\hskip 2pt The {\em Noether-Goldbach couple} of $2n$, $n>1$, exists iff the order of $\mathbb{Z}^{\times}_{2n-1}$ is $2(n-1)$.

\end{lemma}
\begin{proof}
In fact, for any Goldbach couple $(p_1,p_2)$ we can write the condition $2n=p_1+p_2$ in the forma $p_1-n=n-p_2$.

The second proposition follows directly from the previous one.

If $2n-1$ is a prime, then the quasi-Goldbach couple $(2n-1,1)$ becomes a (canonical) Noether-Golbach couple. On the other and a positive integer $m$ is prime iff the order of $\mathbb{Z}_m^{\times}$ is $m-1$, i.e. $\varphi(m)=m-1$. (Lemma \ref{group-units-and-primality}.) Therefore $2n-1$ is prime iff the order of $\mathbb{Z}_{2n-1}^{\times}$ is $2n-1-1=2(n-1)$.
\end{proof}

Even if there is a mirror symmetry in the distribution of the Goldbach-couples, this does not origin from an analogous symmetry in the set of strong generators. In fact, we get the following lemma.

\begin{lemma}[No-mirror symmetry in $\mathbb{Z}^{\blacksquare}_{2n}$]\label{no-mirror-symmetry-in-strong-generators}
The strong generators do not respect the mirror symmetry, in the sense that if there exists a strong generator $b_i\ge n$ of $2n$, does not necessitate that there is also a strong generator $b_j\le n$, such that $b_i-n=n-b_j$.
\end{lemma}
\begin{proof}
This can be proved with a counterexample. For example in $\mathbb{Z}_{10}$, the mirror symmetric of $1$ does not exist. This should be $9$, but it is not prime. Another example could be $556$, where the strong generator $547$ has not a mirror symmetric strong generator. (See Tab. \ref{criterion-to-find-solution-goldbach-conjecture}.) In fact, the absence of mirror symmetry in $\mathbb{Z}^{\blacksquare}_{2n}$ produces quasi-Goldbach couples that are not Goldbach couples.
\end{proof}

\begin{definition}[Noether numbers]\label{noether-numbers}
We call {\em Noether numbers} the even numbers $2n$ such in $\mathbb{Z}_{2n}$ there exists a (canonical) Noether-Goldbach couple.
\end{definition}

\begin{lemma}[Existence of Noether-numbers]\label{existence-noether-numbers}
$2n$ is a Noether number iff the order of $\mathbb{Z}_{2n-1}^{\times}$ is $2(n-1)$.
\end{lemma}
\begin{proof}
This is a by-product of Lemma \ref{mirror-symmetry-and-noether-goldbach-couples} and Definition \ref{noether-numbers}
\end{proof}
\begin{lemma}[Goldbach couples, splitting of the ring $\mathbb{Z}$ and algebraic relations in $\mathbb{Z}$ and $\mathbb{Z}_{2n}$]\label{goldbach-couples-splitting-ring-z}
$\bullet$\hskip 2pt Any non-trivial Goldbach couple $(b_i,b_j)$, $i\not= j$, $b_i,\ b_j\in\mathbb{Z}^{\blacksquare}_{2n}$, gives split representation of the ring $\mathbb{Z}$: $\mathbb{Z}=b_i\mathbb{Z}+b_j\mathbb{Z}$.

This means that hold the equations {\em(\ref{algebraic-relations-in-goldbach-couple})} relating the elements in a same Goldbach couple, but also different elements of different Goldbach couples, and with $2n$.
\begin{equation}\label{algebraic-relations-in-goldbach-couple}
    \left\{
    \begin{array}{l}
      b_i\cdot x+b_j\cdot y=1,\, i\not=j,\ x,y\in\mathbb{Z}\\
      2n\cdot \bar x+b\cdot \bar y=1,\, i\not=j,\ \bar x,\bar y\in\mathbb{Z}\\
    \end{array}
    \right\}\, b_i,\, b_j,\, b\in P[1,2n]^{\blacksquare}
\end{equation}
Above equations {\em(\ref{algebraic-relations-in-goldbach-couple})} can be reinterpreted as equations in $\mathbb{Z}_{2n}$.

\end{lemma}
\begin{proof}
In fact, it is enough to apply Lemma \ref{fundamental-properties-of-ideals-of-integers-numbers}, taking into account that $2n$ is coprime with any $b\in P[1,2n]^{\blacksquare}$.
Furthermore equations (\ref{algebraic-relations-in-goldbach-couple}) can be reinterpreted in $\mathbb{Z}_{2n}$, taking into account the isomorphisms (\ref{algebraic-relations-in-goldbach-couple-a}).
\begin{equation}\label{algebraic-relations-in-goldbach-couple-a}
    \left\{
    \begin{array}{l}
      (b_i\mathbb{Z}+b_j\mathbb{Z})/2n\mathbb{Z}\cong \mathbb{Z}/2n\mathbb{Z}=\mathbb{Z}_{2n}\\
      (2n\mathbb{Z}+b\mathbb{Z})/2n\mathbb{Z}\cong \mathbb{Z}/2n\mathbb{Z}=\mathbb{Z}_{2n}\\
    \end{array}
    \right.
\end{equation}
\end{proof}

\begin{example}
Let us consider the case $2n=22$. See Tab. \ref{example-strong-generators-in-cyclic-groups} for corresponding characterizations of Goldbach couples. Then $p_2=2n-1=21$ is not a prime number, in other word $(1,21)$ is a quasi Goldbach couple. (In fact the canonical Goldbach couple is $(3,19)$.) However, the equation $22-21=1$ says that $22$ is coprime with $21$, hence this equation written in $\mathbb{Z}$, can be rewritten also in $\mathbb{Z}_{22}$, where we can write $19\cdot 19^{-1}=1$. (See Tab. \ref{multiplication-table-in-the-multiplicative-group-integers-mod-22}.) In this way we get the following equation in $\mathbb{Z}_{22}$.  $22\cdot 21-19\cdot 15=1$. This can be rewritten in $\mathbb{Z}$, since $21$ and $19$ are coprimes. We get $-9\cdot 21+10\cdot 19=1$. Instead if we made a similar calculation with $3\cdot 3^{-1}=1$ in $\mathbb{Z}_{22}$ we arrive to the following equation $x\, 21+y\, 3=1$, with $x=22$ and $y=-21\cdot 3^{-1}$. This equation cannot be rewritten in $\mathbb{Z}$, since $21$ is not coprime with $3$.
\end{example}

\begin{lemma}[Strong generators and ring isomorphisms]\label{strong-generators-and-ring-isomorphisms}
Let $\{b_j\}_{1\le j\le s}$ be the strong generators of $\mathbb{Z}_{2n}$. Then one has the ring isomorphism
\begin{equation}\label{strong-generators-and-ring-isomorphisms-a}
    \mathbb{Z}_{b_1\cdots b_s}\cong\prod_{j}\mathbb{Z}_{_{b_j}}.
\end{equation}
\begin{proof}
In fact one has the short exact sequence (\ref{strong-generators-and-ring-isomorphisms-b}).
\begin{equation}\label{strong-generators-and-ring-isomorphisms-b}
    \xymatrix{0\ar[r]&\framebox{$\ker(\phi)=\bigcap_{j}b_j\mathbb{Z}$}\ar[r]&\mathbb{Z}\ar[r]^(0.4){\phi}&
    \prod_{j}(\mathbb{Z}/b_j\mathbb{Z})\ar[r]&0\\}
\end{equation}
The morphism $\phi$ is surjective since the ideals $b_j\mathbb{Z}\subset \mathbb{Z}$ are primes. Therefore one has the isomorphism $$\mathbb{Z}/\ker(\phi)\cong \mathbb{Z}/b_1\cdots b_s\mathbb{Z}=\mathbb{Z}_{b_1\cdots b_s}\cong\prod_{j}(\mathbb{Z}_{b_j}).$$
\end{proof}
\end{lemma}

From above lemmas, and taking into account the criterion in Tab. \ref{criterion-to-find-solution-goldbach-conjecture}, it is clear that since the set $\mathbb{Z}_{2n}^{\blacksquare}$ is finite, and contains prime numbers, even if these do not respect the mirror-symmetry with respect to the middle of the interval $[1,2n]$, (see Proposition \ref{existence-strong-generators-in-cyclic-group}, Lemma \ref{mirror-symmetry-and-noether-goldbach-couples} and Lemma \ref{no-mirror-symmetry-in-strong-generators}), it follows that $p_2+a=2n-(p_1-a)$ must necessarily coincide with a prime number after some finite steps. In fact, in each of this step $p_1-a\ge n$ is taken a strong generator. More precisely, if $\mathbb{Z}^{\times}_{2n-1}$ is of order $2(n-1)$, then $2n$ is a Noether number, hence there is the canonical Noether-Goldbach couple $(1,2n-1)$ of $2n$. Moreover if $n$ is prime, there exists the trivial Goldbach couple $(n,n)$. Other Goldbach couples, when occur,  can be found by considering the $2n-b_j$, with $1<b_j< 2n-1$, strong generators in $\mathbb{Z}_{2n}^{\blacksquare}$. In order to be more explicit in our proof, let us associate to any number $1<2n-b_i=a_j<2n-1$ in our process, the ideal $\mathfrak{a}_i=(2n-b_i)\mathbb{Z}/r\mathbb{Z}\subset \mathbb{Z}_r$, where $r=l.c.m.(a_1,\cdots,a_k)$. Here we denote by $a_i$ the integers in the open interval $]1,2n-1[$ that identify the units of $\mathbb{Z}_{2n}$, and by $b_j$, the $a_i$ that are primes, hence their projections under $\pi:\mathbb{Z}\to\mathbb{Z}_{2n}$, identify the strong generators of $\mathbb{Z}_{2n}$. Then the set $\{\mathfrak{a}_i\}$ of ideals in $\mathbb{Z}_r$, associated to the criterion in Tab. \ref{criterion-to-find-solution-goldbach-conjecture}, must have maximal elements, since $\mathbb{Z}_r$ is Noetherian. (See Lemma \ref{ring-of-integers-modulo-n-as-noetherian-artinian-ring}.) Warn ! We are not interested to a maximal ideal of the set $\{\mathfrak{a}_i\}$, but to ideals in $\{\mathfrak{a}_i\}$ that are maximal ideals of $\mathbb{Z}_r$ ! These exist just for the Noetherian structure of the ring $\mathbb{Z}_r$, and are in a finite number since $\mathbb{Z}_r$ is an Artinian ring. (See Lemma \ref{ring-of-integers-modulo-n-as-noetherian-artinian-ring}.)  On the other hand, any maximal ideal $\mathfrak{m}$ in $\mathbb{Z}_r$ is of the type $\mathfrak{m}=b\mathbb{Z}/r\mathbb{Z}$, with $b\not=1$ a prime of the interval $[1,2n]$, identifying a strong generator of $\mathbb{Z}_{2n}^{\blacksquare}$. By looking to maximal ideals of $\mathbb{Z}_r$, in the set $\{\mathfrak{a}_i\}$, we are sure that these are of the type $b\mathbb{Z}/r\mathbb{Z}$, with $b$ some prime $b=2n-b_i\not=1$.\footnote{Warn ! In general $\mathbb{Z}_r$ contains a finite number of maximal ideals, since it is a semi-local ring. (See Lemma \ref{local-rings-and-semi-local-rings}.) Thus, we can identify by means of such maximal ideals all the possible Goldbach couples, when they occur in $\{\mathfrak{a}_i\}$. However, two different maximal ideals can identify the same Goldbach couple for effect of the mirror symmetry. (See Lemma \ref{mirror-symmetry-and-noether-goldbach-couples}.)}

\begin{lemma}[Relation between $\mathbb{Z}_r$ and maximal ideals]\label{relation-between-zr-and-maximal-ideals}
Let us denote $\mathfrak{m}_j=\frac{b_j\mathbb{Z}}{r\mathbb{Z}}$, $1\le j\le s$ be the maximal ideals in $\mathbb{Z}_r$. One has the short exact sequence {\rm(\ref{short-exact-sequence-for-zr})}.
\begin{equation}\label{short-exact-sequence-for-zr}
    \xymatrix{0\ar[r]&\framebox{$\ker(\phi)=\bigcap_{j}\mathfrak{m}_j$}\ar[r]&\mathbb{Z}_r\ar[r]^(0.4){\phi}&
    \prod_{j}(\mathbb{Z}/\mathfrak{m}_j)\ar[r]&0\\}
\end{equation}
and therefore one has the following isomorphisms: $\mathbb{Z}_r/\mathfrak{m}_j\cong\mathbb{Z}_{b_j}$, and
$$\mathbb{Z}_r/\ker(\phi)\cong\mathbb{Z}_{b_1\cdots b_s}\cong\prod_j\mathbb{Z}/\mathfrak{m}_j\cong\prod_j\mathbb{Z}_{b_j}.$$
\end{lemma}
\begin{proof}
The proof is similar to the one of Lemma \ref{strong-generators-and-ring-isomorphisms}.
\end{proof}
\begin{lemma}
Any two maximal ideals $\mathfrak{m}_1$, $\mathfrak{m}_2$ in $\mathbb{Z}_r$ are coprimes, i.e., $\mathfrak{m}_1+ \mathfrak{m}_2=\mathbb{Z}_r$.\footnote{Warn ! Do not confuse the sum with the direct sum. See Lemma \ref{relations-between-direct-sum-intersection-and-sum-of-z-modules}.}
\end{lemma}
\begin{proof}
In fact, $\mathfrak{m}_1+ \mathfrak{m}_2=\frac{b_1\mathbb{Z}}{r\mathbb{Z}}+\frac{b_2\mathbb{Z}}{r\mathbb{Z}}=
\frac{b_1\mathbb{Z}+b_2\mathbb{Z}}{r\mathbb{Z}}
=\frac{g.c.d.(b_1,b_2)\mathbb{Z}}{r\mathbb{Z}}=\frac{\mathbb{Z}}{r\mathbb{Z}}=\mathbb{Z}_r.$
\end{proof}

\begin{lemma}[Rapresentation of $\mathbb{Z}_r$ by means of local rings]\label{rapresentation-of-zr-by-means-of-local-rings}
$\mathbb{Z}_r$ is isomorphic to the direct product of a finite number of local artin rings. Furthermore,
one has the canonical isomorphism {\em(\ref{isomorphirsm-zr-local-rings})}.
\begin{equation}\label{isomorphirsm-zr-local-rings}
    \mathbb{Z}_r\cong\prod_{x\in Spec(\mathbb{Z}_r}(\mathbb{Z}_r)_x
\end{equation}
where $(\mathbb{Z}_r)_x$ is $\mathbb{Z}_r$ localized at $x$.
\end{lemma}
\begin{proof}
In fact, $\mathbb{Z}_r$ is an Artinian ring. (See also Lemma \ref{cyclic-group-and-subgroups} and Lemma \ref{ring-of-integers-modulo-n-as-noetherian-artinian-ring}.)  Furthermore, $\mathbb{Z}_r$ as an Artinian ring, has a finite number of maximal ideals, and in $\mathbb{Z}_r$ all prime ideals are maximal ideals too. (This agrees with the fact that in Artinian rings all the prime ideals are maximal ones.) Thus $Spec(\mathbb{Z}_r)=Max(\mathbb{Z}_r)$. i.e., the prime spectrum coincides with the maximal spectrum. Taking into account that $\mathbb{Z}_r$ is also a Noetherian ring, one has that $Spec(\mathbb{Z}_r)$ is a finite Hausdorff reducible Noetherian topological space consisting of a finite number of points. These points are closed and open in the Zariski topology, i.e., $Spec(\mathbb{Z}_r)$ is a discrete topological space. One has the short exact sequence (\ref{short-exact-sequence-local-rings}).\footnote{$\mathcal{O}_{Spec(\mathbb{Z}_r)}$ is the sheaf over $Spec(\mathbb{Z}_r)$, identified by $\mathbb{Z}_r$, since $\mathcal{O}_{Spec(\mathbb{Z}_r)}(Spec(\mathbb{Z}_r))=\mathbb{Z}_r$. If $x=\mathfrak{b}$ is a point of $Spec(\mathbb{Z}_r)$, then $\mathop{\lim}\limits_{\overrightarrow{U\ni x}}\mathcal{O}_{Spec(\mathbb{Z}_r)}(U)\cong(\mathbb{Z}_r)_{\mathfrak{b}}$, where the limit is made by means of the restriction homomorphism. (See, e.g., \cite{ATIYAH-MACDONALD}.)}

\begin{equation}\label{short-exact-sequence-local-rings}
    \scalebox{0.8}{$\xymatrix{0\ar[r]&\framebox{$\mathbb{Z}_r=\Gamma(Spec(\mathbb{Z}_r),\mathcal{O}_{Specc(\mathbb{Z}_r)})$}\ar[r]^(0.4){\phi}&\framebox{$\prod_{x\in Spec(\mathbb{Z}_r)}(\mathbb{Z}_r)_x=\prod_{x\in Spec(\mathbb{Z}_r)}\mathcal{O}_{Spec(\mathbb{Z}_r),x}$}\ar[r]&0\\}$}
\end{equation}
In fact $\phi$ is naturally injective. Furthermore, since each point $x$ is also open, then $(\mathbb{Z}_r)_x=\Gamma(\{x\},\mathcal{O}_{Spec(\mathbb{Z}_r)})$, and $\{x\}\bigcap\{y\}=\varnothing$ if $x\not= y$. As a by product, it follows that a section $s\in\Gamma(Spec(\mathbb{Z}_r),\mathcal{O}_{Spec(\mathbb{Z}_r)})$ can be built by a collection of sections $s(x)\in\Gamma(\{x\},\mathcal{O}_{Spec(\mathbb{Z}_r)})$, for $x\in Spec(\mathbb{Z}_r)$. Therefore $\phi$ is surjective too.
\end{proof}

\begin{lemma}[Spectral properties of $\mathbb{Z}_r$ and existence of Goldbach couples]\label{spectral-properties-of-zr-and-existence-of-goldbach-couples}
The points $\mathfrak{b}_j=\frac{b_j\mathbb{Z}}{r\mathbb{Z}}$ in $Spec(\mathbb{Z}_r)$, with $b_j\not= 1$, identifying strong generators in $\mathbb{Z}_{2n}$, are not accumulation points for the ideals $\mathfrak{a}_i=\frac{(2n-b_i)\mathbb{Z}}{r\mathbb{Z}}$, when $(2n-b_i)=a_i$ is not prime.\footnote{This property is important because it is an a-priori motivation to consider the criterion of Tab. \ref{criterion-to-find-solution-goldbach-conjecture} well found. In fact it shows that the ideals $\mathfrak{b}_j$ cannot be considered on the same footing with respect to the ideals $\mathfrak{a}_i$. In other words the mirror symmetry is necessary to understand how the ideals $\mathfrak{a}_i$ ''converge'' to the ideals $\mathfrak{b}_j$.}
\end{lemma}
\begin{proof}
In fact, each point $\mathfrak{b}_j\in Spec(\mathbb{Z}_r)$ has all the neighborhoods of the type $U=Spec(\mathbb{Z}_r)\setminus \mathfrak{m}$ for some maximal ideal $\mathfrak{m}\in Spec(\mathbb{Z}_r)$. Then if $\mathfrak{a}_i$ is not a maximal ideal it must be contained in some intersection of maximal ideals of $Spec(\mathbb{Z}_r)$, (see the next Lemma \ref{ideals-ai-and-maximal-ideals-in-zr}), say $\mathfrak{a}_i\subset \mathfrak{m}_1\bigcap\cdots\bigcap \mathfrak{m}_k$. Then $\mathfrak{m}$ is an accumulation point of $\mathfrak{a}_i$ if any neighborhood of $\mathfrak{m}$ contains all the ideals $\mathfrak{m}_i$, $1\le i\le k$. But this is impossible ! In fact, for example the neighborhood $U_1=Spec(\mathbb{Z}_r)\setminus \mathfrak{m}_1$ of $\mathfrak{m}$ cannot contain $\mathfrak{m}_1$. \end{proof}
\begin{lemma}[Decompositions of ideals $\mathfrak{a}_i$ in irreducible components]\label{decompositions-ideals-ai-in-irreducible-components}
Each ideal $\mathfrak{a}_i$ in the above set $\{\mathfrak{a}_i\}$, admits an irreducible decomposition into primary ideals of $\mathbb{Z}_r$. If $(2n-b_i)=b_1^{r_1}\cdots b_m^{r_m}$, is the prime factorization of $(2n-b_i)$, then one has the representation (\ref{formula-decompositions-ideals-ai-in-irreducible-components}).
\begin{equation}\label{formula-decompositions-ideals-ai-in-irreducible-components}
  \mathfrak{a}_i=\frac{b_1^{r_1}\mathbb{Z}}{r\mathbb{Z}}\bigcap\cdots\bigcap\frac{b_m^{r_m}\mathbb{Z}}{r\mathbb{Z}}.
\end{equation}
Furthermore, one has $\mathfrak{r}(\mathfrak{a}_i)=b_1\cdots b_m\mathbb{Z}/r\mathbb{Z}$.

If $\mathfrak{a}_i$ is primary then $\mathfrak{r}(\mathfrak{a}_i)$ is prime. If $\mathfrak{a}_i$ is maximal then $\mathfrak{r}(\mathfrak{a}_i)=\mathfrak{a}_i$, i.e., $\mathfrak{a}_i$ is a radical ideal. (The converse is in general false.)
\end{lemma}
\begin{proof}
Let us recall that an ideal $\mathfrak{q}\subset R$ of a ring $R$, is {\em primary} if $xy\in\mathfrak{q}$ implies either $x\in\mathfrak{q}$ or $y^n\in\mathfrak{q}$, for some $n>0$. This is equivalent to say that $R/\mathfrak{q}\not=0$ and every zero-divisor in $R/\mathfrak{q}$ is nilpotent. If $\mathfrak{q}$ is primary then $\mathfrak{r}(\mathfrak{q})$ is the smallest prime ideal containing $\mathfrak{q}$.\footnote{For example if $R=\mathbb{Z}$ the only primary ideals are $<0>$ and $<p^n>$, with $p$ prime.} Furthermore an ideal $\mathfrak{q}\subset R$ is called {\em irreducible} if $\mathfrak{q}=\mathfrak{c}\bigcap \mathfrak{d}$ then $\mathfrak{q}=\mathfrak{c}$ or $\mathfrak{q}=\mathfrak{d}$. Since in a Noetherian ring every ideal is a finite intersection of irreducible ideals, it follows also that each ideals $\mathfrak{a}_i$ admits this decomposition. In fact, if $(2n-b_i)=b_1^{r_1}\cdots b_m^{r_m}$, is the prime factorization of $(2n-b_i)$, then $\mathfrak{a}_i=\frac{(2n-b_i)\mathbb{Z}}{r\mathbb{Z}}$ has the natural decomposition (\ref{formula-decompositions-ideals-ai-in-irreducible-components}), where each ideal $\frac{(b_j^{r_j})\mathbb{Z}}{r\mathbb{Z}}$ is a primary ideal in $\mathbb{Z}_r$. Finally in a Noetherian ring every ideal contains a power of its radical. Therefore, $\mathfrak{r}(\mathfrak{a}_i)\supseteqq \mathfrak{a}_i$ and $\mathfrak{r}(\mathfrak{a}_i)^n\subseteqq \mathfrak{a}_i\subseteqq \mathfrak{r}(\mathfrak{a}_i)$, for some $n>0$. If $\mathfrak{a}_i$ is primary, i.e., $\mathfrak{a}_i=\frac{b^s\mathbb{Z}}{r\mathbb{Z}}$, then $\mathfrak{r}(\mathfrak{a}_i)=\frac{b\mathbb{Z}}{r\mathbb{Z}}=\mathfrak{m}$, a maximal ideal of $\mathbb{Z}_r$ and $\mathfrak{m}\supset\mathfrak{a}_i$. In fact, $\mathfrak{a}_i=\frac{\mathbb{Z}}{r'\mathbb{Z}}$ with $r'b^s=r$ and $\mathfrak{m}=\frac{\mathbb{Z}}{r''\mathbb{Z}}$ with $r''b=r$. Then $r''>r'$ and $r''|r'$, since $r'b^s=r''b$. Therefore, $\mathfrak{a}_i=\mathbb{Z}_{r'}\subset \mathbb{Z}_{r''}=\mathfrak{m}$.
\end{proof}

It is useful in these calculations to utilize the following lemma.
\begin{lemma}[Relations between direct sum, intersection and sum of $\mathbb{Z}$-modules]\label{relations-between-direct-sum-intersection-and-sum-of-z-modules}
Let us consider $\mathfrak{a}_i$ as sub-$\mathbb{Z}$-modules of the $\mathbb{Z}$-module $\mathbb{Z}_r$. Then one has the short exact sequence {\em(\ref{short-exact-sequence-relations-between-direct-sum-intersection-and-sum-of-z-modules})}.
\begin{equation}\label{short-exact-sequence-relations-between-direct-sum-intersection-and-sum-of-z-modules}
   \xymatrix{0\ar[r]&\mathfrak{a}_i\bigcap\mathfrak{a}_j\ar[r]^(0.5){f}&\mathfrak{a}_i\bigoplus\mathfrak{a}_j\ar[r]^(0.5){h_i-h_j}
   &\mathfrak{a}_i+\mathfrak{a}_j\ar[r]&0\\},\, i\not=j
\end{equation}
where $h_i:\mathfrak{a}_i\to\mathfrak{a}_i+\mathfrak{a}_j$ and $h_j:\mathfrak{a}_i\to\mathfrak{a}_i+\mathfrak{a}_j$ are the canonical injections and $f(x)=(x,x)\in\mathfrak{a}_i\bigoplus\mathfrak{a}_j$. Then $(h_i-h_j)(x,y)=x-y$. Furthermore, one has also the short exact sequence {\em(\ref{short-exact-sequence-relations-between-direct-sum-intersection-and-sum-of-z-modules-a})}.

\begin{equation}\label{short-exact-sequence-relations-between-direct-sum-intersection-and-sum-of-z-modules-a}
   \xymatrix{0\ar[r]&\mathbb{Z}_r/(\mathfrak{a}_i\bigcap\mathfrak{a}_j)\ar[r]&(\mathbb{Z}_r/\mathfrak{a}_i)\bigoplus
   (\mathbb{Z}_r/\mathfrak{a}_j)\ar[r]^(0.5){p_i-p_j}
   &\mathbb{Z}_r/(\mathfrak{a}_i+\mathfrak{a}_j)\ar[r]&0\\},\, i\not=j
\end{equation}
where $p_i:\mathbb{Z}_r/\mathfrak{a}_i\to\mathbb{Z}_r/(\mathfrak{a}_i+\mathfrak{a}_j)$ and $p_j:\mathbb{Z}_r/\mathfrak{a}_j\to\mathbb{Z}_r/(\mathfrak{a}_i+\mathfrak{a}_j)$are the canonical projections
\end{lemma}
\begin{proof}
This lemma is a direct application of some standard results in commutative algebra. (See, e.g., \cite{BOURBAKI-2}.)
\end{proof}

\begin{example}
For example by considering the next Example \ref{example-2n-28} relative to the case $\mathbb{Z}_{28}$, hence $\mathbb{Z}_{r=11\cdot 13\cdot 17\cdot 19\cdot 23\cdot 5^2\cdot 3^3}$ one has for $\mathfrak{a}_i=\frac{5\mathbb{Z}}{r\mathbb{Z}}$ and $\mathfrak{a}_j=\frac{9\mathbb{Z}}{r\mathbb{Z}}$, the following $\mathbb{Z}$-modules:
$$\mathfrak{a}_i\bigcap\mathfrak{a}_j=\frac{5\cdot 3^2\mathbb{Z}}{r\mathbb{Z}},\,
\mathfrak{a}_i+\mathfrak{a}_j=\mathbb{Z}_r,\,
\mathfrak{a}_i\bigoplus\mathfrak{a}_j=\frac{5\mathbb{Z}}{r\mathbb{Z}}\bigoplus\frac{9\mathbb{Z}}{r\mathbb{Z}}
$$
to which corresponds the short exact sequence {\em(\ref{short-exact-sequence-relations-between-direct-sum-intersection-and-sum-of-z-modules-b})}.
\begin{equation}\label{short-exact-sequence-relations-between-direct-sum-intersection-and-sum-of-z-modules-b}
  \scalebox{0.8}{$ \xymatrix{0\ar[r]&\mathbb{Z}_{11\cdot 13\cdot 17\cdot 19\cdot 23\cdot 5\cdot 3}\ar@{=}[d]\ar[r]&\mathbb{Z}_{11\cdot 13\cdot 17 \cdot 19\cdot 23 \cdot 5\cdot 3^2}\bigoplus\mathbb{Z}_{11\cdot 13\cdot 17 \cdot 19\cdot 23 \cdot 5^2\cdot 3}\ar@{=}[d]\ar[r]
   &\mathbb{Z}_{11\cdot 13\cdot 17 \cdot 19\cdot 23 \cdot 5^2\cdot 3^3}\ar@{=}[d]\ar[r]&0\\
   0\ar[r]&\mathfrak{a}_i\bigcap\mathfrak{a}_j\ar[r]^(0.5){f}&\mathfrak{a}_i\bigoplus\mathfrak{a}_j\ar[r]^(0.5){h_i-h_j}
   &\framebox{$\mathfrak{a}_i+\mathfrak{a}_j\cong\mathbb{Z}_r$}\ar[r]&0\\}$}
\end{equation}
The sequence {\em(\ref{short-exact-sequence-relations-between-direct-sum-intersection-and-sum-of-z-modules-b})} means that $\mathfrak{a}_i\bigoplus\mathfrak{a}_j$ is larger than $\mathfrak{a}_i+\mathfrak{a}_j\cong \mathbb{Z}_r$. In fact, $\mathfrak{a}_i=\frac{5\mathbb{Z}}{r\mathbb{Z}}$ and $\mathfrak{a}_j=\frac{9\mathbb{Z}}{r\mathbb{Z}}$ are coprime ideals, (situation similar to Lemma \ref{fundamental-properties-of-ideals-of-integers-numbers}), but the corresponding modules $\mathbb{Z}_{r'}$ and $\mathbb{Z}_{r''}$ have $g.c.d.(r',r'')=15\not =1$, hence $r'$ is not coprime of $r''$. In other words $\mathbb{Z}_{r'}\bigoplus\mathbb{Z}_{r''}\ncong\mathbb{Z}_{r'\cdot r''}$. (Lemma \ref{cyclic-group-and-subgroups}.) The kernel of the homomorphism $\mathfrak{a}_i\bigoplus\mathfrak{a}_j\to \mathfrak{a}_i+\mathfrak{a}_j$ is just the intersection $\mathfrak{a}_i\bigcap\mathfrak{a}_j$ that is an ideal of $\mathbb{Z}_r$. The application of the short exact sequence {\em(\ref{short-exact-sequence-relations-between-direct-sum-intersection-and-sum-of-z-modules-b})} gives {\em(\ref{short-exact-sequence-relations-between-direct-sum-intersection-and-sum-of-z-modules-c})}.
\begin{equation}\label{short-exact-sequence-relations-between-direct-sum-intersection-and-sum-of-z-modules-c}
   \xymatrix{0\ar[r]&\mathbb{Z}_{5\cdot 3^2}\ar@{=}[d]\ar[r]&\mathbb{Z}_{5}\bigoplus\mathbb{Z}_{9}\ar@{=}[d]\ar[r]&0\ar@{=}[d]\ar[r]&0\\
   0\ar[r]&\mathbb{Z}_r/(\mathfrak{a}_i\bigcap\mathfrak{a}_j)\ar[r]&(\mathbb{Z}_r/\mathfrak{a}_i)\bigoplus
   (\mathbb{Z}_r/\mathfrak{a}_j)\ar[r]^(0.55){p_i-p_j}
   &\mathbb{Z}_r/(\mathfrak{a}_i+\mathfrak{a}_j)\ar[r]&0\\}
\end{equation}
This just means that $\mathbb{Z}_{45}\cong\mathbb{Z}_{5}\bigoplus\mathbb{Z}_{9}$, i.e., it agrees with Lemma \ref{cyclic-group-and-subgroups} since $5$ and $9$ are coprimes.
\end{example}

\begin{lemma}[Ideals $\mathfrak{a}_i$ and maximal ideals in $\mathbb{Z}_r$]\label{ideals-ai-and-maximal-ideals-in-zr}
Each ideal $\mathfrak{a}_i$ is contained into the intersection of some maximal ideals $\mathfrak{m}_i\subset \mathbb{Z}_r$, hence is contained in some maximal ideal of $\mathbb{Z}_r$.
\end{lemma}
\begin{proof}
Let $\mathfrak{a}_i=\frac{(2n-b_i)\mathbb{Z}}{r\mathbb{Z}}$. If $(2n-b_i)=b_j$, then $\mathfrak{a}_i$ is maximal ! Let us assume that $(2n-b_i)=a_i$ is not prime. Then one has that $\mathfrak{a}_i$ is contained in the intersection of primary ideals. This follows from the fact that $\mathbb{Z}_r$ is a Noetherian ring. However, let us look directly this property. In fact, if $(2n-b_i)=a_i$ admits the following prime factorization $(2n-b_i)=a_i=b_1^{r_1}\cdots b_k^{r_k}$, then one has the following isomorphisms
$$\mathfrak{a}_i=\frac{b_1^{r_1}\cdots b_k^{r_k}\mathbb{Z}}{r\mathbb{Z}}=\frac{b_1^{r_1}\mathbb{Z}}{r\mathbb{Z}}\bigcap\cdots
\bigcap\frac{b_k^{r_k}\mathbb{Z}}{r\mathbb{Z}}=\mathfrak{p}_1\bigcap\cdots\bigcap\mathfrak{p}_k$$
where $\mathfrak{p}_m=\frac{b_m^{r_m}\mathbb{Z}}{r\mathbb{Z}}$, $1\le m\le k$, are primary ideals in $\mathbb{Z}_r$. One has also the following inclusions and isomorphisms:
$$\mathfrak{a}_i\subset\mathfrak{r}(\mathfrak{a}_i)=\frac{b_1\cdots b_k\mathbb{Z}}{r\mathbb{Z}}=
\frac{b_1\mathbb{Z}}{r\mathbb{Z}}\bigcap\cdots
\bigcap\frac{b_k\mathbb{Z}}{r\mathbb{Z}}=\mathfrak{m}_1\bigcap\cdots\bigcap\mathfrak{m}_k.$$
\end{proof}
\begin{example}
Let us look to the case $2n=220$, reported in Tab. \ref{criterion-to-find-solution-goldbach-conjecture}. Let us consider only the ideals $\mathfrak{a}_i=\frac{(2n-b_i)\mathbb{Z}}{r\mathbb{Z}}$ corresponding to the steps reported in Tab. \ref{criterion-to-find-solution-goldbach-conjecture}. Then one has the relation between ideals $\mathfrak{a}_i$ and maximal ideals in $\mathbb{Z}_r$ reported in Tab. \ref{some-exammples-maximal-ideals-containing-ideals-ai-case-2n-220}.
\begin{table}[h]
\caption{Some examples of maximal ideals containing some ideals $\mathfrak{a}_i$ in the case $2n=220$.}
\label{some-exammples-maximal-ideals-containing-ideals-ai-case-2n-220}
\begin{tabular}{|l|}
\hline
{\rm{\footnotesize $\mathfrak{a}_1=(220-211)\mathbb{Z}/r\mathbb{Z}=3^2\mathbb{Z}/r\mathbb{Z}\subset 3\mathbb{Z}/r\mathbb{Z}=\mathfrak{m}_1=\mathfrak{r}(\mathfrak{a}_1)$}}\hfill\\
\hline
{\rm{\footnotesize $\mathfrak{a}_2=(220-199)\mathbb{Z}/r\mathbb{Z}=3\cdot 7\mathbb{Z}/r\mathbb{Z}= (3\mathbb{Z}/r\mathbb{Z})\bigcap(7\mathbb{Z}/r\mathbb{Z})=\mathfrak{m}_1\bigcap\mathfrak{m}_2=
\mathfrak{r}(\mathfrak{a}_2)$}}\hfill\\
\hline
{\rm{\footnotesize $\mathfrak{a}_3=(220-197)\mathbb{Z}/r\mathbb{Z}=23\mathbb{Z}/r\mathbb{Z}= \mathfrak{m}_6=\mathfrak{r}(\mathfrak{a}_3)$}}\hfill\\
\hline
\multicolumn{1}{l}{\rm{\footnotesize See also Tab. \ref{criterion-to-find-solution-goldbach-conjecture}.}}\hfill\\
\multicolumn{1}{l}{\rm{\footnotesize $r=l.c.m.(a_i)$.}}\hfill\\
\multicolumn{1}{l}{\rm{\footnotesize $\mathbb{Z}_{220}^{\times}=\{1,a_i\, |\, 1<a_i<220=2^2\cdot 5\cdot 11,\, g.c.d.(220,a_i)=1\}$.}}\hfill\\
\multicolumn{1}{l}{\rm{\footnotesize Maximal ideals in $\mathbb{Z}_{r}$: $\{\mathfrak{m}_i\}=\{
\frac{3\mathbb{Z}}{r\mathbb{Z}},\frac{7\mathbb{Z}}{r\mathbb{Z}},\frac{13\mathbb{Z}}{r\mathbb{Z}},\frac{17\mathbb{Z}}{r\mathbb{Z}},
\frac{19\mathbb{Z}}{r\mathbb{Z}},\frac{23\mathbb{Z}}{r\mathbb{Z}},\cdots\}$.}}\hfill\\
\end{tabular}
\end{table}
\end{example}
\begin{example}
In Tab. \ref{exammples-maximal-ideals-containing-ideals-ai-case-2n-28} are reported the relations between ideals $\mathfrak{a}_i$ and maximal ideals $\mathfrak{m}_j$ in $\mathbb{Z}_r$, for the case $2n=28$, considered also in the next Example \ref{example-2n-28}.
\begin{table}[h]
\caption{Examples of maximal ideals containing ideals $\mathfrak{a}_i$ in the case $2n=28$.}
\label{exammples-maximal-ideals-containing-ideals-ai-case-2n-28}
\begin{tabular}{|l|}
\hline
{\rm{\footnotesize $\mathfrak{a}_1=(2n-23)\mathbb{Z}/r\mathbb{Z}=5\mathbb{Z}/r\mathbb{Z}= \mathfrak{m}_2=\mathfrak{r}(\mathfrak{a}_1)$}}\hfill\\
\hline
{\rm{\footnotesize $\mathfrak{a}_2=(2n-19)\mathbb{Z}/r\mathbb{Z}=3^2\mathbb{Z}/r\mathbb{Z}\subset \mathfrak{m}_1=\mathfrak{r}(\mathfrak{a}_2)$}}\hfill\\
\hline
{\rm{\footnotesize $\mathfrak{a}_3=(2n-17)\mathbb{Z}/r\mathbb{Z}=11\mathbb{Z}/r\mathbb{Z}= \mathfrak{m}_3=\mathfrak{r}(\mathfrak{a}_3)$}}\hfill\\
\hline
{\rm{\footnotesize $\mathfrak{a}_4=(2n-13)\mathbb{Z}/r\mathbb{Z}=3\cdot 5\mathbb{Z}/r\mathbb{Z}= \mathfrak{m}_1\bigcap\mathfrak{m}_2=\mathfrak{r}(\mathfrak{a}_4)$}}\hfill\\
\hline
{\rm{\footnotesize $\mathfrak{a}_5=(2n-11)\mathbb{Z}/r\mathbb{Z}=17\mathbb{Z}/r\mathbb{Z}= \mathfrak{m}_5=\mathfrak{r}(\mathfrak{a}_5)$}}\hfill\\
\hline
{\rm{\footnotesize $\mathfrak{a}_6=(2n-5)\mathbb{Z}/r\mathbb{Z}=23\mathbb{Z}/r\mathbb{Z}= \mathfrak{m}_7=\mathfrak{r}(\mathfrak{a}_6)$}}\hfill\\
\hline
{\rm{\footnotesize $\mathfrak{a}_7=(2n-3)\mathbb{Z}/r\mathbb{Z}=5^2\mathbb{Z}/r\mathbb{Z}\subset \mathfrak{m}_2=\mathfrak{r}(\mathfrak{a}_7)$}}\hfill\\
\hline
\multicolumn{1}{l}{\rm{\footnotesize See also Example \ref{example-2n-28}.}}\hfill\\
\multicolumn{1}{l}{\rm{\footnotesize $r=l.c.m.(a_i)=11\cdot 13\cdot 17\cdot 19\cdot 23\cdot 25\cdot 27$.}}\hfill\\
\multicolumn{1}{l}{\rm{\footnotesize $\mathbb{Z}_{28}^{\times}=\{1,a_i\, |\, 1<a_i<28=2^2\cdot 7,\, g.c.d.(28,a_i)=1\}$.}}\hfill\\
\multicolumn{1}{l}{\rm{\footnotesize Maximal ideals in $\mathbb{Z}_{r}$:}}\hfill\\
\multicolumn{1}{l}{\rm{\footnotesize $\{\mathfrak{m}_i\}=\{
\frac{3\mathbb{Z}}{r\mathbb{Z}},\frac{5\mathbb{Z}}{r\mathbb{Z}},\frac{11\mathbb{Z}}{r\mathbb{Z}},\frac{13\mathbb{Z}}{r\mathbb{Z}},
\frac{17\mathbb{Z}}{r\mathbb{Z}},\frac{19\mathbb{Z}}{r\mathbb{Z}},\frac{23\mathbb{Z}}{r\mathbb{Z}}\}$.}}\hfill\\
\end{tabular}
\end{table}

\end{example}

In order to conclude this proof, i.e., to assure that in the set $\{\mathfrak{a}_i\}$ are also included maximal ideals of the type $\mathfrak{m}=b\mathbb{Z}/r\mathbb{Z}\subset\mathbb{Z}_r$, the following lemma give the definite answer.

\begin{lemma}[Goldbach bordism and Goldbach couples]\label{goldbach-bordism-and-goldbach-couples}
$\bullet$\hskip 2pt For $n=1,2,3$, one has $\{\mathfrak{a}_i\}=\varnothing$.

$\bullet$\hskip 2pt When the following conditions occur:

{\em(i)} $n\not=$ prime;

{\em(ii)} $2n-1\not=$ prime;

{\em(iii)} $n>3$;

the set of ideals $\{\mathfrak{a}_i\}$ contains a maximal ideal of $\mathbb{Z}_r$ at least.

$\bullet$\hskip 2pt  For any $n>3$, all Goldbach couples, other eventual trivial and Noether-Goldbach ones, can be identified by means of maximal ideals occuring in $\{\mathfrak{a}_i\}$.
\end{lemma}

\begin{proof}
$\bullet$\hskip 2pt For $n=1,2,3$, one has $\{\mathfrak{a}_i\}=\varnothing$, since the only Goldbach couples are trivial and Noether-Goldbach ones. (See Example \ref{n=1}, Example \ref{n=2} and Example \ref{n=3}.)

$\bullet$\hskip 2pt Let us consider the natural embeddings $\mathbb{Z}\to\mathbb{R}\to\mathbb{R}^2$, $n\mapsto n\mapsto(n,0)$. Then we say that a couple of points $a,\, b\, \in\mathbb{R}^2$ are {\em $2n$-Goldbach bording} if there exists a smooth curve $\gamma:[0,1]\to\mathbb{R}^2$ such that $\gamma(0)=a$, $\gamma(1)=b$ and $\gamma$ intersects the segment $[\gamma(0),\gamma(1)]\subset\mathbb{R}$, contained in the straight-line $\mathbb{R}\subset\mathbb{R}^2$, identified by the two points $a,\, b\, \in\mathbb{R}^2$, into a couple of integers $\bar a,\, \bar b\subset\mathbb{Z}\subset\mathbb{R}$, such that $(\bar a,\bar b)$ is a Goldbach couple with respect to an even integer $2n$, up to diffeomorphisms of $\mathbb{R}^2$. In particular $(\bar a,\bar b)$ can be also a trivial Goldbach couple, i.e., $(\bar a,\bar a)$, with $\bar a=n$ prime. Let us denote by ${}^{2n}\Omega_{GB}$ the {\em $2n$-Goldbach bordism group}. Let us prove that ${}^{2n}\Omega_{GB}=\mathbb{Z}_2$, i.e., any two points in $\mathbb{R}^2$ are $2n$-Goldbach bording, for $n\ge 1$. We shall use the following Lemma.
\begin{lemma}\label{zero-unoriented-bordism-group-plane}
The unoriented $0$-bordism group of $\mathbb{R}^2$, is $\Omega_0(\mathbb{R}^2)\cong\mathbb{Z}_2$. Furthermore, any diffeomorphism $f:\mathbb{R}^2\to\mathbb{R}^2$ induces an isomorphism $f_*:\Omega_0(\mathbb{R}^2)\cong\Omega_0(\mathbb{R}^2)$.
\end{lemma}
\begin{proof}
In fact, one has
$$\Omega_0(\mathbb{R}^2)\cong H_0(\mathbb{R}^2;\mathbb{Z}_2)\bigotimes_{\mathbb{Z}_2}\Omega_0
\cong\mathbb{Z}_2\bigotimes_{\mathbb{Z}_2}\mathbb{Z}_2\cong\mathbb{Z}_2.$$
$\Omega_0$ denotes the unoriented $0$-bordism group.\footnote{For details on the (co)bordism group theory see, e.g., Refs. \cite{ATIYAH,MADSEN-MILGRAM,PONTRJAGIN,QUILLEN,RUDYAK,STONG,SWITZER,THOM,WALL}.} Furthermore, for any diffeomorphism $f:\mathbb{R}^2\to\mathbb{R}^2$, we get the induced homomorphism $f_*:\Omega_0(\mathbb{R}^2)\to\Omega_0(f(\mathbb{R}^2))$, given by $f_*:[a]\mapsto[f(a)]$. On the other hand
$$\Omega_0(f(\mathbb{R}^2))\cong H_0(f(\mathbb{R}^2);\mathbb{Z}_2)\bigotimes_{\mathbb{Z}_2}\Omega_0
\cong H_0(\mathbb{R}^2;\mathbb{Z}_2)\bigotimes_{\mathbb{Z}_2}\mathbb{Z}_2\cong\Omega_0(\mathbb{R}^2).$$
\end{proof}

In the short exact sequence (\ref{short-exact-sequence-2-gb-non-oriented-bordism}) is reported the relation between ${}^{2}\Omega_{GB}$ and non-oriented $0$-bordism in $\mathbb{R}^2$.
\begin{equation}\label{short-exact-sequence-2-gb-non-oriented-bordism}
  \xymatrix{ 0\ar[r]&\ker({}^{2}b)\ar[r]&{}^{2}\Omega_{GB}\ar[r]^(0.4){{}^{2}b}&\framebox{$\Omega_0(\mathbb{R}^2)\cong\mathbb{Z}_2$}\ar[r]&0  \\}
 \end{equation}
 In fact, given two points $a=(x_1,y_1),\, b=(x_2,y_2)\in\mathbb{R}^2$, we can assume that they are on the $x$-axis and by a further transformation to assume that $a=(0,0)$ and $b=(2,0)$. Then the curve $y=\sin(\hat x)$ with $\hat x=x/\pi$, passes for $\hat x=0$ for $a$ and for $\hat x=2$ for $b$ and has a further zero at $\hat x=1$. This identifies the trivial Goldbach couple $(1,1)$ for the even number $2n=2$. Therefore the two points $a$ and $b$ are $2$-Goldbach bording. Since these are arbitrary points in $\mathbb{R}^2$, it follows that $\ker({}^{2}b)=0$, hence ${}^{2}\Omega_{GB}\cong\Omega_0(\mathbb{R}^2)$. This conclusion can be generalized to any even number $2n$. In fact one can consider the mapping ${}^{2n}f:x\mapsto \hat x=x/n\pi$, that transforms the interval $[0,2n\pi]$ into $[0,2]$ and the points $a=0$ and $b=2n\pi$ into the points $\hat x=0$ and $\hat x=2$ respectively. Furthermore the point $x=\pi n$ is transformed into $\hat x=1$. Then the curve $y=\sin(\hat x)$ realizes again the $2$-Goldbach bordism. On the other hand one has the induced homomorphism ${}^{2n}\Omega_{GB} \to {}^{2}\Omega_{GB}$ on the Goldbach bordism groups.
 More precisely one has the exact commutative diagram (\ref{commutative-exact-diagram-2n-gb-non-oriented-bordism}).\footnote{Warn ! A priori we cannot know the structure of the $2n$-Goldbach bordism group, but the mapping ${}^{2n}f$ allows us to relate it to the $2$-Goldbach bordism group. In fact, the mapping ${}^{2n}f:\mathbb{R}^{2}\to\mathbb{R}^{2}$ is a diffeomorphism, since it is represented by the functions $({}^{2n}f^1,{}^{2n}f^2):(x,y)=(x^1,x^2)\mapsto (\hat x=\frac{x}{\pi n},y)$. The corresponding jacobian matrix is $(\partial x_k.{}^{2n}f^j)=\scalebox{0.7}{$\left(
                              \begin{array}{cc}
                                1/\pi n & 0 \\
                                0 & 1 \\
                              \end{array}
                            \right)$}$ with $\det(\partial x_k.{}^{2n}f^j)=1/\pi n\not=0$. Thus we can consider the inverse diffeomorphism ${}^{2n}f^{-1}:\mathbb{R}^{2}\to\mathbb{R}^{2}$ and we necessarily get ${}^{2n}\Omega_{GB} \cong {}^{2n}f^{-1}_*( {}^{2}\Omega_{GB})$. }

 \begin{equation}\label{commutative-exact-diagram-2n-gb-non-oriented-bordism}
    \xymatrix{&0&0&\\
    0\ar[r]&{}^{2}\Omega_{GB}\ar[u]\ar[r]^(0.4){{}^{2}b}&\framebox{$\Omega_0(\mathbb{R}^2) \cong\mathbb{Z}_2$}\ar[u]\ar[r]&0\\
  0\ar[r]&{}^{2n}\Omega_{GB}\ar[u]^(0.5){{}^{2n}f_*}\ar[r]^(0.4){{}^{2n}b}&\framebox{$\Omega_0(\mathbb{R}^2) \cong\mathbb{Z}_2$}\ar[u]^(0.5){{}^{2n}f_*}\ar[r]&0\\
  &0\ar[u]&0\ar[u]&\\}
\end{equation}
The homomorphism ${}^{2n}f_*:{}^{2n}\Omega_{GB}\to{}^{2}\Omega_{GB}$ is an isomorphism for construction and the homomorphism ${}^{2n}b:{}^{2n}\Omega_{GB}\to \Omega_0(\mathbb{R}^2)$ is also an isomorphism since ${}^{2}\Omega_{GB}\cong\Omega_0(\mathbb{R}^2) \mathop{\cong}\limits^{{}^{2n}f_*}\Omega_0(\mathbb{R}^2) $ are isomorphisms. Therefore, one has the isomorphism ${}^{2n}\Omega_{GB}\cong\mathbb{Z}_2$. The identification of all Goldbach bordism groups ${}^{2n}\Omega_{GB}$, $n> 1$, with ${}^{2}\Omega_{GB}$, it could erroneously focus the importance of the Goldbach splitting $2=1+1$. In fact, one can also choice $2m$, with $m=2,\, 3,\, 4,\, \dots$, where one can easily calculate ${}^{2m}\Omega_{GB}$ and see that ${}^{2m}\Omega_{GB}\cong\mathbb{Z}_2$, for some $m>1$, and $m$ prime. Then, we can prove that ${}^{2(m+r)}\Omega_{GB}\cong\mathbb{Z}_2$, $r>0$, too. In fact, since the $2(m+r)$-Goldbach bordism group is determined up to diffeomorphisms of $\mathbb{R}^2$, we can deform the curve $y=\sin(\hat x)$, $\hat x\in[0,2(m+r)]\subset\mathbb{R}$, into a curve on the interval $[0,2m]\subset\mathbb{R}$, that intersects the $x$-axis at the integer $m$. Since $m$ is assumed prime, it follows that $(m,m)$ is the trivial Goldbach couple for $2m$. This is equivalent to say that ${}^{2(m+r)}\Omega_{GB}$ is isomorphic to ${}^{2m}\Omega_{GB}\cong\mathbb{Z}_2$. The situation is resumed in the commutative diagram (\ref{relation-goldbach-bordism-groups-n-prime}), where $m>1$, prime and $r>0$.
\begin{equation}\label{relation-goldbach-bordism-groups-n-prime}
    \xymatrix{{}^{2(m+r)}\Omega_{GB}\ar@{=}[d]^{\wr}\ar@{=}[r]^(0.7){\sim}&\cdots\ar@{=}[r]^(0.4){\sim}&{}^{2m}\Omega_{GB}\ar@{=}[d]^{\wr}\ar@{=}[r]^(0.6){\sim}&\cdots\ar@{=}[r]^{\sim}&
    {}^{6}\Omega_{GB}\ar@{=}[d]^{\wr}\ar@{=}[r]^{\sim}&{}^{4}\Omega_{GB}\ar@{=}[d]^{\wr}\ar@{=}[r]^{\sim}&{}^{2}\Omega_{GB}\ar@{=}[d]^{\wr}\\
    \mathbb{Z}_2\ar@{=}[r]&\cdots\ar@{=}[r]&\mathbb{Z}_2\ar@{=}[r]&\cdots\ar@{=}[r]&\mathbb{Z}_2\ar@{=}[r]&\mathbb{Z}_2
    \ar@{=}[r]&\mathbb{Z}_2\\}
\end{equation}

The existence of the group ${}^{2n}\Omega_{GB}\cong \mathbb{Z}_2$ proves that for any fixed $2n$ there exists at least a Goldbach couple. Then, as a by product, we get that in the set of ideals $\{\mathfrak{a}_i\}$, and under the above conditions (i)-(ii)-(iii),  there exists a maximal ideal at least, $\mathfrak{m}=b\mathbb{Z}/r\mathbb{Z}$, with $b=2n-b_i$ prime and $b_i$ a prime of the interval $]1,2n[\subset\mathbb{N}$ identifying a strong generator in $\mathbb{Z}_{2n}$. In fact, under conditions (i)-(ii)-(iii), a Goldbach couple cannot be a trivial one and neither a Noether-Goldbach couple. Furthermore, since we have assumed $n>3$, $\{\mathfrak{a}_i\}$ cannot be empty. The remaining of the lemma directly follows from above results and previous lemmas.
\end{proof}

In order to avoid any possible confusion, let us apply the proof to some significative examples.

\begin{example}[$\framebox{$2n=2$}$]\label{n=1}
In this case there exists the canonical Noether-Goldbach couple $(1,2n-1)=(1,1)$. This is a trivial Goldbach couple, since $n=1$ is prime. The set of ideals $\{\mathfrak{a}_i\}=\varnothing$.
\end{example}

\begin{example}[$\framebox{$2n=4$}$]\label{n=2}
In this case the only Goldbach couples are the canonical Noether-Goldbach couple $(1,2n-1)=(1,3)$, and the trivial Goldbach couple $(2,2)$ since $n=2$ is prime. The set of ideals $\{\mathfrak{a}_i\}=\varnothing$.
\end{example}

\begin{example}[$\framebox{$2n=6$}$]\label{n=3}
In this case there exists the canonical Noether-Goldbach couple $(1,2n-1)=(1,5)$. Since $n=3$ is prime there exists also the trivial Goldbach couple $(3,3)$. Then $\{\mathfrak{a}_i\}=\varnothing$. Do not exist other Goldbach couples.
\end{example}

\begin{example}[$\framebox{$2n=8$}$]
In this case there exists the canonical Noether-Goldbach couple $(1,2n-1)=(1,7)$. Instead does not exist the trivial Goldbach couple since $n=4$ is even. $r=l.c.m.(3,5)=15$.
 $$\{\mathfrak{a}_i\}=\{\frac{(2n-5)\mathbb{Z}}{15\mathbb{Z}},\frac{(2n-3)\mathbb{Z}}{15\mathbb{Z}}\}=
 \{\frac{3\mathbb{Z}}{15\mathbb{Z}},\frac{5\mathbb{Z}}{15\mathbb{Z}}\}.$$
 The set of maximal ideals is $$\{\mathfrak{m}_{i}\}=\{\mathfrak{m}_{1}=\frac{3\mathbb{Z}}{15\mathbb{Z}},\mathfrak{m}_{2}=\frac{5\mathbb{Z}}{15\mathbb{Z}}\},$$ to which corresponds the same Goldbach couple $(3,5)$. By summarizing, the Goldbach couples in $\mathbb{Z}_{8}$ are $(1,7)$ and $(3,5)$.
\end{example}

\begin{example}[$\framebox{$2n=10$}$]
In this case does not exist the canonical Noether-Goldbach couple, since $2n-1=9$ is not a prime number, but there exists the trivial Goldbach couple $(5,5)$, since $n$ is prime. $r=l.c.m.(3,7)=21$.
 $$\{\mathfrak{a}_i\}=\{\frac{(2n-7)\mathbb{Z}}{21\mathbb{Z}},\frac{(2n-3)\mathbb{Z}}{21\mathbb{Z}}\}=
 \{\frac{3\mathbb{Z}}{21\mathbb{Z}},\frac{7\mathbb{Z}}{21\mathbb{Z}}\}.$$
 The set of maximal ideals is $$\{\mathfrak{m}_{i}\}=\{\mathfrak{m}_{1}=\frac{3\mathbb{Z}}{21\mathbb{Z}},\mathfrak{m}_{2}=\frac{7\mathbb{Z}}{21\mathbb{Z}}\}.$$ To $\mathfrak{m}_1$ and $\mathfrak{m}_2$ corresponds the same Goldbach couple $(3,7)$. By summarizing the Goldbach couple in $\mathbb{Z}_{10}$ is $(3,7)$. To this must be added the trivial Goldbach couple $(5,5)$ that does not come from units in $\mathbb{Z}_{10}$.
\end{example}

\begin{example}[$\framebox{$2n=28$}$]\label{example-2n-28}
In this case there does not exist a canonical Noether-Goldbach couple, since $2n-1=27$ is not a prime number, and neither does there exist a trivial Goldbach couple, since $n=14$ is not prime. $$r=l.c.m.(3,5,9,11,13,15,17,19,23,25,27)=11\cdot 13\cdot 17\cdot 19\cdot 23\cdot 25\cdot 27=717084225.$$
 $$\left\{
 \begin{array}{ll}
   \{\mathfrak{a}_i\}&=\{\frac{(2n-23)\mathbb{Z}}{r\mathbb{Z}},
 \frac{(2n-19)\mathbb{Z}}{r\mathbb{Z}},\frac{(2n-17)\mathbb{Z}}{r\mathbb{Z}},
 \frac{(2n-13)\mathbb{Z}}{r\mathbb{Z}},
 \frac{(2n-11)\mathbb{Z}}{r\mathbb{Z}},
 \frac{(2n-5)\mathbb{Z}}{r\mathbb{Z}},
 \frac{(2n-3)\mathbb{Z}}{r\mathbb{Z}}\}\\
 &\\
   &= \{\frac{5\mathbb{Z}}{r\mathbb{Z}},
 \frac{9\mathbb{Z}}{r\mathbb{Z}},\frac{11\mathbb{Z}}{r\mathbb{Z}},
 \frac{15\mathbb{Z}}{r\mathbb{Z}},
 \frac{17\mathbb{Z}}{r\mathbb{Z}},
 \frac{23\mathbb{Z}}{r\mathbb{Z}},
 \frac{25\mathbb{Z}}{r\mathbb{Z}}\}.\\
 \end{array}
 \right.$$

 The set of maximal ideals of $\mathbb{Z}_r$, belonging to the set $\{\mathfrak{a}_i\}$ is the following:\footnote{Compare with Tab. \ref{exammples-maximal-ideals-containing-ideals-ai-case-2n-28}.} $$\{\mathfrak{m}_{i}\}=\{\mathfrak{m}_{2}=\frac{5\mathbb{Z}}{r\mathbb{Z}},
 \mathfrak{m}_{3}=\frac{11\mathbb{Z}}{r\mathbb{Z}},
 \mathfrak{m}_{5}=\frac{17\mathbb{Z}}{r\mathbb{Z}},
 \mathfrak{m}_{7}=\frac{23\mathbb{Z}}{r\mathbb{Z}}\}.$$ To $\mathfrak{m}_2$ and $\mathfrak{m}_7$ corresponds the same Goldbach couple $(5,23)$ and to $\mathfrak{m}_{3}$ and $\mathfrak{m}_{5}$ there corresponds the Goldbach couple $(11,17)$.\end{example}

In conclusion, in the ring $\mathbb{Z}_{2n}$ there exists a canonical Goldbach couple, and this can be found by means of the criterion in Tab. \ref{criterion-to-find-solution-goldbach-conjecture}. The same criterion allows us to find also all the other Goldbach couples. Therefore, the proof of the Theorem \ref{goldbach-strong-generators-cyclic-group-2n} is done !
\end{proof}

After Theorem \ref{goldbach-strong-generators-cyclic-group-2n} we have the following corollaries.

\begin{cor}[Goldbach Conjecture]\label{goldbach-conjecture}
Any even integer $2n$, $n\ge 1$, can be split into the sum of two primes $p_1$ and $p_2$: $2n=p_1+p_2$.\footnote{Let us emphasize that $n$ can be any integer $\ge 1$. In fact, if $n$ is a prime number, it is trivial that $2n$ is the sum of two primes: $2n=n+n$.}
\end{cor}

\begin{cor}[Restricted Goldbach Conjecture]\label{restricted-goldbach-conjecture}
Any even integer $2n$, $n>1$, can be split into the sum of two primes $p_1$ and $p_2$: $2n=p_1+p_2$.\footnote{Let us underline that the GC in its original form considers $1$ as a prime number ! More recently, a restricted version of GC, (say RGC), has been proposed by assuming the restricted prime numbers set $P^\bullet=P\setminus\{1\}$ only, and even numbers $2n$, with $n>1$. With this respect, let us emphasize that our proof of the GC can be adapted also to prove the RGC. In fact, this is the GC with the additional restriction that the Noether-Goldbach couples cannot be considered as acceptable solutions. However, Noether-Goldbach couples are found simply by looking to the fact if $2n-1$ are primes or not. So, identified at the beginning, in the process of the criterion of Tab. \ref{criterion-to-find-solution-goldbach-conjecture}. Really, this criterion becomes interesting just when do not exist the canonical Noether-Goldbach couples. In fact, the maximal ideal $\mathfrak{m}=b\mathbb{Z}/r\mathbb{Z}$ considered in the proof of the GC, associated to the set of ideals $\mathfrak{a}_i=(2n-b_i)\mathbb{Z}/r\mathbb{Z}\subset \mathbb{Z}_r$, necessarily has $b\not=1$.}
\end{cor}

\begin{cor}[Criterion to find the Goldbach couples for any fixed $2n$, $n\ge 1$]\label{criterion}
The following steps give us a criterion to find all the Goldbach couples for any fixed integer $2n$, $n\ge 1$.

{\em 1)} If $n$ is prime there exists the trivial Goldbach couple $(n,n)$.

{\em 2)} If $2n-1$ is prime there exists the Noether-Goldbach couple $(1,2n-1)$.

{\em 3)} If $n>3$ all the other possible Goldbach couples are identified by means of maximal ideals of $\mathbb{Z}_r$ belonging to the set $\{\mathfrak{a}_i\}$.\footnote{Here the meaning of symbols are the ones just considered in Lemma \ref{goldbach-bordism-and-goldbach-couples}.}
\end{cor}

\section[Applications]{\bf Applications}

In this section we give some applications interesting the classical Euclidean geometry and the quantum algebra in the sense introduced by A. Pr\'astaro. (See \cite{PRAS1,PRAS2} and related Pr\'astaro's works quoted therein.)

\begin{proposition}[Goldbach triangle]\label{goldbach-triangle}
In a circle $\Gamma$ of radius $n\in \mathbb{N}$, there exists an inscribed right triangle $ABC$, with hypothenus $AB$ passing for the centre $O$ of $\Gamma$, such that the projection $H$ of the vertex $C$ on $AB$, divides $AB$ into two segments $AH$ and $HB$ of length respectively $p_1$ and $p_2$, prime numbers.\footnote{For details on this geometric interpretation of the GC see \cite{NAMBIAR}, where it is emphasized the equivalence of the GC and the solution of the following Diophantine equation: $n^2=a^2+b^2$, where $n$, $a$ and $b$ are three integers such that $a=p_1\, p_2$ and $2b=p_2-p_1$, with $p_1$ and $p_2$, prime numbers. This relates the GC to a {\em Fermat like theorem}. Let us recall that in 1900, David Hilbert proposed the solvability of all Diophantine problems as the tenth of his celebrated problems. However, after 70 years has been published a result in mathematical logic that in general Diophantine problems are unsolvable. (Matiyasevich's theorem \cite{MATIYASEVICH}.)
Therefore, this proof of the Goldbach's conjecture is also an encouragement for mathematicians to solve problems, even if their solutions could have fat chance according to some general statement in mathematical logic ! (See also \cite{HAZEWINKEL} for general information on Diophantine equations and \cite{SUN} for the undecibility of these equations.)

Warn! Since the proof of the GC (or RGC) given in this paper, is obtained by methods of the Commutative Algebra and Algebraic Topology, one could consider yet this conjecture as an example of the G\"odel's incompleteness theorem. In other words, the GC (or RGC) is a true proposition in Arithmetic, but not provable in Arithmetic. The criterion in Tab. \ref{criterion-to-find-solution-goldbach-conjecture} works well in Arithmetic, but its proof is beyond Arithmetic.}
\end{proposition}

In the following we give an application of the GC to the quantum algebra, in the sense of A. Pr\'astaro \cite{PRAS2}.

\begin{theorem}[Quantum algebraic interpretation of the Goldbach conjecture]\label{quantum-algebraic-interpretation-goldbach-conjecture}
$\bullet$\hskip 2pt Let $A$ be a quantum algebra in the sense of A. Pr\'astaro, then there exists the canonical homomorphism {\em(\ref{canonical-homomorphism-quantum-algebra-c-c})}, {\em(quantum-Goldbach-homomorphism)}.
\begin{equation}\label{canonical-homomorphism-quantum-algebra-c-c}
\left\{\begin{array}{l}
 g_*:2\mathbb{Z}\bigotimes_{\mathbb{Z}} A\to \mathbb{Z}_{2}\bigotimes_{\mathbb{Z}} A\bigoplus \mathbb{Z}_{2}\bigotimes_{\mathbb{Z}} A\\
 g_*(2n\otimes a)= ([p_1]\otimes a,[p_2]\otimes a) \in(1\otimes a,1\otimes a)\\
\end{array}\right.
\end{equation}
where $(p_1,p_2)$ is the Goldbach couple identified by the criterion reported in Tab. \ref{criterion-to-find-solution-goldbach-conjecture} and codified by Theorem \ref{goldbach-strong-generators-cyclic-group-2n}.
We call $2\mathbb{Z}\bigotimes_{\mathbb{Z}} A$ the (additive) {\em group of quantum extended even-numbers}.
Furthermore one has the commutative diagram {\em(\ref{commutative-diagram-g-c-quantum-algebra})}, with exact vertical line.
\begin{equation}\label{commutative-diagram-g-c-quantum-algebra}
    \xymatrix{&0\ar[d]\\
        &2\mathbb{Z}\bigotimes_{\mathbb{Z}}A\ar[d]^{b}\ar[dl]^{g_*}\\
    {\framebox{$\mathbb{Z}_{2}\bigotimes_{\mathbb{Z}} A\bigoplus \mathbb{Z}_{2}\bigotimes_{\mathbb{Z}} A$}}\ar[dr]_{+}&{\framebox{$\mathbb{Z}\bigotimes_{\mathbb{Z}}A\cong A$}}\ar[d]^{c}\\
    &\mathbb{Z}_2\bigotimes_{\mathbb{Z}}A\ar[d]\\
    &0\\}
\end{equation}
One has the canonical isomorphisms reported in {\em(\ref{canonical-isomorphim-in-quantum-algebra})}.
\begin{equation}\label{canonical-isomorphim-in-quantum-algebra}
\left\{\begin{array}{l}
  {\rm im}(b)\cong 2\mathbb{Z}\bigotimes_{\mathbb{Z}}A\cong\ker(c)\\
  {\rm im}(c)\cong \mathbb{Z}\bigotimes_{\mathbb{Z}}A/2\mathbb{Z}\bigotimes_{\mathbb{Z}}A\cong \mathbb{Z}_2\bigotimes_{\mathbb{Z}}A\\
\end{array}
\right.
\end{equation}
$\bullet$\hskip 2pt The quantum-Goldbach-homomorphism gives a relation between number theory, crystallographic groups and integral bordism groups of PDEs and quantum PDEs.
\end{theorem}

\begin{proof}
Let us first consider the following free resolution of the $\mathbb{Z}$-module $\mathbb{Z}_2$:
$$\xymatrix{0\ar[r]&\mathbb{Z}\ar[r]^{2}&\mathbb{Z}\ar[r]&\mathbb{Z}_2\ar[r]&0}$$
By tensoring this sequence with a quantum algebra $A$, considered as a $\mathbb{Z}$-module by means of the embeddings $\mathbb{Z}\to \mathbb{K}\to A$, where $\mathbb{K}=\mathbb{R}$, or $\mathbb{K}=\mathbb{C}$, we get the exact sequence (\ref{exact-sequence-induced-via-tor-from-exact-sequence}).
\begin{equation}\label{exact-sequence-induced-via-tor-from-exact-sequence}
   \scalebox{0.7}{$\xymatrix{0\ar[r]&Tor^{\mathbb{Z}}(A;\mathbb{Z})\ar@{=}[d]\ar[r]&Tor^{\mathbb{Z}}(A;\mathbb{Z})\ar@{=}[d]\ar[r]&
   Tor^{\mathbb{Z}}(A;\mathbb{Z}_2)\ar@{=}[d]\ar[r]
   &A\bigotimes_{\mathbb{Z}}\mathbb{Z}\ar@{=}[d]\ar[r]&
   A\bigotimes_{\mathbb{Z}}\mathbb{Z}\ar[r]\ar@{=}[d]&A\bigotimes_{\mathbb{Z}}\mathbb{Z}_2\ar[r]\ar@{=}[d]&0\\
   0\ar[r]&0\ar[r]&0\ar[r]&Tor^{\mathbb{Z}}(A;\mathbb{Z}_2)\ar[r]&A\ar[r]^(0.5){2}&A\ar[r]&A\bigotimes_{\mathbb{Z}}\mathbb{Z}_2\ar[r]&0\\}
$}\end{equation}
From the bottom horizontal line, we can calculate $Tor^{\mathbb{Z}}(A;\mathbb{Z}_2)=\ker(\xymatrix{A\ar[r]^{2}&A\\})$. Since $A$ is a $\mathbb{K}$-vector space, it follows that $\ker(\xymatrix{A\ar[r]^{2}&A\\})=\{0\}$.
Similarly, by working with the following free resolution of $\mathbb{Z}$-module $\mathbb{Z}_2$:
$$\xymatrix{0\ar[r]&2\mathbb{Z}\ar[r]&\mathbb{Z}\ar[r]&\mathbb{Z}_2\ar[r]&0}$$
we get the vertical exact sequence in (\ref{commutative-diagram-g-c-quantum-algebra}). This is connected with the quantum-Goldbach-homomorphism. In fact, we have $+\circ g_*=c\circ b$. Then the isomorphisms reported in (\ref{canonical-isomorphim-in-quantum-algebra}) are directly obtained from standard algebraic considerations on the vertical exact sequence in (\ref{commutative-diagram-g-c-quantum-algebra}).

Finally the quantum Goldbach homomorphism allows us to represent the group of quantum extended even-numbers into a quantum extension of the crystallographic group $p4m=\mathbb{Z}^2\, \rtimes\, D_4$. In fact, $D_4=\mathbb{Z}_2\times\mathbb{Z}_2$ is the point group of $p4m$. On the other hand $\mathbb{Z}_2\bigoplus\mathbb{Z}_2$ can be interpreted also as integral bordism groups of some PDEs. (See \cite{PRAS1,PRAS2} and some Pr\'astaro's works, quoted therein on the relation between integral bordism groups of PDEs and quantum PDEs and crystallographic groups.)
\end{proof}

\end{document}


\setcounter{page}{\pagenumber}
\newcommand{\T}{\mathbb{T}}
\newcommand{\R}{\mathbb{R}}
\newcommand{\Q}{\mathbb{Q}}
\newcommand{\N}{\mathbb{N}}
\newcommand{\Z}{\mathbb{Z}}
\newcommand{\tx}[1]{\quad\mbox{#1}\quad}
\parindent=0pt
\def\SRA{\hskip 2pt\hbox{$\joinrel\mathrel\circ\joinrel\to$}}
\def\tbox{\hskip 1pt\frame{\vbox{\vbox{\hbox{\boldmath$\scriptstyle\times$}}}}\hskip 2pt}
\def\circvert{\vbox{\hbox to 8.9pt{$\mid$\hskip -3.6pt $\circ$}}}
\def\IM{\hbox{\rm im}\hskip 2pt}
\def\COIM{\hbox{\rm coim}\hskip 2pt}
\def\COKER{\hbox{\rm coker}\hskip 2pt}
\def\TR{\hbox{\rm tr}\hskip 2pt}
\def\GRAD{\hbox{\rm grad}\hskip 2pt}
\def\RANK{\hbox{\rm rank}\hskip 2pt}
\def\MOD{\hbox{\rm mod}\hskip 2pt}
\def\DEN{\hbox{\rm den}\hskip 2pt}
\def\DEG{\hbox{\rm deg}\hskip 2pt}

\title[The Landau's problems.II: Landau's problems solved]{THE LANDAU's PROBLEMS.II:\\ LANDAU's PROBLEMS SOLVED}

\author{Agostino Pr\'astaro}
\maketitle
\vspace{-.5cm}
{\footnotesize
\begin{center}
Department SBAI - Mathematics, University of Rome La Sapienza, Via A.Scarpa 16,
00161 Rome, Italy. \\
E-mail: {\tt agostino.prastaro@uniroma1.it}
\end{center}
}
\vskip 0.5cm
\centerline{\em This paper is dedicated to G. Cantor, S. Ramanujan and R. Thom.}
\vskip 0.5cm
\begin{abstract}
Three of the well known {\em four Landau's problems} are solved in this paper. (In \cite{PRAS1} the proof of the Goldbach's conjecture has been already given.)
\end{abstract}

\vskip 0.5cm

\noindent {\bf AMS Subject Classification:} 11R04; 11T30; 11D99; 11U05; 55N22; 55N20; 81R50; 81T99; 20H15.

\vspace{.08in} \noindent \textbf{Keywords}: Landau's problems. Goldbach's conjecture; Algebraic number theory; Diophantine equations; Cobordism groups.

\section[Introduction]{\bf Introduction}\label{introduction-section}
\vskip 0.5cm

This is the second part of a work devoted to the four problems in Number Theory, called {\em Landau's problems}, listed by Edmund Landau at the 1912 International Congress of Mathematicians. These problems are the following. 1. Goldbach's conjecture. 2. Twin prime conjecture. 3. Legendre's conjecture. 4. Are there infinitely many primes $p$ such that $p-1$ is a perfect square ? In \cite{PRAS1} the proof of the Goldbach's conjecture has been already given. In this paper we show that by utilizing some algebraic topologic methods introduced in \cite{PRAS1}, some Landau's problems can be proved too. Furthermore, for the above fourth Landau's problem a Euler-Riemann zeta function estimate is given and settled the problem negatively by evaluating the cardinality of the set of solutions of a suitable Diophantine equation of Ramanujan-Nagell-Lebesgue type.
\vskip 0.5cm

\section[The twin prime conjecture]{\bf The twin prime conjecture proved}\label{section-twin}
\vskip 0.5cm

In this section we shall prove that the twin conjecture is true, as a particular case of the more general de Polignac's conjecture. This has been possible thanks to our proof of the Goldbach's conjecture given in \cite{PRAS1}.\footnote{For complementary information on (co)bordism group theory and Diophantine equations see, e.g., \cite{ATIYAH, HAZEWINKEL, STONG, SWITZER, THOM, WALL}.}

\begin{definition}\label{twin-prime}
A {\em twin prime} is a prime number $p$ that differs from another prime number $q$ by two: $p-q=2$.
\end{definition}
\begin{conj}\label{twin-conjecture}
The {\em twin prime conjecture} states that there are many primes $q$ such that $q+2$ is prime too.
\end{conj}

\begin{conj}[de Polignac's conjecture (1849)]\label{de-polignac-conjecture}
For every natural number $n\in\mathbb{N}$, there are infinitely many prime pairs $p$ and $q$ such that  $p-q=2n$.

(For the case $n=1$ the de Polignac's conjecture reduces to the twin conjecture.)
\end{conj}

In the following we list some well-known theorems about twin conjecture and de Polignac's conjecture. (For details see, e.g., \href{http://en.wikipedia.org/wiki/Twin_prime_conjecture}{\tt wikipedia/Twin-prime-conjecture} and references quoted therein.)

\begin{proposition}[Brun's theorem (1915)]\label{brun-theorem}
The sum of reciprocals of the twin primes is convergent. More precisely the number of twin primes less than $N$ does not exceed $CN/(\lg N)^2$, for some absolute constant $C>0$.
\end{proposition}

\begin{proposition}[Some estimates]\label{some-estimates}

{\em(1) (Paul Erd\"os (1940))} There is a constant $c<1$ and infinitely many primes $p$ such that $(q-p)<(c\ln p)$, where $q$ denotes the next prime after $p$.

{\em(2) (Helmut Maier (1986))} The constant $c<0.25$.

{\em(3) (Daniel Goldston and Cem Yidmm (2004))} The constant $c=0.085786$.

{\em(4) (Daniel Goldston, J\'anos Pintz and Cem Yidmm (2005))} The constant can be chosen to be arbitrarily small, i.e. $\mathop{\lim}\limits_{n\to\infty}\inf \frac{p_{n+1}-p_n}{\lg p_n}=0$.

\end{proposition}

\begin{definition}\label{isolated-prime}
An {\em isolated prime} is a prime number $p$ such that neither $p-2$ or $p+2$ is prime.
\end{definition}

\begin{example}\label{examples-isolated-prime}
The first {\em isolated primes} is listed in the following:
$$2,\, 23,\, 37,\, 47,\, 53,\, 67,\, 79,\, 83,\, 89,\, 97,\, ...$$
\end{example}
\begin{theorem}[The pre-Polignac's conjecture proved]\label{pre-polignac-conjecture-proved}
For every even positive integer $2n$, $n\ge 1$, we can find a prime $0< q<2n$ and a prime $p>2n$ such that $p-q=2n$.
\end{theorem}

\begin{proof}
Let us consider the following lemmas.
\begin{lemma}\label{lemma-a}
If $p$ and $q$ are two primes that satisfy conditions in Theorem \ref{pre-polignac-conjecture-proved}, then $q$ cannot divide $2n$.
\end{lemma}
\begin{proof}
In fact, if $q$ divides $2n$, then $p=2n+q=\bar q\, q+q=q(\bar q+1)$, hence $p$ could not be prime.
\end{proof}

\begin{lemma}\label{lemma-b}
Under the same hypotheses of Lemma \ref{lemma-a}, $q$ identifies a strong generator of $\mathbb{Z}_{2n}$, i.e., $q\in \mathbb{Z}_{2n}^{\blacksquare}\subset \mathbb{Z}_{2n}^{\times}\subset \mathbb{Z}_{2n}$. (We use the same notation introduced in \cite{PRAS1} and, for abuse of notation, we denote by $q$ yet its projection on $\mathbb{Z}_{2n}$, by means of the canonical epimorphism $\pi:\mathbb{Z}\to \mathbb{Z}_{2n}$.)
\end{lemma}
\begin{proof}
It follows directly from Lemma \ref{lemma-a}.
\end{proof}
Let us, now, prove that in the set $[0,2n]$ there exists at least one prime $q$ such that $2n+q$ is prime too. For this we shall follows a strategy similar to one used in \cite{PRAS1} to analogous situations.
\begin{definition}[de Polignac couples]\label{de-polignac-couples}
 We define {\em $2n$-de Polignac couple} a pair of positive primes $(p,q)$ such that $p-q=2n$.
 \end{definition}
 \begin{definition}[de Polignac bordism]\label{de-polignac-bordism}
We say that two points $a,\, b\in\mathbb{R}^2$, are {\em $2n$-de Polignac bording}, if there exists a smooth curve $\gamma:[0,1]\to\mathbb{R}^2$, such that $\gamma(0)=a$, $\gamma(1)=b$ and $\gamma$ intersects the line $\overline{ab}$ into a $2n$-Polignac couple up to diffeomorphisms of $\mathbb{R}^2$. Let us denote by ${}^{2n}\Omega_{deP}$ the corresponding bordism group that we call {\em $2n$-de Polignac bordism group}.
\end{definition}
\begin{lemma}\label{lemma-c}
One has the canonical isomorphism ${}^{2n}\Omega_{deP}\cong \mathbb{Z}_2$.
\end{lemma}
\begin{proof}
Let us consider first the case $2n=2$, i.e., $n=1$. In such a case we can consider that any couple of points $a,\, b\in\mathbb{R}^2$, can be identified respectively with the points $(0,0)\in\mathbb{R}^2$ and $(4,0)\in\mathbb{R}^2$. Then the curve $y=\sin(\hat x)$, taken over the interval $\hat x=\frac{x}{\pi}\in[0,4]$, intersects the $\hat x-axis$ in all the integers between $0$ and $4$. In particular the points $1$ and $3$ are a $2$-de Polignac couple, since $3-1=2$.\footnote{\hskip 3pt Warn ! In this paper, we consider the integer $1$ as prime. However, our proof works well also considering $2$ as the first prime. In such a case, it is enough to start the proof by considering the function $y=\sin(\hat x)$ on the interval $[0,8]$, (instead than $[0,4]$), in order to prove that ${}^{2}\Omega_{deP}\cong \mathbb{Z}_2$.}\label{warn1}

Let us generalize this situation to any $2n$, $n>1$. Then we can consider the curve $y=\sin(\hat x)$, in the interval $[0,2n\cdot 2]=[0,4n]$. We claim that in such interval there are two prime numbers $q\in[0,2n]$ and $p\in[2n,4n]$, such that $p-q=2n$. In fact, since the curve $y=\sin(\hat x)$, intersects the $\hat x$-axis in all the integers between $0$ and $4n$, it follows that there exists a $2n$-de Polignac couple if ${}^{2n}\Omega_{deP}$ can be identified with ${}^{2}\Omega_{deP}$ by means of a suitable diffeomorphism of $\mathbb{R}^2$. This diffeomorphism is the following: $f:(\hat x,y)\mapsto (\frac{\hat x}{n},y)$. In fact by means of such a diffeomorphism the interval $[0,4n]$ is transformed into $[0,4]$ and the curve $y=\sin(\hat x)$ is deformed in the curve $y=\sin(\frac{\hat x}{n})$. This last intersects the interval $[0,4]$ into the $2$-de Polignac couple $(1,3)$. Therefore the curve $y=\sin(\hat x)$  must necessarily intersect $[0,4n]$ into a $2n$-de Polignac couple, since we get the bijection $ {}^{2n}\Omega_{deP}\cong {}^{2}\Omega_{deP}\cong \mathbb{Z}_2$.\footnote{\hskip 3pt Warn ! Here the proof has been made considering $1$ as prime. However, as just emphasized in the previous note, we can consider $2$ as the first prime. Then the proof works by considering the function $y=\sin(\hat x)$ on the interval $[0,8n]$ (instead that $[0,4n]$). In this way we can relate the bordism group ${}^{2n}\Omega_{deP}$ with ${}^{2}\Omega_{deP}$, by deforming the function $y=\sin(\hat x)$ on the interval $[0,8]$, by means the same deformation used above. Then this deformed curve passes for the two couples of points $\{3,5\}$ and $\{5,7\}$ that are $2$-de Polignac couples.}\label{warn2}
\end{proof}
In this way the proof of the theorem is done.
\end{proof}
\begin{theorem}[The de Polignac's conjecture proved]\label{polignac-conjecture-proved}
The de Polignac conjecture is true.
\end{theorem}

\begin{proof}
After Theorem \ref{pre-polignac-conjecture-proved} it is enough to prove that for any fixed $2n$, we can further find $2n$-de Polignac couples by extending the process followed in the proof of Theorem \ref{pre-polignac-conjecture-proved}, i.e., by considering the curve $y=\sin(\hat x)$ on the interval $[0,2n.2^{m}]$, $m>2$. In fact we can deform such a curve into one on the interval  $[0,2n\cdot 2^{m-1}]$, by considering the diffeomorphism $f:\mathbb{R}^2\to\mathbb{R}^2$, $(\hat x,y)\mapsto (\frac{\hat x}{2},y)$. By starting with $n=2$ we can prove that there exists the $2n$-de Polignac couple $(p,q)$ with $q\in[0,2n\cdot 2]$ and $p\in[2n\cdot 2,2n\cdot 2^2]$, just considering the process illustrated in the proof of Theorem \ref{pre-polignac-conjecture-proved}. Then by iterating on $m$ we can find new $2n$-de Polignac couples. Thus for any $2n$ the set of $2n$-de Polignac couples is necessarily infinite.
\end{proof}

\begin{example}
In Tab. \ref{examples-de-polignac-couples} are reported some examples of $2n$-de Polignac couples $(p,q)$, $q\in[0,2n\cdot 2^m]$, $m\ge 1$, found in correspondence of some increasing values of $n$ and $m$.
\end{example}
\begin{table}[t]
\caption{Examples of $2n$-de Polignac couples $(p,q)$ with $q\in[0,2n\cdot 2^m]$, $1\le n\le 5$, $n=10$ and $1\le m\le 4$.}
\label{examples-de-polignac-couples}
\scalebox{0.6}{$\begin{tabular}{|c|c|c|c|c|c|}
\hline
\hfil{\rm{\footnotesize $n$}}\hfil&\hfil{\rm{\footnotesize $2n$}}\hfil&\hfil{\rm{\footnotesize $m=1$}}\hfil&\hfil{\rm{\footnotesize $m=2$}}\hfil&\hfil{\rm{\footnotesize $m=3$}}\hfil&\hfil{\rm{\footnotesize $m=4$}}\hfil\\
\hline\hline
{\rm{\footnotesize $1$}}\hfil&\hfil{\rm{\footnotesize $2$}}&\hfil{\rm{\footnotesize $(1,3)$}}&{\rm{\footnotesize $(3,5)$, $(5,7)$}}\hfill&{\rm{\footnotesize $(11,13)$}}\hfill&{\rm{\footnotesize $(17,19)$, $(29,31)$}}\hfill\\
\hline
{\rm{\footnotesize $2$}}\hfil&\hfil{\rm{\footnotesize $4$}}&\hfil{\rm{\footnotesize $(3,7)$}}&{\rm{\footnotesize $(7,11)$}}\hfill&{\rm{\footnotesize $(13,17)$, $(19,23)$}}\hfill&{\rm{\footnotesize $(43,47)$, $(37,41)$}}\hfill\\
\hline
{\rm{\footnotesize $3$}}\hfil&\hfil{\rm{\footnotesize $6$}}&\hfil{\rm{\footnotesize $(1,7)$, $(5,11)$}}&{\rm{\footnotesize $(17,23)$, $(13,19)$, $(11,17)$, $(7,13)$}}\hfill&{\rm{\footnotesize $(41,47)$, $(37,43)$, $(31,37)$, $(23,29)$}}\hfill&{\rm{\footnotesize $(83,89)$, $(73,79)$, $(61,67)$, $(53,59)$, $(47,53)$}}\hfill\\
\hline
{\rm{\footnotesize $4$}}\hfil&\hfil{\rm{\footnotesize $8$}}&\hfil{\rm{\footnotesize $(3,11)$, $(5,13)$}}&{\rm{\footnotesize $(11,19)$, $(23,31)$}}\hfill&{\rm{\footnotesize $(29,37)$, $(53,61)$}}\hfill&{\rm{\footnotesize $(101,109)$, $(89,97)$, $(71,79)$, $(59,67)$}}\hfill\\
\hline
{\rm{\footnotesize $5$}}\hfil&\hfil{\rm{\footnotesize $10$}}&\hfil{\rm{\footnotesize $(1,11)$, $(3,13)$, $(7,17)$}}&{\rm{\footnotesize $(19,29)$, $(13,23)$}}\hfill&{\rm{\footnotesize $(61,71)$, $(43,53)$, $(37,47)$}}\hfill&{\rm{\footnotesize $(139,149)$, $(127,137)$}}\hfill\\
&&&&&{\rm{\footnotesize $(103,113)$, $(79,89)$, $(73,83)$}}\hfill\\
\hline
{\rm{\footnotesize $10$}}\hfil&\hfil{\rm{\footnotesize $20$}}&\hfil{\rm{\footnotesize $(3,23)$, $(11,31)$, $(17,37)$}}&{\rm{\footnotesize $(41,61)$, $(47,67)$, $53,73)$, $59,79)$}}\hfill&{\rm{\footnotesize $(83,103)$, $(107,127)$, $(131,151)$, $(137,157)$}}\hfill&{\rm{\footnotesize $(173,193)$, $(191,211)$, $(251,271)$}}\hfill\\
&&&&&{\rm{\footnotesize $(257,277)$, $(263,283)$, $(293,313)$}}\hfill\\
\hline
\multicolumn{6}{l}{\rm{\footnotesize The $2n$-case with $n=1$ corresponds to the twin conjecture.}}\\
\multicolumn{6}{l}{\rm{\footnotesize Infinite $2n$-de Polignac couples are obtained since we can take any integer $m\ge 1$.}}\\
\multicolumn{6}{l}{\rm{\footnotesize For $n=10$ we can see that continuing on $m>4$, there is also the $20$-de Polignac couple $(317,337)$, where $317\in[0,320]$, considered for $1\le m\le 4$.}}\\
\multicolumn{6}{l}{\rm{\footnotesize Warn ! In this paper we consider $1$ as prime. However our proof works well also by considering $2$ as the first prime. (See footnotes at page \pageref{warn1} and \pageref{warn2}.)}}\\
\end{tabular}$}
\end{table}

\begin{cor}[The twin conjecture proved]\label{twin-conjecture-proved}
The twin conjecture is true.
\end{cor}

\section{\bf The Legendre's conjecture proved}\label{section-legendre-conjecture-proveed}

In this section we shall prove that the Legendre's conjecture is true. The method utilized follows the same line adopted to prove the Goldbach's conjecture in \cite{PRAS1}.

\begin{conj}\label{legendre-conjecture}
Legendre's conjecture, proposed by Adrien-Marie Legendre, states that there is a prime number between $n^2$ and $(n+1)^2$ for every positive integer $n$.
\end{conj}

In this section we shall prove the Conjecture \ref{legendre-conjecture}, i.e., we have the following theorem.

\begin{theorem}[The Legendre's conjecture]\label{legendre-conjecture-proved}
There is a prime number between $n^2$ and $(n+1)^2$, for any positive integer $n$.
\end{theorem}

\begin{proof}
Let us, first, remark that if there exists a prime between $n^2+1$ and $(n+1)^2$, there exists in advance a prime between $n^2$ and $(n+1)^2$, hence in order to prove the Legendre conjecture, it is enough to prove that there exists a prime between $n^2+1$ and $(n+1)^2$. Let us emphasize that $(n+1)^2-(n^2+1)=2n$.  Now, from our proof of the Goldbach conjecture in \cite{PRAS1}, we get that there exists at least a prime in $[0,2n]$. With this respect, let us say that two points $a,\, b\in\mathbb{R}^2$ are {\em $n$-Legendre bording} if there exists a smooth curve $\gamma:[0,1]\to\mathbb{R}^2$, such that $\gamma(0)=a$, $\gamma(1)=b$ and this curve intersects the line $\overline{ab}$ into a prime $p$ contained between $n^2$ and $(n+1)^2$, up to diffeomorphisms of $\mathbb{R}^2$. We shall utilize the following lemma.

\begin{lemma}\label{diffeomorphisms-and-integer-distances}
Two generic points $a,\, b\in\mathbb{R}^2$ can be identified respectively with the points $n^2+1$ and $(n+1)^2$ of the $x$-axis, by means suitable diffeomorphisms of $\mathbb{R}^2$.
\end{lemma}
\begin{proof}
In fact, by means of an affine rigid transformation of $f_1:\mathbb{R}^2\to \mathbb{R}^2$, identified by an affine transformation of coordinates $(x,y)=(x^i)\mapsto(\tilde x^i=A^i_kx^k+u^i)$, we can identify the point $(n^2+1)=(n^2+1,0)$ with $a$, (hence the point $(n+1)^2=((n+1)^2,0)$ into a point of the line $\overline{ab}$). This affine transformation includes also a translation in the direction $\overline{ab}$ sending in $a$ the origin of the $\tilde x^i$-coordinates. Therefore in such a set of coordinates $(\tilde x^i)$, one has the following identification: $a=(0,0)=(n^2+1)$, $((n+1)^2)=(2n,0)$ and $b=(\bar b,0)$, where $\bar b$ denotes the Euclidean canonical distance of $b$ from $a$, that after the affine rigid transformation $f_1$ is conserved. The next diffeomorphism of $\mathbb{R}^2$ is a contraction $f_2$ sending the point $b$ onto the point $((n+1)^2)$. $f_2(\tilde x^1,\tilde x^2)=(\bar x^1,\bar x^2)$, with $\bar x^1=\tilde x^1.\gamma$, $\gamma=\frac{2n}{{\bar b}}$. Therefore in $\bar x^k$-coordinates $a$ and $b$ have the coordinates $(0,0)$ and $(2n,0)$ respectively. In this way the distance between $a$ and $b$ is the even integer $2n$, characterizing the original distance between $(n^2+1)$ and $((n+1)^2)$.
\end{proof}

Since by Lemma \ref{diffeomorphisms-and-integer-distances} we can transform, by suitable diffeomorphisms of $\mathbb{R}^2$,  $(n^2+1)$ into $a$ and $(n+1)^2$ into $b$, it follows the curve $\bar y=\sin (\hat x)$ with $\hat x=\bar x/\pi$, considered over the interval $[a,b]$, identified with $[0,2n]$, intersects such interval in at least one primes. This follows from our proof of the Goldbach conjecture \cite{PRAS1}. As a by product there exists a prime between $n^2$ and $(n+1)^2$. In other words the $n$-Legendre bordism group ${}^n\Omega_{Leg}$ of $\mathbb{R}^2$ can be identified with $\mathbb{Z}_2$, i.e., one has the isomorphism ${}^n\Omega_{Leg}\cong \mathbb{Z}_2$. This can be seen also step by step, considering ${}^n\Omega_{Leg}$ for $n=2,\, 3,\, 4,\, \cdots$, where it is easy to calculate this type of bordism group. For example from Tab. \ref{examples-legendre-conjecture} we can see that ${}^1\Omega_{Leg}\cong \mathbb{Z}_2$ since any couple of points $a,\, b\in\mathbb{R}^2$ can be identified with $(n^2+1=2)$ and $((n+1)^2=4)$ respectively, distance $2$. The interval $[2,4]$ can be deformed into $[0,2]$, on the which we can consider the curve $y=\sin(\hat x)$, that intersects the $\hat x$-axis at the prime $\hat x=1$. Therefore in the $[2,4]$ there exists a prime too. In fact there is the prime $3\in[2,4]$. As a by product there is also a prime in the interval $[1,4]$: this is $3$ (other than $2$). Next, by passing to the case $n=2$, we can identify $a,\, b\in\mathbb{R}^2$ with $(n^2+1=5)$ and $((n+1)^2=9)$ respectively, distance $4$. The interval $[5,9]$ can be deformed into $[0,4]$, on the which we can consider the curve $y=\sin(\hat x)$, that intersects the $\hat x$-axis at the prime $\hat x=2$ and $\hat x=3$. Therefore in the $[5,9]$ there exists a prime too. In fact there is the prime $7\in[5,9]$. As a by product there is also a prime in the interval $[4,9]$: this is $3$ (other than $5$ and $7$).
Iterating this process on $n$, we can get the commutative diagram (\ref{n-legendre-bordism-groups-commutative-diagram}).\footnote{Warn ! Even if in this paper we consider $1$ as prime, the proof of the Goldbach conjecture proved in (\cite{PRAS1}), works well also considering $2$ as the first prime. With this respect, we can state that our proof of the Legendre's conjecture is true also by omitting $1$ in the set of primes. (Compare with results reported in Tab. \ref{examples-legendre-conjecture}.)}\label{warn3}

\begin{equation}\label{n-legendre-bordism-groups-commutative-diagram}
    \xymatrix{{}^{n}\Omega_{Leg}\ar@{=}[d]^{\wr}\ar@{=}[r]^(0.6){\sim}&\cdots\ar@{=}[r]^{\sim}&
    {}^{3}\Omega_{Leg}\ar@{=}[d]^{\wr}\ar@{=}[r]^{\sim}&{}^{2}\Omega_{Leg}\ar@{=}[d]^{\wr}\ar@{=}[r]^{\sim}&{}^{1}\Omega_{Leg}\ar@{=}[d]^{\wr}\\
   \mathbb{Z}_2\ar@{=}[r]&\cdots\ar@{=}[r]&\mathbb{Z}_2\ar@{=}[r]&\mathbb{Z}_2
    \ar@{=}[r]&\mathbb{Z}_2\\}
\end{equation}
\end{proof}

\begin{example}
In Tab. \ref{examples-legendre-conjecture} are reported some examples of primes $p\in[n^2,(n+1)^2]$.
\end{example}

\begin{table}[t]
\caption{Examples of primes $p\in[n^2,(n+1)^2]$.}
\label{examples-legendre-conjecture}
\scalebox{0.8}{$\begin{tabular}{|c|c|c|c|}
\hline
\hfil{\rm{\footnotesize $n$}}\hfil&\hfil{\rm{\footnotesize $n^2$}}\hfil&\hfil{\rm{\footnotesize $(n+1)^2$}}\hfil&\hfil{\rm{\footnotesize prime $p\in[n^2,(n+1)^2]$}}\hfil\\
\hline\hline
{\rm{\footnotesize $1$}}\hfil&\hfil{\rm{\footnotesize $1$}}&\hfil{\rm{\footnotesize $4$}}&{\rm{\footnotesize $1$, $2$, $3$}}\hfill\\
\hline
{\rm{\footnotesize $2$}}\hfil&\hfil{\rm{\footnotesize $4$}}&\hfil{\rm{\footnotesize $9$}}&{\rm{\footnotesize $5$, $7$}}\hfill\\
\hline
{\rm{\footnotesize $3$}}\hfil&\hfil{\rm{\footnotesize $9$}}&\hfil{\rm{\footnotesize $16$}}&{\rm{\footnotesize $11$, $13$}}\hfill\\
\hline
{\rm{\footnotesize $4$}}\hfil&\hfil{\rm{\footnotesize $16$}}&\hfil{\rm{\footnotesize $25$}}&{\rm{\footnotesize $17$, $19$, $23$}}\hfill\\
\hline
{\rm{\footnotesize $5$}}\hfil&\hfil{\rm{\footnotesize $25$}}&\hfil{\rm{\footnotesize $36$}}&{\rm{\footnotesize $29$, $31$}}\hfill\\
\hline
{\rm{\footnotesize $6$}}\hfil&\hfil{\rm{\footnotesize $36$}}&\hfil{\rm{\footnotesize $49$}}&{\rm{\footnotesize $37$, $41$, $43$, $47$}}\hfill\\
\hline
{\rm{\footnotesize $7$}}\hfil&\hfil{\rm{\footnotesize $49$}}&\hfil{\rm{\footnotesize $64$}}&{\rm{\footnotesize $53$, $59$, $61$}}\hfill\\
\hline
{\rm{\footnotesize $8$}}\hfil&\hfil{\rm{\footnotesize $64$}}&\hfil{\rm{\footnotesize $81$}}&{\rm{\footnotesize $67$, $71$, $73$, $79$}}\hfill\\
\hline
{\rm{\footnotesize $9$}}\hfil&\hfil{\rm{\footnotesize $81$}}&\hfil{\rm{\footnotesize $100$}}&{\rm{\footnotesize $83$, $89$, $97$}}\hfill\\
\hline
{\rm{\footnotesize $10$}}\hfil&\hfil{\rm{\footnotesize $100$}}&\hfil{\rm{\footnotesize $121$}}&{\rm{\footnotesize $101$, $103$, $107$, $109$, $113$}}\hfill\\
\hline
\multicolumn{4}{l}{\rm{\footnotesize Warn ! In this paper we consider $1$ as prime.}}\\
\multicolumn{4}{l}{\rm{\footnotesize However our proof works well also by considering $2$ as the first prime. (See footnote at page \pageref{warn3}).}}\\
\end{tabular}$}
\end{table}

Theorem \ref{legendre-conjecture-proved} can be reformulated in terms of triangle numbers.
\begin{definition}
One defines {\em triangle numbers} the following ones $$T_n=\sum_{1\le k\le n}k=1+2+3+\cdots+n=\frac{n(n+1)}{2}=\binom{n+1}{n},\, n\ge 1.$$
\end{definition}
\begin{lemma}[Properties of triangle numbers]\label{properties-triangle-numbers}
One has the following properties:

{\em(a)} $T_{a+b}=T_a+T_b+a\cdot b$;

{\em(b)} $T_{a\cdot b}=T_a\cdot T_b+T_{a-1}\cdot T_{b-1}$;

{\em(c)} $T_{n}+T_{n-1}=(T_n-T_{n-1})^2=n^2$;\footnote{Note that we can also write $n^2=\sum_{1\le k\le n}(2k-1)$. In fact $T_{n}+T_{n-1}$ is equal to the sum of all the first $n$ odd numbers.}

{\em(d) (Gauss's Eureka result)(1796)} Any number $n\in\mathbb{N}$ can be written as a sum of at most three triangular numbers.\footnote{For example $n^2+1=T_{n}+T_{n-1}+T_1$, as it results from Lemma \ref{properties-triangle-numbers}(c) and from the fact that $T_1=1$.}

{\em(e) (Triangle numbers and Bernoulli's numbers)} $T_n=\frac{1}{2}(B_0\, n^2+2B_1\, n^1)=S_1(n)$, with the following Bernoulli numbers $B_0=1$ and $B_1=\frac{1}{2}$ and $S_m(n)=\sum_{1\le k\le m}k^m$, given by the {\em Bernoulli's formula}
$$S_m(n)=\frac{1}{m+1}\sum_{0\le k\le m}\binom{m+1}{k}\, B_k\, n^{m+1-k}.$$

{\em(f) (Triangle numbers that are also square)} The following recurrence relation gives a non-exhaustive generation of triangular numbers that are also square: $$S(n+1)=4\cdot S(n)\cdot[8\cdot S(n)+1],\, S(1)=1.$$

\end{lemma}

\begin{theorem}[Legendre conjecture and triangle numbers]\label{legendre-conjecture-triangle-numbers}
Between $T_{n-1}+T_{n}$ and $T_{n}+T_{n+1}$ there is a prime $p\in P$.
\end{theorem}
\begin{proof}
In fact, $T_{n-1}+T_{n}=n^2$ and $T_{n}+T_{n+1}=(n+1)^2$, hence from Theorem \ref{legendre-conjecture-proved} it follows soon the proof.
\end{proof}

\section[On the infinite set of ghost right-triangles and parabolic primes]{\bf On the set of ghost right-triangles and parabolic primes}\label{the-set-ghost-right-triangles-parabolic-primes-section}

In this section we shall consider the set $P_{\square}$ of all primes $p\in P$ such that $p-1$ is a perfect square. We identify $P_{\square}$ with a set of particular right-triangles, that we call {\em ghost right-triangles}. A goal is a characterization of $P_{\square}$ by means of the Euler-Riemann zeta function, $\zeta$, that offers good chances in order to consider $P_{\square}$ an infinite set. Furthermore it is proved that $\sharp(P_{\square})<\aleph_0$ as set of solutions of a suitable Diophantine equation. Therefore, the main result of this section is the solution of the fourth Landau's problem answering in the negative to the question there contained.

Let us start with some metric characterizations of $\mathbb{N}^2$.

\begin{proposition}[The discrete metric spaces $\mathbb{Z}^2$ and $\mathbb{N}^2$]\label{some-discrete-metric-spaces}
$\mathbb{Z}^2$ endowed with the induced metric topology of $\mathbb{R}^2$, by means of the inclusion $\mathbb{Z}^2\hookrightarrow \mathbb{R}^2$, is a discrete metric topologic space. There is the Euclidean metric $d:\mathbb{Z}^2\to \mathbb{R}$, given by $d((a,b),(\bar a,\bar b))=\sqrt{(\bar a-a)^2+(\bar b-b)^2}\in\mathbb{R}$. This topology restricted to $\mathbb{N}^2\subset\mathbb{Z}^2$, induces an analogous metric topology on $\mathbb{N}^2$.
\end{proposition}
\begin{definition}\label{metric-regular-line}
We call {\em metric-regular line} a line $\gamma\subset \mathbb{R}^2$ bording two points $A,\, B\in\mathbb{N}^2$ if its length in the Euclidean metric of $\mathbb{R}^2$, is an integer: $d_{\gamma}(A,B)\in\mathbb{N}$.
\end{definition}
\begin{definition}\label{ghost-line}
We call {\em ghost line} a line $\gamma\subset \mathbb{R}^2$ bording two points $A,\, B\in\mathbb{N}^2$ if its length in the Euclidean metric of $\mathbb{R}^2$, is not an integer: $d_{\gamma}(A,B)\not\in\mathbb{N}$.
\end{definition}

\begin{example}
Let us recall that a {\em Pythagorean right-triangle} is a right-triangle with integers side-lengths, identified by a {\em Pythagorean triple} $(a,b,c)$, such that the Diophantine equation: $a^2+b^2=c^2$, $a,\, b,\, c\in\mathbb{N}$. The set of Pythagorean right-triangles is infinite, since if $(a,b,c)$ is a Pythagorean triple, then also $(m\cdot a,m\cdot b,m\cdot c)$ is a Pythagorean triple for any positive integer $m\in\mathbb{N}$.
In a Pythagorean right-triangle and in a primitive Pythagorean triangle, two consecutive vertices are connected by a metric-regular line.\footnote{A {\em primitive Pythagorean triangle} is a Pythagorean right-triangle identified by a triple $(a,b,c)$, where $a$, $b$ and $c$ are pairwise coprime.}
\end{example}

\begin{example}
Are always metric-regular lines the ones lying on the lattice identified by $\mathbb{N}^2$.
\end{example}

\begin{definition}[Ghost triangles in $\mathbb{N}^2$]\label{ghost-triangles}
We call {\em ghost triangles} in $\mathbb{N}^2$, triangles $ABC\subset\mathbb{R}^2$, having vertexes in $\mathbb{N}^2$, such that at least one side is a ghost line.
\end{definition}

\begin{definition}[$n$-Ghost right-triangles]\label{n-ghost-right-triangles}
We call {\em $n$-ghost right-triangle} a right-triangle $ABC$, with legs $\overline{AC}=n\in \mathbb{N}$ and $\overline{CB}=1$ respectively, such that the Diophantine equation: $n^2+1=p$ is satisfied with $p\in P$, where $P\subset\mathbb{N}$ is the subset of primes. We denote
$$P_{\square}=\{p\in P\, |\, p=n^2+1;\, n\in\mathbb{N}\}\subset P$$
the set of ghost right-triangles.\footnote{$P_{\square}=\{2,5,17,\cdots\}\subset P$ coincides with the subset of primes $p$ such that $p-1$ is a perfect square. One can also define $P_{\square\bullet}=\{n\in \mathbb{N}\, |\, n^2+1=p\in P\}\subset \mathbb{N}$. This is not a new set, since $P_{\square\bullet}\leftrightarrow P_{\square}$. Therefore, for abuse of notation we shall denote $P_{\square\bullet}$ with $P_{\square}$ yet.}
\end{definition}

\begin{table}[t]
\caption{$n$-ghost-right-triangles: $1\le n\le 60$, $n^2+1=p\in P$.}
\label{tab-examples-ghost-right-triangles}
\scalebox{0.6}{$\begin{tabular}{|c|c|c|c|c|}
\hline
\hfil{\rm{\footnotesize $n$}}\hfil&\hfil{\rm{\footnotesize $p=n^2+1$}}\hfil&\hfil{\rm{\footnotesize ($\square$)=ghost right-triangle}}\hfil&\hfil{\rm{\footnotesize primes between two consecutive ($\square$)}}\hfill\\
\hline\hline
{\rm{\footnotesize $1$}}\hfil&\hfil{\rm{\footnotesize $p=1+1=2$}}\hfill&\hfil{\rm{\footnotesize ($\square$)}}\hfil&\hfill\\
\hline
&&&\hfil{\rm{\footnotesize $3$}}\hfill\\
\hline
{\rm{\footnotesize $2$}}\hfil&\hfil{\rm{\footnotesize $p=4+1=5$}}\hfill&\hfil{\rm{\footnotesize ($\square$)}}\hfil&\hfill\\
\hline
{\rm{\footnotesize $3$}}\hfil&\hfil{\rm{\footnotesize $p=9+1=10$}}\hfil&&\hfil{\rm{\footnotesize $7,\, 11,\, 13$}}\hfill\\
\hline
{\rm{\footnotesize $4$}}\hfil&\hfil{\rm{\footnotesize $p=16+1=17$}}\hfil&\hfil{\rm{\footnotesize ($\square$)}}\hfil&\hfill\\
\hline
{\rm{\footnotesize $5$}}\hfil&\hfil{\rm{\footnotesize $p=25+1=26$}}\hfil&&\hfil{\rm{\footnotesize $19,\, 23,\, 29,\, 31$}}\hfill\\
\hline
\hfil{\rm{\footnotesize $6$}}\hfil&\hfil{\rm{\footnotesize $p=36+1=37$}}\hfil&\hfil{\rm{\footnotesize ($\square$)}}\hfil&\hfill\\
\hline
{\rm{\footnotesize $7$}}\hfil&\hfil{\rm{\footnotesize $p=49+1=50$}}\hfil&&\hfil{\rm{\footnotesize $41,\, 43,\, 47,\, 53$}}\hfill\\
\hline
{\rm{\footnotesize $8$}}\hfil&\hfil{\rm{\footnotesize $p=64+1=65$}}\hfil&&\hfil{\rm{\footnotesize $59,\, 61,\, 67,\, 71$}}\hfill\\
\hline
{\rm{\footnotesize $9$}}\hfil&\hfil{\rm{\footnotesize $p=81+1=82$}}\hfill&&\hfil{\rm{\footnotesize $73,\, 79,\, 83,\, 89,\, 97$}}\hfil\\
\hline
{\rm{\footnotesize $10$}}\hfil&\hfil{\rm{\footnotesize $p=100+1=101$}}\hfil&\hfil{\rm{\footnotesize ($\square$)}}\hfil&\hfill\\
\hline
{\rm{\footnotesize $11$}}\hfil&\hfil{\rm{\footnotesize $p=121+1=122$}}\hfil&&\hfil{\rm{\footnotesize $103,\, 107,\, 109,\, 113,\, 127,\, 131,\, 137,\, 139$}}\hfill\\
\hline
\hfil{\rm{\footnotesize $12$}}\hfil&\hfil{\rm{\footnotesize $p=144+1=145$}}\hfil&&\hfil{\rm{\footnotesize $149,\, 151,\, 157,\, 163,\, 167,\, 173$}}\hfill\\
\hline
{\rm{\footnotesize $13$}}\hfil&\hfil{\rm{\footnotesize $p=169+1=170$}}\hfil&&\hfil{\rm{\footnotesize $179,\, 181,\, 191,\, 193$}}\hfill\\
\hline
{\rm{\footnotesize $14$}}\hfil&\hfil{\rm{\footnotesize $p=196+1=197$}}\hfil&\hfil{\rm{\footnotesize ($\square$)}}\hfil&\hfill\\
\hline
{\rm{\footnotesize $15$}}\hfil&\hfil{\rm{\footnotesize $p=225+1=226$}}\hfil&&\hfil{\rm{\footnotesize $199,\, 211\, 223,\, 227,\, 229,\, 233,\, 239,\, 241,\, 245$}}\hfill\\
\hline
{\rm{\footnotesize $16$}}\hfil&\hfil{\rm{\footnotesize $p=256+1=257$}}\hfil&\hfil{\rm{\footnotesize ($\square$)}}&\hfill\\
\hline
{\rm{\footnotesize $17$}}\hfil&\hfil{\rm{\footnotesize $p=289+1=290$}}\hfil&&\hfil{\rm{\footnotesize $263,\, 269,\, 271,\, 277,\, 281,\, 283,\, 293,\, 307,\, 311,\, 313$}}\hfill\\
\hline
{\rm{\footnotesize $18$}}\hfil&\hfil{\rm{\footnotesize $p=324+1=325$}}\hfil&&\hfil{\rm{\footnotesize $317,\, 331,\, 337,\, 347,\, 349,\, 353,\, 359$}}\hfill\\
\hline
{\rm{\footnotesize $19$}}\hfil&\hfil{\rm{\footnotesize $p=361+1=362$}}\hfil&&\hfil{\rm{\footnotesize $367,\, 373,\, 379,\, 383,\, 389,\, 397$}}\hfill\\
\hline
{\rm{\footnotesize $20$}}\hfil&\hfil{\rm{\footnotesize $p=400+1=401$}}\hfil&\hfil{\rm{\footnotesize ($\square$)}}&\hfill\\
\hline
{\rm{\footnotesize $21$}}\hfil&\hfil{\rm{\footnotesize $p=441+1=442$}}\hfil&&\hfil{\rm{\footnotesize $409,\, 419,\, 421,\, 431,\, 433,\, 439,\, 443,\, 449,\, 457,\, 461$}}\hfill\\
\hline
{\rm{\footnotesize $22$}}\hfil&\hfil{\rm{\footnotesize $p=484+1=485$}}\hfil&&\hfil{\rm{\footnotesize $463,\, 467,\, 479,\, 487,\, 491,\, 499,\, 503,\, 509,\, 521,\, 523$}}\hfill\\
\hline
{\rm{\footnotesize $23$}}\hfil&\hfil{\rm{\footnotesize $p=529+1=530$}}\hfil&&\hfil{\rm{\footnotesize $541,\, 547,\, 557,\, 563,\, 569,\, 571$}}\hfill\\
\hline
{\rm{\footnotesize $24$}}\hfil&\hfil{\rm{\footnotesize $p=576+1=577$}}\hfil&\hfil{\rm{\footnotesize ($\square$)}}&\hfill\\
\hline
{\rm{\footnotesize $25$}}\hfil&\hfil{\rm{\footnotesize $p=625+1=626$}}\hfil&&\hfil{\rm{\footnotesize $587,\, 593,\, 599,\, 601,\, 607,\, 613,\, 617,\, 619,\, 631,\, 641,\, 643,\, 647,\, 653,\, 659,\, 661,\, 673$}}\hfill\\
\hline
{\rm{\footnotesize $26$}}\hfil&\hfil{\rm{\footnotesize $p=676+1=677$}}\hfil&\hfil{\rm{\footnotesize ($\square$)}}&\hfill\\
\hline
{\rm{\footnotesize $27$}}\hfil&\hfil{\rm{\footnotesize $p=729+1=730$}}\hfil&&\hfil{\rm{\footnotesize $683,\, 691,\, 701,\, 709,\, 719,\, 727,\, 733,\, 739,\, 739,\, 743,\, 751$}}\hfill\\
\hline
{\rm{\footnotesize $28$}}\hfil&\hfil{\rm{\footnotesize $p=784+1=785$}}\hfil&&\hfil{\rm{\footnotesize $757,\, 761,\, 769,\, 773,\, 787,\, 797,\, 809,\, 811,\, 821,\, 823$}}\hfill\\
\hline
{\rm{\footnotesize $29$}}\hfil&\hfil{\rm{\footnotesize $p=841+1=842$}}\hfil&&\hfil{\rm{\footnotesize $827,\, 829,\, 839,\, 853,\, 857,\, 859,\, 863,\, 877,\, 881,\, 883$}}\hfill\\
\hline
{\rm{\footnotesize $30$}}\hfil&\hfil{\rm{\footnotesize $p=900+1=901$}}\hfil&&\hfil{\rm{\footnotesize $887,\, 907,\, 911,\, 919,\, 929,\, 937,\, 941,\, 947,\, 953,\, 967,\, 971$}}\hfill\\
\hline
{\rm{\footnotesize $31$}}\hfil&\hfil{\rm{\footnotesize $p=961+1=962$}}\hfil&&\hfil{\rm{\footnotesize $977,\, 983,\, 991,\, 997,\, 1009,\, 1013,\, 1019,\, 1021,\, 1031,\, 1033$}}\hfill\\
\hline
{\rm{\footnotesize $32$}}\hfil&\hfil{\rm{\footnotesize $p=1024+1=1025$}}\hfil&&\hfil{\rm{\footnotesize $1039,\, 1049,\, 1051,\, 1061,\, 1063,\, 1069,\, 1087,\, 1091,\, 1093,\, 1097$}}\hfill\\
\hline
{\rm{\footnotesize $33$}}\hfil&\hfil{\rm{\footnotesize $p=1089+1=1090$}}\hfil&&\hfil{\rm{\footnotesize $1103,\, 1109,\, 1117,\, 1123,\, 1129,\, 1151,\, 1153,\, 1163,\, 1171,\, 1181$}}\hfill\\
\hline
{\rm{\footnotesize $34$}}\hfil&\hfil{\rm{\footnotesize $p=1156+1=1157$}}\hfil&&\hfil{\rm{\footnotesize $1187,\, 1193,\, 1201,\, 1213,\, 1217,\, 1223,\, 1229,\, 1231,\, 1237,\, 1249$}}\hfill\\
\hline
{\rm{\footnotesize $35$}}\hfil&\hfil{\rm{\footnotesize $p=1225+1=1226$}}\hfil&&\hfil{\rm{\footnotesize $1259,\, 1277,\, 1279,\, 1283,\, 1289,\, 1291$}}\hfill\\
\hline
{\rm{\footnotesize $36$}}\hfil&\hfil{\rm{\footnotesize $p=1296+1=1297$}}\hfil&\hfil{\rm{\footnotesize ($\square$)}}&\hfill\\
\hline
{\rm{\footnotesize $37$}}\hfil&\hfil{\rm{\footnotesize $p=1369+1=1370$}}\hfil&&\hfil{\rm{\footnotesize $
1301,\, 1303,\, 1307,\, 1319,\, 1321,\, 1327,\, 1361,\, 1367,\, 1373,\, 1381,\, 1399,\, 1409,\, 1423,\, 1427,\,
$}}\hfill\\
\hline
{\rm{\footnotesize $38$}}\hfil&\hfil{\rm{\footnotesize $p=1444+1=1445$}}\hfil&&\hfil{\rm{\footnotesize $
1429,\, 1433,\, 1439,\, 1447,\, 1451,\, 1453,\, 1459,\, 1471,\, 1481,\, 1483,\, 1487,\, 1489,\, 1493,\, 1499
$}}\hfill\\
\hline
{\rm{\footnotesize $39$}}\hfil&\hfil{\rm{\footnotesize $p=1521+1=1522$}}\hfil&&\hfil{\rm{\footnotesize $
1511,\, 1523,\, 1531,\, 1543,\, 1549,\, 1553,\, 1559,\, 1567,\, 1571,\, 1579,\, 1583,\, 1597
$}}\hfill\\
\hline
{\rm{\footnotesize $40$}}\hfil&\hfil{\rm{\footnotesize $p=1600+1=1601$}}\hfil&\hfil{\rm{\footnotesize ($\square$)}}&\hfill\\
\hline
{\rm{\footnotesize $42$}}\hfil&\hfil{\rm{\footnotesize $p=1764+1=1765$}}\hfil&&\hfill\\
\hline
\hline{\rm{\footnotesize $44$}}\hfil&\hfil{\rm{\footnotesize $p=1936+1=1937$}}\hfil&&\hfill\\
\hline{\rm{\footnotesize $46$}}\hfil&\hfil{\rm{\footnotesize $p=2116+1=2117$}}\hfil&&\hfill\\
\hline{\rm{\footnotesize $48$}}\hfil&\hfil{\rm{\footnotesize $p=2304+1=2305$}}\hfil&&\hfill\\
\hline{\rm{\footnotesize $50$}}\hfil&\hfil{\rm{\footnotesize $p=2500+1=2501$}}\hfil&&\hfill\\
\hline{\rm{\footnotesize $52$}}\hfil&\hfil{\rm{\footnotesize $p=2704+1=2705$}}\hfil&&\hfill\\
\hline{\rm{\footnotesize $54$}}\hfil&\hfil{\rm{\footnotesize $p=2916+1=2917$}}\hfil&\hfil{\rm{\footnotesize ($\square$)}}&\hfill\\
\hline{\rm{\footnotesize $56$}}\hfil&\hfil{\rm{\footnotesize $p=3136+1=3137$}}\hfil&\hfil{\rm{\footnotesize ($\square$)}}&\hfill\\
\hline{\rm{\footnotesize $58$}}\hfil&\hfil{\rm{\footnotesize $p=3364+1=3365$}}\hfil&&\hfill\\
\hline{\rm{\footnotesize $60$}}\hfil&\hfil{\rm{\footnotesize $p=3600+1=3601$}}\hfil&&\hfill\\
\hline
\multicolumn{4}{l}{\rm{\footnotesize For $n\ge 40$ are reported only calculations for even numbers $n=2m$, according to Lemma \ref{necessity-even}.}}\\
\multicolumn{4}{l}{\rm{\footnotesize Furthermore for $n\ge 40$ are omitted primes non-parabolic prime.}}\\
\end{tabular}$}
\end{table}

\begin{example}[Some examples of $n$-ghost right triangles]\label{examples-ghost-right-triangles}
It is easy to find natural numbers $n\in\mathbb{N}$ satisfying the condition $n^2+1=p\in P$. In Tab. \ref{tab-examples-ghost-right-triangles} are reported the $n$-ghost right-triangles with $1\le n\le 60$. There are reported also, for $1\le n\le 40$, the primes between consecutive primes $p$, such that $p-1$ is a perfect square. In order to represent right-triangles in $\mathbb{N}^2$, we can adopt the assumption to draw their hypotenuses by full-lines for Pythagorean triangles and dot-lines for $n$-ghost right-triangles. In {\em(\ref{pythagorean-and-ghost-right-triangles-in-n-2})} are represented the Pythagorean triangle $(3,4,5)$ and $(6,8,10)$, and the $n$-ghost right-triangles with $n=6$ and $n=10$. The hypotenuses of these last right-triangles have length respectively $d((0,6),(1,0))=\sqrt{37}\notin\mathbb{N}$, $d((0,10),(1,0))=\sqrt{101}\notin\mathbb{N}$, since $37$ and $101$ are prime. Instead for the hypotenuses of the reported Pythagorean triangles one has respectively $d((0,8),(6,0))=\sqrt{100}=10\in\mathbb{N}$ and $d((0,4),(3,0))=\sqrt{25}=5\in\mathbb{N}$. See figure in the left in {\em(\ref{pythagorean-and-ghost-right-triangles-in-n-2})} where full-line segments represent integer numbers (i.e., are metric-regular lines) and dot-line segments do not represent integer numbers (i.e., are ghost lines).

Let us emphasize that the set $P_{\square}$ can be identified with a discrete subset of the branch of parabola $\Gamma$, of equation $y=x^2+1$ in $\mathbb{R}^2$, $x\ge 0$, when in the positive $x$-axis there is embedded $\mathbb{N}\subset\mathbb{R}$ and in the positive $y$-axis there is embedded the set of primes $P\subset\mathbb{R}$. Then the points belonging to $P_{\square}\subset \Gamma$, are identified as intersection of {\em $P$-horizontal lines} (i.e., horizontal lines starting from the primes on the $y$-axis) that meet {\em $\mathbb{N}$-vertical lines}, (i.e.,vertical lines starting from the integers on the $x$-axis), on the the above parabola $\Gamma$. (See figure on the right in {\em(\ref{pythagorean-and-ghost-right-triangles-in-n-2})}. $P_{\square}$ is identified with the set $\Gamma_{\square}\subset\Gamma$, partially reported therein by some $\square$-points on $\Gamma$.) For this reason we shall call the set $P_{\square}$ also the {\em set of parabolic primes}, according we call {\em parabolic prime} a prime $p\in P_{\square}$.
\end{example}

\begin{proposition}[$n$-ghost right-triangles and Pythagorean triangles]\label{pythagorean-ghost-right-triangles}
$\bullet$\hskip 2pt In a $n$-ghost right-triangle the legs are metric-regular lines. Instead the hypotenus is a ghost line.

$\bullet$\hskip 2pt A $n$-ghost right-triangle cannot be a Pythagorean triangle and neither a primitive Pythagorean triangle.
\end{proposition}

\begin{proof}
$\bullet$\hskip 2pt In fact $\overline{AC}=n,\, \overline{CB}=1\in \mathbb{N}$. Instead $\overline{AB}=\sqrt{p}\not\in\mathbb{N}$.

$\bullet$\hskip 2pt A $n$-ghost right-triangle $ABC$ has the length of the hypotenus $\overline{AB}=q$ such that $q^2=p\in P$. But this is impossible since $p$ is prime. A $n$-ghost right-triangle cannot be neither a primitive Pythagorean triangle. In fact $n$ cannot divide $p$.
\end{proof}
\begin{equation}\label{pythagorean-and-ghost-right-triangles-in-n-2}
  \scalebox{0.6}{$ \xymatrix{&&&&&&&\\
   *+[o][F]{10}\ar@{.}[ddddddddddr]\ar[u]^{b}\ar@{-}[d]\ar@{-}[r]&\circ\ar@{-}[d]\ar@{-}[r]&\circ\ar@{-}[d]\ar@{-}[r]&
   \circ\ar@{-}[d]\ar@{-}[r]&\circ\ar@{-}[d]\ar@{-}[r]&\circ\ar@{-}[d]\ar@{-}[r]&
   \circ\ar@{-}[d]&\\
    \circ\ar@{-}[d]\ar@{-}[r]&\circ\ar@{-}[d]\ar@{-}[r]&\circ\ar@{-}[d]\ar@{-}[r]&
   \circ\ar@{-}[d]\ar@{-}[r]&\circ\ar@{-}[d]\ar@{-}[r]&\circ\ar@{-}[d]\ar@{-}[r]&
   \circ\ar@{-}[d]&\\
 *+[o][F]{8}\ar@{-}[ddddddddrrrrrr]\ar@{-}[d]\ar@{-}[r]&\circ\ar@{-}[d]\ar@{-}[r]&\circ\ar@{-}[d]\ar@{-}[r]&
   \circ\ar@{-}[d]\ar@{-}[r]&\circ\ar@{-}[d]\ar@{-}[r]&\circ\ar@{-}[d]\ar@{-}[r]&
   \circ\ar@{-}[d]&\\
   \circ\ar@{-}[d]\ar@{-}[r]&\circ\ar@{-}[d]\ar@{-}[r]&\circ\ar@{-}[d]\ar@{-}[r]&
   \circ\ar@{-}[d]\ar@{-}[r]&\circ\ar@{-}[d]\ar@{-}[r]&\circ\ar@{-}[d]\ar@{-}[r]&
   \circ\ar@{-}[d]&\\
   *+[o][F]{6}\ar@{.}[ddddddr]\ar@{-}[d]\ar@{-}[r]&\circ\ar@{-}[d]\ar@{-}[r]&\circ\ar@{-}[d]\ar@{-}[r]&
   \circ\ar@{-}[d]\ar@{-}[r]&\circ\ar@{-}[d]\ar@{-}[r]&\circ\ar@{-}[d]\ar@{-}[r]&
   \circ\ar@{-}[d]&\\
   \circ\ar@{-}[d]\ar@{-}[r]&\circ\ar@{-}[d]\ar@{-}[r]&\circ\ar@{-}[d]\ar@{-}[r]&
   \circ\ar@{-}[d]\ar@{-}[r]&\circ\ar@{-}[d]\ar@{-}[r]&\circ\ar@{-}[d]\ar@{-}[r]&
   \circ\ar@{-}[d]&\\
   *+[o][F]{4}\ar@{.}[ddddr]\ar@{-}[ddddrrr]\ar@{-}[d]\ar@{-}[r]&\circ\ar@{-}[d]\ar@{-}[r]&\circ\ar@{-}[d]\ar@{-}[r]&
   \circ\ar@{-}[d]\ar@{-}[r]&\circ\ar@{-}[d]\ar@{-}[r]&\circ\ar@{-}[d]\ar@{-}[r]&
   \circ\ar@{-}[d]&\\
   \circ\ar@{-}[d]\ar@{-}[r]&\circ\ar@{-}[d]\ar@{-}[r]&\circ\ar@{-}[d]\ar@{-}[r]&
   \circ\ar@{-}[d]\ar@{-}[r]&\circ\ar@{-}[d]\ar@{-}[r]&\circ\ar@{-}[d]\ar@{-}[r]&
   \circ\ar@{-}[d]&\\
   *+[o][F]{2}\ar@{.}[ddr]\ar@{-}[d]\ar@{-}[r]&\circ\ar@{-}[d]\ar@{-}[r]&\circ\ar@{-}[d]\ar@{-}[r]&
   \circ\ar@{-}[d]\ar@{-}[r]&\circ\ar@{-}[d]\ar@{-}[r]&\circ\ar@{-}[d]\ar@{-}[r]&
   \circ\ar@{-}[d]&\\
   *+[o][F]{1}\ar@{.}[dr]\ar@{-}[d]\ar@{-}[r]&\circ\ar@{-}[d]\ar@{-}[r]&\circ\ar@{-}[d]\ar@{-}[r]&
   \circ\ar@{-}[d]\ar@{-}[r]&\circ\ar@{-}[d]\ar@{-}[r]&\circ\ar@{-}[d]\ar@{-}[r]&
   \circ\ar@{-}[d]&\\
   *+[o][F]{0}\ar@{-}[r]&*+[o][F]{1}\ar@{-}[r]&\circ\ar@{-}[r]&
   *+[o][F]{3}\ar@{-}[r]&\circ\ar@{-}[r]&\circ\ar@{-}[r]&*+[o][F]{6}\ar[r]_{a}&\\}$}\,
\hskip 5pt\scalebox{0.3}{$ \xymatrix{&&&&&&&&&&&&\\
   \scalebox{2.0}{$P$}\ar@{^{(}->}[d]&*+[o][F]{\scalebox{2.0}{$101$}}\ar[u]^{\scalebox{2.0}{$y$}}\ar@{-}[d]\ar@{-}[r]&
   \circ\ar@{-}[d]\ar@{-}[r]&\circ\ar@{-}[d]\ar@{-}[r]&
   \circ\ar@{-}[d]\ar@{-}[r]&\circ\ar@{-}[d]\ar@{-}[r]&\circ\ar@{-}[d]\ar@{-}[r]&\circ\ar@{-}[d]\ar@{-}[r]&
   \circ\ar@{-}[d]\ar@{-}[r]&\circ\ar@{-}[d]\ar@{-}[r]&\circ\ar@{-}[r]\ar@{-}[d]&\scalebox{2.0}{$\square$}\ar@{-}[d]&&\\
    \scalebox{2.0}{$\mathbb{R}$}& *+[o][F]{97}\ar@{-}[d]\ar@{-}[r]&\circ\ar@{-}[d]\ar@{-}[r]&\circ\ar@{-}[d]\ar@{-}[r]&
   \circ\ar@{-}[d]\ar@{-}[r]&\circ\ar@{-}[d]\ar@{-}[r]&\circ\ar@{-}[d]\ar@{-}[r]&
   \circ\ar@{-}[d]\ar@{-}[r]&\circ\ar@{-}[d]\ar@{-}[r]&\circ\ar@{-}[d]\ar@{-}[r]&
   \circ\ar@{-}[d]\ar@{-}[r]&\circ\ar@{-}[d]&\\
  & *+[o][F]{89}\ar@{-}[d]\ar@{-}[r]&\circ\ar@{-}[d]\ar@{-}[r]&\circ\ar@{-}[d]\ar@{-}[r]&
   \circ\ar@{-}[d]\ar@{-}[r]&\circ\ar@{-}[d]\ar@{-}[r]&\circ\ar@{-}[d]\ar@{-}[r]&
   \circ\ar@{-}[d]\ar@{-}[r]&\circ\ar@{-}[d]\ar@{-}[r]&\circ\ar@{-}[d]\ar@{-}[r]&
   \circ\ar@{-}[d]\ar@{-}[r]&\circ\ar@{-}[d]&\\
   & *+[o][F]{83}\ar@{-}[d]\ar@{-}[r]&\circ\ar@{-}[d]\ar@{-}[r]&\circ\ar@{-}[d]\ar@{-}[r]&
   \circ\ar@{-}[d]\ar@{-}[r]&\circ\ar@{-}[d]\ar@{-}[r]&\circ\ar@{-}[d]\ar@{-}[r]&
   \circ\ar@{-}[d]\ar@{-}[r]&\circ\ar@{-}[d]\ar@{-}[r]&\circ\ar@{-}[d]\ar@{-}[r]&
   \circ\ar@{-}[d]\ar@{-}[r]&\circ\ar@{-}[d]&\\
   & *+[o][F]{79}\ar@{-}[d]\ar@{-}[r]&\circ\ar@{-}[d]\ar@{-}[r]&\circ\ar@{-}[d]\ar@{-}[r]&
   \circ\ar@{-}[d]\ar@{-}[r]&\circ\ar@{-}[d]\ar@{-}[r]&\circ\ar@{-}[d]\ar@{-}[r]&
   \circ\ar@{-}[d]\ar@{-}[r]&\circ\ar@{-}[d]\ar@{-}[r]&\circ\ar@{-}[d]\ar@{-}[r]&
   \circ\ar@{-}[d]\ar@{-}[r]&\circ\ar@{-}[d]&\\
   & *+[o][F]{73}\ar@{-}[d]\ar@{-}[r]&\circ\ar@{-}[d]\ar@{-}[r]&\circ\ar@{-}[d]\ar@{-}[r]&
   \circ\ar@{-}[d]\ar@{-}[r]&\circ\ar@{-}[d]\ar@{-}[r]&\circ\ar@{-}[d]\ar@{-}[r]&
   \circ\ar@{-}[d]\ar@{-}[r]&\circ\ar@{-}[d]\ar@{-}[r]&\circ\ar@{-}[d]\ar@{-}[r]&
   \circ\ar@{-}[d]\ar@{-}[r]&\circ\ar@{-}[d]&\\
  & *+[o][F]{71}\ar@{-}[d]\ar@{-}[r]&\circ\ar@{-}[d]\ar@{-}[r]&\circ\ar@{-}[d]\ar@{-}[r]&
   \circ\ar@{-}[d]\ar@{-}[r]&\circ\ar@{-}[d]\ar@{-}[r]&\circ\ar@{-}[d]\ar@{-}[r]&
   \circ\ar@{-}[d]\ar@{-}[r]&\circ\ar@{-}[d]\ar@{-}[r]&\circ\ar@{-}[d]\ar@{-}[r]&
   \circ\ar@{-}[d]\ar@{-}[r]&\circ\ar@{-}[d]&\\
   & *+[o][F]{67}\ar@{-}[d]\ar@{-}[r]&\circ\ar@{-}[d]\ar@{-}[r]&\circ\ar@{-}[d]\ar@{-}[r]&
   \circ\ar@{-}[d]\ar@{-}[r]&\circ\ar@{-}[d]\ar@{-}[r]&\circ\ar@{-}[d]\ar@{-}[r]&
   \circ\ar@{-}[d]\ar@{-}[r]&\circ\ar@{-}[d]\ar@{-}[r]&\circ\ar@{-}[d]\ar@{-}[r]&
   \circ\ar@{-}[d]\ar@{-}[r]&\circ\ar@{-}[d]&\\
   & *+[o][F]{61}\ar@{-}[d]\ar@{-}[r]&\circ\ar@{-}[d]\ar@{-}[r]&\circ\ar@{-}[d]\ar@{-}[r]&
   \circ\ar@{-}[d]\ar@{-}[r]&\circ\ar@{-}[d]\ar@{-}[r]&\circ\ar@{-}[d]\ar@{-}[r]&
   \circ\ar@{-}[d]\ar@{-}[r]&\circ\ar@{-}[d]\ar@{-}[r]&\circ\ar@{-}[d]\ar@{-}[r]&
   \circ\ar@{-}[d]\ar@{-}[r]&\circ\ar@{-}[d]&\\
  & *+[o][F]{59}\ar@{-}[d]\ar@{-}[r]&\circ\ar@{-}[d]\ar@{-}[r]&\circ\ar@{-}[d]\ar@{-}[r]&
   \circ\ar@{-}[d]\ar@{-}[r]&\circ\ar@{-}[d]\ar@{-}[r]&\circ\ar@{-}[d]\ar@{-}[r]&
   \circ\ar@{-}[d]\ar@{-}[r]&\circ\ar@{-}[d]\ar@{-}[r]&\circ\ar@{-}[d]\ar@{-}[r]&
   \circ\ar@{-}[d]\ar@{-}[r]&\circ\ar@{-}[d]&\\
   & *+[o][F]{53}\ar@{-}[d]\ar@{-}[r]&\circ\ar@{-}[d]\ar@{-}[r]&\circ\ar@{-}[d]\ar@{-}[r]&
   \circ\ar@{-}[d]\ar@{-}[r]&\circ\ar@{-}[d]\ar@{-}[r]&\circ\ar@{-}[d]\ar@{-}[r]&
   \circ\ar@{-}[d]\ar@{-}[r]&\circ\ar@{-}[d]\ar@{-}[r]&\circ\ar@{-}[d]\ar@{-}[r]&
   \circ\ar@{-}[d]\ar@{-}[r]&\circ\ar@{-}[d]&\\
   & *+[o][F]{47}\ar@{-}[d]\ar@{-}[r]&\circ\ar@{-}[d]\ar@{-}[r]&\circ\ar@{-}[d]\ar@{-}[r]&
   \circ\ar@{-}[d]\ar@{-}[r]&\circ\ar@{-}[d]\ar@{-}[r]&\circ\ar@{-}[d]\ar@{-}[r]&
   \circ\ar@{-}[d]\ar@{-}[r]&\circ\ar@{-}[d]\ar@{-}[r]&\circ\ar@{-}[d]\ar@{-}[r]&
   \circ\ar@{-}[d]\ar@{-}[r]&\circ\ar@{-}[d]&\\
   & *+[o][F]{43}\ar@{-}[d]\ar@{-}[r]&\circ\ar@{-}[d]\ar@{-}[r]&\circ\ar@{-}[d]\ar@{-}[r]&
   \circ\ar@{-}[d]\ar@{-}[r]&\circ\ar@{-}[d]\ar@{-}[r]&\circ\ar@{-}[d]\ar@{-}[r]&
   \circ\ar@{-}[d]\ar@{-}[r]&\circ\ar@{-}[d]\ar@{-}[r]&\circ\ar@{-}[d]\ar@{-}[r]&
   \circ\ar@{-}[d]\ar@{-}[r]&\circ\ar@{-}[d]&\\
   & *+[o][F]{41}\ar@{-}[d]\ar@{-}[r]&\circ\ar@{-}[d]\ar@{-}[r]&\circ\ar@{-}[d]\ar@{-}[r]&
   \circ\ar@{-}[d]\ar@{-}[r]&\circ\ar@{-}[d]\ar@{-}[r]&\circ\ar@{-}[d]\ar@{-}[r]&
   \circ\ar@{-}[d]\ar@{-}[r]&\circ\ar@{-}[d]\ar@{-}[r]&\circ\ar@{-}[d]\ar@{-}[r]&
   \circ\ar@{-}[d]\ar@{-}[r]&\circ\ar@{-}[d]&\\
   & *+[o][F]{\scalebox{2.0}{$37$}}\ar@{-}[d]\ar@{-}[r]&\circ\ar@{-}[d]\ar@{-}[r]&\circ\ar@{-}[d]\ar@{-}[r]&
   \circ\ar@{-}[d]\ar@{-}[r]&\circ\ar@{-}[d]\ar@{-}[r]&\circ\ar@{-}[d]\ar@{-}[r]&
   \scalebox{2.0}{$\square$}\ar@{.}[uuuuuuuuuuuuuurrrr]\ar@{-}[d]\ar@{-}[r]&
   \circ\ar@{-}[d]\ar@{-}[r]&\circ\ar@{-}[d]\ar@{-}[r]&
   \circ\ar@{-}[d]\ar@{-}[r]&\circ\ar@{-}[d]&\\
   & *+[o][F]{31}\ar@{-}[d]\ar@{-}[r]&\circ\ar@{-}[d]\ar@{-}[r]&\circ\ar@{-}[d]\ar@{-}[r]&
   \circ\ar@{-}[d]\ar@{-}[r]&\circ\ar@{-}[d]\ar@{-}[r]&\circ\ar@{-}[d]\ar@{-}[r]&
   \circ\ar@{-}[d]\ar@{-}[r]&\circ\ar@{-}[d]\ar@{-}[r]&\circ\ar@{-}[d]\ar@{-}[r]&
   \circ\ar@{-}[d]\ar@{-}[r]&\circ\ar@{-}[d]&\\
   & *+[o][F]{29}\ar@{-}[d]\ar@{-}[r]&\circ\ar@{-}[d]\ar@{-}[r]&\circ\ar@{-}[d]\ar@{-}[r]&
   \circ\ar@{-}[d]\ar@{-}[r]&\circ\ar@{-}[d]\ar@{-}[r]&\circ\ar@{-}[d]\ar@{-}[r]&
   \circ\ar@{-}[d]\ar@{-}[r]&\circ\ar@{-}[d]\ar@{-}[r]&\circ\ar@{-}[d]\ar@{-}[r]&
   \circ\ar@{-}[d]\ar@{-}[r]&\circ\ar@{-}[d]&\\
   & *+[o][F]{23}\ar@{-}[d]\ar@{-}[r]&\circ\ar@{-}[d]\ar@{-}[r]&\circ\ar@{-}[d]\ar@{-}[r]&
   \circ\ar@{-}[d]\ar@{-}[r]&\circ\ar@{-}[d]\ar@{-}[r]&\circ\ar@{-}[d]\ar@{-}[r]&
   \circ\ar@{-}[d]\ar@{-}[r]&\circ\ar@{-}[d]\ar@{-}[r]&\circ\ar@{-}[d]\ar@{-}[r]&
   \circ\ar@{-}[d]\ar@{-}[r]&\circ\ar@{-}[d]&\\
   & *+[o][F]{19}\ar@{-}[d]\ar@{-}[r]&\circ\ar@{-}[d]\ar@{-}[r]&\circ\ar@{-}[d]\ar@{-}[r]&
   \circ\ar@{-}[d]\ar@{-}[r]&\circ\ar@{-}[d]\ar@{-}[r]&\circ\ar@{-}[d]\ar@{-}[r]&
   \circ\ar@{-}[d]\ar@{-}[r]&\circ\ar@{-}[d]\ar@{-}[r]&\circ\ar@{-}[d]\ar@{-}[r]&
   \circ\ar@{-}[d]\ar@{-}[r]&\circ\ar@{-}[d]&\\
   & *+[o][F]{\scalebox{2.0}{$17$}}\ar@{-}[d]\ar@{-}[r]&\circ\ar@{-}[d]\ar@{-}[r]&\circ\ar@{-}[d]\ar@{-}[r]&
   \circ\ar@{-}[d]\ar@{-}[r]&\scalebox{2.0}{$\square$}\ar@{.}[uuuuurr]\ar@{-}[d]\ar@{-}[r]&\circ\ar@{-}[d]\ar@{-}[r]&
   \circ\ar@{-}[d]\ar@{-}[r]&\circ\ar@{-}[d]\ar@{-}[r]&\circ\ar@{-}[d]\ar@{-}[r]&
   \circ\ar@{-}[d]\ar@{-}[r]&\circ\ar@{-}[d]&\\
   & *+[o][F]{13}\ar@{-}[d]\ar@{-}[r]&\circ\ar@{-}[d]\ar@{-}[r]&\circ\ar@{-}[d]\ar@{-}[r]&
   \circ\ar@{-}[d]\ar@{-}[r]&\circ\ar@{-}[d]\ar@{-}[r]&\circ\ar@{-}[d]\ar@{-}[r]&
   \circ\ar@{-}[d]\ar@{-}[r]&\circ\ar@{-}[d]\ar@{-}[r]&\circ\ar@{-}[d]\ar@{-}[r]&
   \circ\ar@{-}[d]\ar@{-}[r]&\circ\ar@{-}[d]&\\
   & *+[o][F]{11}\ar@{-}[d]\ar@{-}[r]&\circ\ar@{-}[d]\ar@{-}[r]&\circ\ar@{-}[d]\ar@{-}[r]&
   \circ\ar@{-}[d]\ar@{-}[r]&\circ\ar@{-}[d]\ar@{-}[r]&\circ\ar@{-}[d]\ar@{-}[r]&
   \circ\ar@{-}[d]\ar@{-}[r]&\circ\ar@{-}[d]\ar@{-}[r]&\circ\ar@{-}[d]\ar@{-}[r]&
   \circ\ar@{-}[d]\ar@{-}[r]&\circ\ar@{-}[d]&\\
   & *+[o][F]{7}\ar@{-}[d]\ar@{-}[r]&\circ\ar@{-}[d]\ar@{-}[r]&\circ\ar@{-}[d]\ar@{-}[r]&
   \circ\ar@{-}[d]\ar@{-}[r]&\circ\ar@{-}[d]\ar@{-}[r]&\circ\ar@{-}[d]\ar@{-}[r]&
   \circ\ar@{-}[d]\ar@{-}[r]&\circ\ar@{-}[d]\ar@{-}[r]&\circ\ar@{-}[d]\ar@{-}[r]&
   \circ\ar@{-}[d]\ar@{-}[r]&\circ\ar@{-}[d]&\\
   & *+[o][F]{\scalebox{2.0}{$5$}}\ar@{-}[d]\ar@{-}[r]&\circ\ar@{-}[d]\ar@{-}[r]&
  \scalebox{2.0}{$\square$}\ar@{.}[uuuurr]\ar@{-}[d]\ar@{-}[r]&
   \circ\ar@{-}[d]\ar@{-}[r]&\circ\ar@{-}[d]\ar@{-}[r]&\circ\ar@{-}[d]\ar@{-}[r]&
   \circ\ar@{-}[d]\ar@{-}[r]&\circ\ar@{-}[d]\ar@{-}[r]&\circ\ar@{-}[d]\ar@{-}[r]&
   \circ\ar@{-}[d]\ar@{-}[r]&\circ\ar@{-}[d]&\\
   & *+[o][F]{3}\ar@{-}[d]\ar@{-}[r]&\circ\ar@{-}[d]\ar@{-}[r]&\circ\ar@{-}[d]\ar@{-}[r]&
   \circ\ar@{-}[d]\ar@{-}[r]&\circ\ar@{-}[d]\ar@{-}[r]&\circ\ar@{-}[d]\ar@{-}[r]&
   \circ\ar@{-}[d]\ar@{-}[r]&\circ\ar@{-}[d]\ar@{-}[r]&\circ\ar@{-}[d]\ar@{-}[r]&
   \circ\ar@{-}[d]\ar@{-}[r]&\circ\ar@{-}[d]&\\
   & *+[o][F]{\scalebox{2.0}{$2$}}\ar@{-}[d]\ar@{-}[r]&\scalebox{2.0}{$\square$}\ar@{.}[uur]\ar@{-}[d]\ar@{-}[r]&
   \circ\ar@{-}[d]\ar@{-}[r]&
   \circ\ar@{-}[d]\ar@{-}[r]&\circ\ar@{-}[d]\ar@{-}[r]&\circ\ar@{-}[d]\ar@{-}[r]&
   \circ\ar@{-}[d]\ar@{-}[r]&\circ\ar@{-}[d]\ar@{-}[r]&\circ\ar@{-}[d]\ar@{-}[r]&
   \circ\ar@{-}[d]\ar@{-}[r]&\circ\ar@{-}[d]&\\
   & *+[o][F]{1}\ar@/_/[ur]\ar@{-}[d]\ar@{-}[r]&\circ\ar@{-}[d]\ar@{-}[r]&\circ\ar@{-}[d]\ar@{-}[r]&
   \circ\ar@{-}[d]\ar@{-}[r]&\circ\ar@{-}[d]\ar@{-}[r]&\circ\ar@{-}[d]\ar@{-}[r]&
   \circ\ar@{-}[d]\ar@{-}[r]&\circ\ar@{-}[d]\ar@{-}[r]&\circ\ar@{-}[d]\ar@{-}[r]&
   \circ\ar@{-}[d]\ar@{-}[r]&\circ\ar@{-}[d]&\\
   &*+[o][F]{0}\ar@{-}[r]&*+[o][F]{\scalebox{2.0}{$1$}}\ar@{-}[r]&*+[o][F]{\scalebox{2.0}{$2$}}\ar@{-}[r]&
   *+[o][F]{3}\ar@{-}[r]&*+[o][F]{\scalebox{2.0}{$4$}}\ar@{-}[r]&*+[o][F]{5}\ar@{-}[r]&*+[o][F]{\scalebox{2.0}{$6$}}\ar@{-}[r]
   &*+[o][F]{7}\ar@{-}[r]&*+[o][F]{8}\ar@{-}[r]&*+[o][F]{9}\ar@{-}[r]&*+[o][F]{\scalebox{2.0}{$10$}}\ar[r]_{\scalebox{2.0}{$x$}}&\\
   &&&&&&&&&&\scalebox{2.0}{$\mathbb{N}$}\ar@{^{(}->}[r]& \scalebox{2.0}{$\mathbb{R}$}&\\}$}
 \end{equation}

\begin{remark}
The justification of the Definition \ref{ghost-line}, Definition \ref{ghost-triangles} and Definition \ref{n-ghost-right-triangles} is to ascribe to the fact that for lines, legs and hypotenuses, that are not on the lattice $\mathbb{N}^2$, is not assured that their lengths are integers. Furthermore, when these lines are not on the lattice, the condition that should be metric-regular lines, means that they can be isometrically deformed on corresponding lines on the lattice. Therefore, ghost lines cannot be considered isometrically equivalent to lines on the lattice. For example by looking to the examples reported in {\em(\ref{pythagorean-and-ghost-right-triangles-in-n-2})}, we can see that the $6$-ghost right triangle has hypotenuse of length $\sqrt{37}\approx 6.08276253$. Thus this hypotenuse cannot be isometrically deformed on the lattice, e.g., on the legs, since their union has length $6+1=7$, namely an integer. Similarly happens for the $10$-ghost right triangle where the hypotenuse has length $\sqrt{101}\approx 10.04987562$. The union of its legs has length $10+1=11$. Therefore one cannot isometrically deform the hypotenuse on the union of the legs. Instead for the Pythagorean triple $(6,8,10)$ the hypotenuse has length $10$, hence can be isometrically deformed on the union of the legs that has length $6+8=14$. For example sending the hypotenuse on the leg of length $8$ and continuously the other part of length $2$, on the other leg (of length $4$). In such a way the line of the hypotenuse remains of length $10$, but lies exactly on the lattice. Similarly for the other Pythagorean triangle $(3,4,5)$, one has the hypotenuse of length $5$, that can be isometrically deformed on the union of its legs that have full length $4+3=7$.
\end{remark}

\begin{conj}[Conjecture on the infiniteness of $P_{\square}$]\label{the-infinite-set-ghost-rectangular}
There are infinite many parabolic primes $p$, or equivalently there are infinite $n$-ghost right-triangles.
\end{conj}

We get the following theorem that gives us a chance to answer in the affirmative to the Conjecture \ref{the-infinite-set-ghost-rectangular}.

\begin{theorem}[Euler-Riemann zeta function estimate of $P_{\square}$]\label{zeta-function-estimate-set-ghost-right-triangles}
One cannot exclude that $P_{\square}$ is infinite, since one has the limitation {\em(\ref{limitation-p-1-perfect-squares})}.
\begin{equation}\label{limitation-p-1-perfect-squares}
    1<\sum_{p\in P_{\square}}\frac{1}{p-1}\le\zeta(2)=\frac{\pi^2}{6}
\end{equation}
where $\zeta(\alpha)$ is the {\em Euler-Riemann zeta function}. We call {\em Euler-Riemann zeta function estimate of $P_{\square}$} the inequality {\em(\ref{limitation-p-1-perfect-squares})}.
\end{theorem}

\begin{proof}
Let us try to prove the infiniteness of $P_{\square}$ by means of the {\em Euler-Riemann zeta function} $\zeta(\alpha)$. Really one can consider the series $\sum_{1\le n\le\infty}\frac{1}{n^2}$, since each prime $p$ such that $p-1$ is a perfect square must be just of the type $p-1=n^2$. Now, it is well known that the series $\sum_{1\le n\le\infty}\frac{1}{n^2}$ is convergent. In fact, $\sum_{1\le n\le\infty}\frac{1}{n^2}=\zeta(2)=\frac{\pi^2}{6}$.\footnote{This has been first proved by Euler (1735) solving the so-called{\em Basel problem}.} Therefore, one cannot exclude that $P_{\square}$ is finite. In fact in such a case it surely should be satisfied the condition $\sum_{p\in P_{\square}}\frac{1}{p-1}\le\frac{\pi^2}{6}$. At the same time, the estimate  $\sum_{p\in P_{\square}}\frac{1}{p-1}\le \frac{\pi^2}{6}$ neither excludes the infinite alternative for the set $P_{\square}$, since in any case $\sum_{p\in P_{\square}}\frac{1}{p-1}$ must converge since it admits as majorant the convergent series $\sum_{1\le n\le\infty}\frac{1}{n^2}$. Therefore, the Euler-Riemann zeta function estimate of $P_{\square}$ (\ref{limitation-p-1-perfect-squares}) does not exclude that $P_{\square}$ is infinite.
\end{proof}

Let us conclude this section answering to the Conjecture \ref{the-infinite-set-ghost-rectangular}.

\begin{theorem}[The Conjecture \ref{the-infinite-set-ghost-rectangular} is false]\label{not-infinite-set-ghost-rectangular}
The set $P_{\square}$ is finite.
\end{theorem}
\begin{proof}
Let us first consider the following lemmas.
\begin{lemma}\label{necessity-even}
A $n$-ghost triangle,  $n>1$, must necessarily have $n$ even. (This condition is not sufficient !)
\end{lemma}
\begin{proof}
In fact, if $n=2m$, $n>1$, then $n^2$ is necessarily an even number $\le 4$, hence $n^2+1$ is an odd number $\ge 5$, hence can be a prime $p\ge 5$. Instead, if $n=2m+1$, $n>1$, then $n^2$ is again an odd number and $n^2+1$ is necessarily an even number, hence it cannot be a prime number $p>2$. Therefore in order that $n^2+1$ should be a prime number $p>2$ it is necessary that $n$ should be an even number $\ge 2$. To prove that the condition in the lemma is not sufficient it is enough to look to Tab. \ref{tab-examples-ghost-right-triangles} where there are many examples with $n$ even but with $n^2+1\not\in P$. For example $n=8$, where $n^2+1=65\not\in P$.
\end{proof}

\begin{lemma}\label{iff-condition}
A right triangle $ABC$ in $\mathbb{R}^2$ is a $n$-ghost triangle,  $n>1$, with $\overline{AC}=n$, and $\overline{BC}=1$, iff the Diophantine equation {\em(\ref{diophantine-equation-iff-condition})}
\begin{equation}\label{diophantine-equation-iff-condition}
  \varphi((2m)^2+1)=(2m)^2,
\end{equation}
is satisfied, where $\varphi:\mathbb{N}\to\mathbb{N}$ is the Euler's totient function.
\end{lemma}
\begin{proof}
\begin{equation}\label{exact-commutative-diagram-ghost-right-triangles}
   \xymatrix{&&\ar @{} [drr] |{\scalebox{2}{$\circlearrowright$}}\mathbb{N}\ar[r]^{\bullet}&\mathbb{N}\ar[r]^{+1}&\mathbb{N}\ar[d]^{\varphi}\\
   0\ar[r]&{}^{>1}P_{\square}\ar[r]&\ar @{} [drr] |{\scalebox{2}{$\circlearrowleft$}}\mathbb{N}\ar[u]^{2}\ar[rd]_{2}\ar@<2pt>[rr]^{b}\ar@<-2pt>[rr]_{a}&&\mathbb{N}\\
   &&&\mathbb{N}\ar[ur]_{\bullet}&\\}
\end{equation}

In fact a positive integer $m>1$ is prime iff $\varphi(m)=m-1$. This follows soon from the well-known {\em Euler product-formula}: $\varphi(m)=m(1-\frac{1}{a_1})\cdots (1-\frac{1}{a_k})$, when
$m=a_1^{r_1}\cdots a_k^{r_k}$, is the prime factorization of $m$.\footnote{Note that equation (\ref{diophantine-equation-iff-condition}) is not in contradiction with Lemma \ref{necessity-even} since $\varphi(r)$ is always even for $r\ge 3$.}
\end{proof}

From above lemmas we get that theorem is proved iff the set ${}^{>1}P_{\square}$ identified with (\ref{diophantine-equation-iff-condition-a})
\begin{equation}\label{diophantine-equation-iff-condition-a}
  {}^{>1}P_{\square}\equiv\{m\in\mathbb{N}\, |\, \varphi((2m)^2+1)=(2m)^2,\, m\ge 1\}
\end{equation}
is infinite. In fact $P_{\square}\, \leftrightarrow\, {}^{>1}P_{\square}$. Therefore, we can also write $$P_{\square}=\ker_{a}b=\{n\in\mathbb{N}\, |\,a(n)=b(n)\}\subset\mathbb{N} ,$$ where $a$ and $b$ are defined by the exact commutative diagram (\ref{exact-commutative-diagram-ghost-right-triangles}). We get that
$P_{\square}$ is infinite iff $\IM(a)\bigcap\IM(b)\subset P\subset\mathbb{N}$ is infinite. On the other hand one has $\sharp(\IM{(a)})=\aleph_0$, since there is the bijection $f_a:\mathbb{N}\to \IM{(a)}$, given by $f_a(n)=(2n)^2$. The structure of $\IM{(b)}$ is given by $$\IM(b)=\{|\mathbb{Z}^{\times}_{(2n)^2+1}|=|\mathbb{Z}^{\times}_{a_1^{r_1}}|\cdot|\mathbb{Z}^{\times}_{a_2^{r_2}}|\cdots |\mathbb{Z}^{\times}_{a_k^{r_k}}|\}_{n\in\mathbb{N}},$$ when $(2n)^2+1$ admits the prime factorization $(2n)^2+1=a_1^{r_1}\cdot a_2^{r_2}\cdots a_k^{r_k}$. Therefore we get also $\sharp(\IM(b))=\aleph_0$. From Tab. \ref{tab-examples-ghost-right-triangles} we know that $\IM(a)\bigcap\IM(b)\not=\varnothing$, but we understand also that $P_{\square}$ can be identified with a proper subset of $P$. Let us note that solutions set $\mathcal{S}ol(\mathbb{N})$ of the Diophantine equation (\ref{diophantine-equation-primes})
\begin{equation}\label{diophantine-equation-primes}
  \varphi(n)=n-1,\, n\in\mathbb{N}
\end{equation}

is the infinite subset $P\subset \mathbb{N}$ of primes, hence $\sharp(\mathcal{S}ol(\mathbb{N}))=\sharp(\mathbb{N})=\aleph_0$. This means that cardinality of the set of solutions of equation (\ref{diophantine-equation-primes}) must coincide with the cardinality of $\mathbb{N}$. This does not necessitate remain true for other infinite subsets of $P$. In fact, we have the following lemma.

\begin{lemma}[Cardinality and Diophantine equation]\label{cardinality-diophantine-equation}
Let fixed $X\subset \mathbb{N}$ be a subset $X$ of $\mathbb{N}$. Let us denote $\mathcal{S}ol(X)\subset \mathbb{N}$, the set of solutions of the Diophantine equation {\em(\ref{diophantine-equation-primes-a})}
\begin{equation}\label{diophantine-equation-primes-a}
  \varphi(m)=m-1,\, m\in X\subset \mathbb{N}.
\end{equation}
We shall prove that $\sharp(\mathcal{S}ol(\IM(b)))<\aleph_0$.
\end{lemma}
\begin{proof}
One has the commutative diagram (\ref{commutative-diagram-bijections}) of bijections and inclusions.\footnote{Warn ! We just know that $\mathcal{S}ol(\IM(b))$ is a proper subset of $P$. In fact must be $P_{\square}\subset P$. Therefore the commutativity of diagram (\ref{commutative-diagram-bijections}) does not mean that $P_{\square}$ is infinite.}
\begin{equation}\label{commutative-diagram-bijections}
   \xymatrix{&\mathbb{N}\ar[d]_{h}\ar[r]^{f}& \ar[l]\IM(b)&\\
P\ar[r]&\ar[l]\mathcal{S}ol(\mathbb{N})\ar[u]&\mathcal{S}ol(\IM(b))\ar@{_{(}->}[l]\ar@{_{(}->}[u]\ar[r]&P_{\square}\ar[l]\\}
\end{equation}
In order to prove this lemma let us try to use the same algebraic topologic approach followed in the previous sections. Then we say that two points $a,\, b\in \mathbb{R}^2$ are $n-\square$-bording if there exists a smooth curve $\gamma:[0,1]\to \mathbb{R}^2$, such that $\gamma(0)=a$, $\gamma(1)=b$ and such that it intersects the segment $[a,b]$ into a point $p$ that is a parabolic prime $p=n^2+1$, $n\in \mathbb{N}$, up to diffeomorphisms of $\mathbb{R}^2$. We denote by ${}^n\Omega_{\square}$ the corresponding bordism group. We put ${}^n\Omega_{\square}=\varnothing$ if $n^2+1\not\in P$. By using Lemma \ref{necessity-even} we get that ${}^{n=2m+1}\Omega_{\square}=\varnothing$, $m\ge 1$.

Let us prove that ${}^{n=1}\Omega_{\square}\cong\mathbb{Z}_2$. In fact any two points $a,\, b\in \mathbb{R}^2$ can be identified with the ends of the interval $[0,(n+1)^2=4]\subset\mathbb{R}$ of the $x$-axis. Now, we know from Theorem \ref{legendre-conjecture-proved} that in the interval $[n^2=1,(n+1)^2=4]$ there exists at least a prime. In particular we see that $n^2+1=2$ is a parabolic prime. The curve $\gamma:\, \{y=\sin(\hat x)\, \}_{\hat x=\frac{x}{\pi}}$, over the interval $[0,4]$, bords $0$ with $4$ and intersects the interval $[0,4]$ in the point $p=2$. Hence the points $0$ and $4$ on the $x$-axis $1-\square$-bord. Since any arbitrary couple of points $a,\, b\in \mathbb{R}^2$ can be identified with the points $0$ and $4$ of the $x$-axis, by means of suitable diffeomorphisms of $\mathbb{R}^2$, it follows that ${}^{n=1}\Omega_{\square}\cong\mathbb{Z}_2$. In general we can say that ${}^{n}\Omega_{\square}\cong\mathbb{Z}_2$ iff $n^2+1=p\in P$. In fact, iff $n^2+1=p\in P$ one has that $m_1\in[m_2]\in {}^{n}\Omega_{\square}$, for any $m_1,\, m_2\in \mathbb{N}$, such that $m_1<n^2+1=p<m_2$. Then since any arbitrary couple of points $a,\, b\in \mathbb{R}^2$ can be identified with the points $m_1$ and $m_2$ of the $x$-axis, by means of suitable diffeomorphisms of $\mathbb{R}^2$, it follows that ${}^{n}\Omega_{\square}\cong\mathbb{Z}_2$ iff $n^2+1=p\in P$. We can realize the isomorphism $f_*:{}^{n}\Omega_{\square}\cong {}^{1}\Omega_{\square}$, by a diffeomorphism $f:\mathbb{R}^2\to\mathbb{R}^2$, $(x,y)\mapsto
(\bar x=x\frac{2}{n^2+1},y)$. In fact by considering the curve $y=\sin(\hat x)$, $\hat x=\frac{x}{\pi}$, over the interval $[0,2.(n^2+1)]$, that intersects this interval in the parabolic prime $p=n^2+1\in P_{\square}$, after the above diffeomorphism the interval $[0,4]$ is intersect by the deformed curve in the point $2$ that is just a parabolic prime, ( $1^2+1=2\in P_{\square}$), corresponding to the bordism group ${}^{1}\Omega_{\square}$.

\begin{example}\label{examples-isomorphism-square-bordism-induced-diffeomorphisms}

$\bullet$\hskip 2pt The isomorphism $f_*:{}^{2}\Omega_{\square}\to {}^{1}\Omega_{\square}$, can be realized considering the diffeomorphism $(x,y)\mapsto (\bar x=x\frac{2}{5},y)$. In fact we can consider the curve $y=\sin(\hat x)$, $\hat x=\frac{x}{\pi}$, over the interval $[0,10]$, that intersects this interval in the parabolic prime $p=5\in P_{\square}$. After the above diffeomorphism the interval $[0,10]$ is deformed into $[0,4]$ that is intersect by the deformed curve in the point $2$ that is just a parabolic prime identifying the bordism group ${}^{1}\Omega_{\square}$.

$\bullet$\hskip 2pt Similarly the isomorphism $f_*:{}^{4}\Omega_{\square}\to {}^{1}\Omega_{\square}$, can be realized considering the diffeomorphism $(x,y)\mapsto (\bar x=x\frac{2}{17},y)$. In fact we can consider the curve $y=\sin(\hat x)$, $\hat x=\frac{x}{\pi}$, over the interval $[0,34]$, that intersects this interval in the parabolic prime $p=17\in P_{\square}$. After the above diffeomorphism the interval $[0,34]$ is deformed into $[0,4]$ that is intersect by the deformed curve in the parabolic prime $2$.

$\bullet$\hskip 2pt Again the isomorphism $f_*:{}^{10}\Omega_{\square}\to {}^{1}\Omega_{\square}$, can be realized considering the diffeomorphism $(x,y)\mapsto (\bar x=x\frac{2}{101},y)$. In fact we can consider the curve $y=\sin(\hat x)$, $\hat x=\frac{x}{\pi}$, over the interval $[0,202]$, that intersects this interval in the parabolic prime $p=101\in P_{\square}$. After the above diffeomorphism the interval $[0,202]$ is deformed into $[0,4]$ that is intersect by the deformed curve in the parabolic prime $2$.
\end{example}
\end{proof}

Let us assume, now, that the set $P_{\square}$ is finite: $P_{\square}=\{p_1=n^2_1-1<\, \cdots <\, <p_k=n^2_k-1\}$. Then a parabolic prime $p>p_k$ should not exist. Therefore for any integer $n>n_k$, we should have ${}^{n}\Omega_{\square}=\varnothing$. On the other hand for any diffeomorphism $f:\mathbb{R}^2\to\mathbb{R}^2$ we can induce an isomorphism $f_*:{}^{1}\Omega_{\square}\to f_*({}^{1}\Omega_{\square})$. In particular we can take $f_*$, induced by the inverse diffeomorphism $(x,y)\mapsto(x\frac{2}{p},y)$, where $p\in P$. Such a diffeomorphism deforms the interval $[0,2p]$ into $[0,4]$ and the curve $y=\sin(\hat x)$ on $[0,2p]$, with intersection $p\in[0,2p]$, into a smooth curve intersecting $[0,4]$ into the point $2$. Since we cannot state that $p$ is a parabolic prime, we cannot state that $f_*({}^{1}\Omega_{\square})$ coincides with some $n-\square$ bordism group ${}^{n}\Omega_{\square}$, with $p-1=n^2,\, p>p_k$. In other words, now the situation is different from the ones considered in the previous sections, since the diffeomorphisms $f$ are useful to prove the conjecture iff $p$ are parabolic primes ! (Instead in the previous sections it was enough to consider diffeomorphisms that reduce curves and intervals only, without any attention to other a priori requests.) Really, what can be obtained from above diffeomorphisms and induced isomorphisms, is the sequence of  isomorphisms reported in (\ref{parabolic-primes-and-bordism-isomorphisms}).
\begin{equation}\label{parabolic-primes-and-bordism-isomorphisms}
  {}^{1}\Omega_{\square}\cong{}^{2}\Omega_{\square} \cong{}^{6}\Omega_{\square}\cong{}^{10}\Omega_{\square}\cong\cdots\cong{}^{p}\Omega_{GB}\cong\mathbb{Z}_2,\, \forall\, p\in P,
\end{equation}

but we are not authorized to state that the set $P_{\square}$ is infinite, since if $p$ is not just a priori a parabolic prime the existence of the smooth curve $y=\sin(\hat x)$ on the interval $[0,2p]$ states only that $2p\in[0]\in{}^{\bar n}\Omega_{\square}$, with $\bar n^2+1=\bar p\in P_{\square}$ and $\bar p<p\in P$ !

Let us emphasize that the parabola $\Gamma:\, \{y=x^2+1\}\subset \mathbb{R}^2$, (see figure on the right in (\ref{pythagorean-and-ghost-right-triangles-in-n-2}), when considered on the lattice $\mathbb{N}^2$, identifies a discrete subset $\Gamma_{\mathbb{N}}\subset \Gamma$, encoded by the Diophantine equation $m=n^2+1$, ($m,\, n\in \mathbb{N}$):
$$\Gamma_{\mathbb{N}}=\{m=n^2+1\}_{n\in\mathbb{N}}\subset \Gamma.$$
Therefore we get also the following bijection and embedding $\Gamma_{\mathbb{N}}\, \leftrightarrow\, \mathbb{N}\, \hookrightarrow\, \mathbb{R}$. As a by product we get the following relations between the corresponding cardinalities. $$\sharp(\Gamma_{\mathbb{N}})=\sharp(\mathbb{N})=\aleph_0\, <\, \sharp(\mathbb{R})=\mathfrak{c}.$$ We can consider also $P$ embedded into $\mathbb{N}$ and $\sharp(P)=\aleph_0$, but this does not necessitate that $P\subset \Gamma_{\mathbb{N}}$. From Tab. \ref{tab-examples-ghost-right-triangles} we know that $P_{\square}\, \leftrightarrow\, P\bigcap\Gamma_{\mathbb{N}}\, \leftrightarrow\, \Gamma_{\square}\not=\varnothing$, but this does not necessitate that $\sharp(P\bigcap\Gamma_{\mathbb{N}})=\aleph_0$. With this respect,
in order to conclude the proof of Theorem \ref{the-infinite-set-ghost-rectangular} , let us consider some lemmas more.
\begin{lemma}\label{lemma-totient-parabolic-function}
The prime $p$, $p>2$, is parabolic iff $\varphi(p)=(2n)^2$, for some $n\ge 1$.\footnote{This lemma further justifies the construction of the figure on the right in (\ref{pythagorean-and-ghost-right-triangles-in-n-2}).}
\end{lemma}
\begin{proof}
In fact if we restrict the Euler's totient function $\varphi:\mathbb{N}\to 2\mathbb{N}$ to the subset of primes $P\subset\mathbb{N}$, we get that $\varphi(p)=p-1$. Therefore, taking also account of Lemma \ref{necessity-even}, $p>2$ is a parabolic prime iff $\varphi(p)=(2n)^2$, for some integer $n\ge 1$.
\end{proof}
\begin{lemma}[Equations of Ramanujan-Nagell type]\label{ramanujan-nagel-equations-type}
The set of the {\em Diophantine equations of Ramanujan-Nagell type} $x^2+A=BC^n$, with $A,\, B,\ C$ fixed and $x,\, n$ variable, admits a finite set of solutions.
\end{lemma}
\begin{proof}
See \cite{LJUNGGREN, RAMANUJAN} and the work by Carl Ludwig Siegel in \cite{SIEGEL, SARADHA-SRINIVASAN}, for a proof of the lemma about the Ramanujan-Nagell type equations.
\end{proof}

\begin{lemma}[Equations of Ramanujan-Lebesgue-Nagell type]\label{lebessgue-nagel-equations-type}
The set of the {\em Diophantine equations of Ramanujan-Lebesgue-Nagell type} $x^2+A=By^n$, with $A,\, B$ fixed and $x,\, y,\, n$ variable, admits a finite set of solutions.
\end{lemma}
\begin{proof}
See \cite{NAGELL1, NAGELL2, SARADHA-SRINIVASAN}.
\end{proof}
Let us now observe that the Diophantine equation in Lemma \ref{lemma-totient-parabolic-function} can be considered as a Ramanujan-Nagell type equation with the additional constraints $A=1=B$, $C=p\in P$ and $n=1$, or as the Ramanujan-Lebesgue-Nagell equation with the additional constraints $A=1$, $B=1$, $y=p\in P$, and $n=1$. Therefore, for example, if the set $\{x,y,n\}\in\mathbb{N}^3$ of solutions of the Ramanujan-Lebesgue-Nagell equation is finite, it follows that the set $P_{\square}$ of solutions of the Diophantine equation in Lemma \ref{lemma-totient-parabolic-function} must necessarily be finite too. As a by product we get $\sharp(\mathcal{S}ol(\IM(b)))<\sharp(\IM(b))=\aleph_0=\sharp(\mathbb{N})$.
In other words the Conjecture \ref{the-infinite-set-ghost-rectangular} is not true.
\end{proof}

\begin{table}[t]
\caption{Examples of correspondences $\IM(b)\, \leftrightarrow\, \mathbb{N}\, \leftrightarrow\,P\, \to\, P_{\square}$.}
\label{example-correspondences}
\scalebox{0.6}{$\begin{tabular}{|c|c|c|c|}
\hline
\hfil{\rm{\footnotesize $P$}}\hfil&\hfil{\rm{\footnotesize $\mathbb{N}$}}\hfil&\hfil{\rm{\footnotesize $\IM(b)\bigcup\{2\}$}}\hfil&\hfil{\rm{\footnotesize $P_{\square}$}}\hfill\\
\hline\hline
{\rm{\footnotesize $2$}}\hfil&\hfil{\rm{\footnotesize $1$}}\hfill&\hfil{\rm{\footnotesize $2$}}\hfil&\hfill{\rm{\footnotesize $2$}}\hfill\\
\hline
{\rm{\footnotesize $3$}}\hfil&\hfil{\rm{\footnotesize $2$}}\hfill&\hfil{\rm{\footnotesize $5$}}\hfil&\hfill{\rm{\footnotesize $5$}}\hfill\\
\hline
{\rm{\footnotesize $5$}}\hfil&\hfil{\rm{\footnotesize $3$}}\hfil&\hfil{\rm{\footnotesize $10$}}\hfil&\hfil{\rm{\footnotesize $17$}}\hfill\\
\hline
{\rm{\footnotesize $7$}}\hfil&\hfil{\rm{\footnotesize $4$}}\hfil&\hfil{\rm{\footnotesize $17$}}\hfil&\hfil{\rm{\footnotesize $37$}}\hfill\\
\hline
{\rm{\footnotesize $11$}}\hfil&\hfil{\rm{\footnotesize $5$}}\hfil&\hfil{\rm{\footnotesize $25$}}\hfil&\hfil{\rm{\footnotesize $101$}}\hfill\\
\hline
\hfil{\rm{\footnotesize $13$}}\hfil&\hfil{\rm{\footnotesize $6$}}\hfil&\hfil{\rm{\footnotesize $37$}}\hfil&\hfil{\rm{\footnotesize $197$}}\hfill\\
\hline
{\rm{\footnotesize $17$}}\hfil&\hfil{\rm{\footnotesize $7$}}\hfil&\hfil{\rm{\footnotesize $50$}}\hfil&\hfill{\rm{\footnotesize $257$}}\hfill\\
\hline
{\rm{\footnotesize $19$}}\hfil&\hfil{\rm{\footnotesize $8$}}\hfil&\hfil{\rm{\footnotesize $65$}}\hfill&\hfil{\rm{\footnotesize $401$}}\hfill\\
\hline
{\rm{\footnotesize $23$}}\hfil&\hfil{\rm{\footnotesize $9$}}\hfil&\hfil{\rm{\footnotesize $82$}}\hfil&\hfil{\rm{\footnotesize $577$}}\hfill\\
\hline
{\rm{\footnotesize $29$}}\hfil&\hfil{\rm{\footnotesize $10$}}\hfil&\hfil{\rm{\footnotesize $101$}}\hfil&\hfil{\rm{\footnotesize $677$}}\hfill\\
\hline
{\rm{\footnotesize $31$}}\hfil&\hfil{\rm{\footnotesize $11$}}\hfil&\hfil{\rm{\footnotesize $122$}}\hfil&\hfil{\rm{\footnotesize $1297$}}\hfill\\
\hline
\hfil{\rm{\footnotesize $37$}}\hfil&\hfil{\rm{\footnotesize $12$}}\hfil&\hfil{\rm{\footnotesize $145$}}&\hfil{\rm{\footnotesize $1601$}}\hfill\\
\hline
{\rm{\footnotesize $41$}}\hfil&\hfil{\rm{\footnotesize $13$}}\hfil&\hfil{\rm{\footnotesize $70$}}\hfil&\hfil{\rm{\footnotesize $2917$}}\hfill\\
\hline
{\rm{\footnotesize $43$}}\hfil&\hfil{\rm{\footnotesize $14$}}\hfil&\hfil{\rm{\footnotesize $197$}}\hfil&\hfil{\rm{\footnotesize $3137$}}\hfill\\
\hline
{\rm{\footnotesize $\cdots$}}\hfil&\hfil{\rm{\footnotesize $\cdots$}}\hfil&\hfil{\rm{\footnotesize $\cdots$}}\hfil&\hfil{\rm{\footnotesize $\cdots$}}\hfill\\
\hline
\end{tabular}$}
\end{table}

\begin{example}
In Tab. \ref{example-correspondences} are reported some valuations of the correspondence $$\IM(b)\, \leftrightarrow\, \mathbb{N}\, \leftrightarrow\,P\, \to\, P_{\square}.$$
There ``$\cdots$" closing the table mean that those correspondences must stop at some point.\footnote{It is a curiosity to know that the largest prime found is $p_{largest}=2^{57885161}-1$, (2013 January 2013). Since $p_{largest}-1=2(2^{57885160}-1)$, $p_{largest}$ cannot be a parabolic prime, i.e., $p_{largest}-1\not=n^2$, for some $n\in \mathbb{N}$. (For more information see e.g., \href{http://en.wikipedia.org/wiki/Largest_known_prime}{\tt Wikipedia/Largest Known Prime Number}.)}
\end{example}